\documentclass[12pt]{article}
\usepackage[a4paper,top=3cm,bottom=3.5cm,left=1.5cm,right=1.5cm]{geometry}

 \usepackage{hyperref}

\usepackage{amsmath}
\usepackage{amssymb}
\usepackage{amsthm}
\usepackage{mathrsfs}
\usepackage{calligra}

\usepackage{comment}
\usepackage{stmaryrd}

\numberwithin{equation}{section}

\theoremstyle{plain}
\newtheorem{proposition}{Proposition}[section]
\newtheorem{definition}[proposition]{Definition}
\newtheorem{lemma}[proposition]{Lemma}
\newtheorem{theorem}[proposition]{Theorem}
\newtheorem{corollary}[proposition]{Corollary}

\theoremstyle{definition}
\newtheorem{remark}[proposition]{Remark}
\newtheorem{example}[proposition]{Example}

\theoremstyle{remark}

\usepackage{graphicx}

\usepackage[bbgreekl]{mathbbol}
\usepackage{pifont, bbm}

\usepackage{color}

\def\Vect{\mathfrak{X}}
\def\VV{\Vect(A)}
\def\OOm{\Omega(A)}
\def\oA{\otimes_A}
\def\dd{\mathrm{d}}
\def\op{\mathrm{op}}
\def\id{\mathrm{id}}
\def\ev{\mathrm{ev}}
\def\coev{\mathrm{coev}}
\def\MMM{\mathscr{M}}

\def\FF{\mathcal{F}}
\def\ra{\triangleright}
\def\RA{\ra}

\def\LLL{\mathscr{L}^{\mathcal{}}}
\def\LL{\mathscr{L}^{\mathcal{}}}

\def\RR{\mathcal{R}}

\def\Tau{\mathcal T}

\def\le{\langle}
\def\re{\rangle}

\def\FF{\mathcal F}

\def\s'O{\stackrel{_{{\displaystyle\st \footnotesize '}}}{_{^{^{\displaystyle\otimes}}}}}

\def\1s{{1_\st }}
\def\3s{{3_\st }}
\def\2s{{2_\st }}
\def\ef1{{1_\FF}}
\def\ef2{{3_\FF}}
\def\ef3{{2_\FF}}

\def\TT{{\mathcal T}}

\def\le{\langle}
\def\re{\rangle}
\def\nn{\nonumber}

\def\ker{{\rm ker}}

\def\AA{A}

\def\dd{{\nabla}}
\def\dif{{\mathrm d}}

\def\Rsf{\mathsf{R}}
\def\Tsf{\mathsf{T}}
\def\Ls{\mathsf{L}}
\def\Rs{\mathsf{R}}

\def\fgf{\mathsf{f}}
\def\g{\mathsf{g}}
\def\gg{\mathfrak{g}}
\def\stt{\!\stackrel{_{_{\displaystyle\st }}}{_{^{_{\triangleleft}}}}\!}

\def\sym{{\rm{sym}}}

\def\ii{\mathrm{i}}
\def\of{{\bar{{\rm{f}}\,}}}

\def\oR{\bar R}
\def\al{\alpha}
\def\be{\beta}
\def\op{\mathrm{op}}

\def\Der{\mathrm{Der}}

\def\op{{{\scalebox{0.67}[0.67]{\mbox{$op$}}}\!}}
\def\cop{{{\scalebox{0.67}[0.67]{\mbox{$\! cop$}}}\!}}
\def\dd{\mathrm{d}}
\def\fp{\rm{fp}}
\def\geom{\rm{geom}}

\def\cat{{{}^H_\AA\MMM_\AA^\sym}}
\def\catm{{{}^H_\AA\MMM^{}_\AA}}
\def\catfp{{{}^H_\AA\MMM_\AA^{\sym, \rm{fp}}}}
\def\Cat{{{}^{H}_{\OM}\MMM_{\OM}^\sym}}
\def\Catfp{{{}^{H}_{\OM}\MMM_{\OM}^{\sym, \rm{fp}}}}
\def\CatM{{{}^{H}_{\OMM}\MMM_{\OMM}^\sym}}
\def\catFF{{{}^{H^\FF}_{\AA^\FF}\MMM_{\AA^\FF}^\sym}}
\def\catfpFF{{{}^{H^\FF}_{\AA^\FF}\MMM_{\AA^\FF}^{\sym, \rm{fp}}}}

\def\dotF{{\,\raisebox{1.5pt}{\scalebox{0.5}[0.5]{\mbox{$\bullet^{\mbox{$\FF$}}$}}}}}

\def\std{{\raisebox{\depth}{\scalebox{1.45}[-1]{\mbox{$\mathbb \Delta$}}}}_{\!}}
\def\dst{{}_{\,}\reflectbox{\mbox{$\std$}}_{\!}}

\def\pstd{{}^{\phantom{J}}_{\raisebox{\depth}{\scalebox{1.1}[-.75]{\mbox{${\mathbb
            \Delta}$}}}}}
\def\pdst{{}_{\,}{\reflectbox{\mbox{$\pstd$}}}}

\def\Dst{\dd_{\!\!}\pdst}
\def\stD{\dd\pstd}

\def\Dstst{\dd_{{}^{*\!}}\pstd}

\def\stDu{\stD{}_{\mbox{$\!{}^{}_u$}}}
\def\stDv{\stD{}_{\mbox{$\!{}^{}_v$}}}

\def\stDalu{\stD{}_{\mbox{$\!{}^{}_{\oR_\al\ra u}\!$}}}
\def\stDalpv{\stD{}_{\mbox{$\!{}^{}_{\oR^\al\ra v}\!$}}}

\def\ddst{{\phantom{\!\!\!\dd}}_{\!\!}\pdst}
\def\stdd{\phantom{\!\!\!\dd}\pstd}
\def\cd{{\rm{cd}}}

\def\OMPu{\Omega^{\bullet +1\!}}
\def\tOMPu{\tilde{\Omega}^{\bullet +1\!}}
\def\OMmu{\Omega^{\bullet -1\!}}
\def\OMP2{\Omega^{\bullet +2\!}}
\def\OMm2{\Omega^{\bullet -2\!}}
\def\OM{{\Omega^{\!\!\;\bullet}\!\!\:(A)}}
\def\tOM{{\tilde{\Omega}^{\!\!\;\bullet}\!\!\:(A)}}
\def\OMM{{\Omega^{\!\!\;\bullet}}}

\def\Con{\mathrm{Con}}
\def\hom{\mathrm{hom}}
\def\Hom{\mathrm{Hom}}
\def\hxi{h}

\def\stG{{}^*\Gamma}
\def\stU{{}^*\Upsilon}
\def\stGU{{}^{*\!}(\Gamma\otimes_A\Upsilon)}

\def\stL{{}^{*\!}L}

\def\sts{{}^*\!s}
\def\stt{{}^{*\!}t}
\def\stu{{}^{*\!}u}
\def\stsigma{{}^{*\!}\sigma}
\def\kg{K_\g}

\def\vv{{\rm{v}}}

\def\uXu{u}

\def\vvv{u}
\def\endex{\hfill \ding{122}}

  \DeclareMathAlphabet{\mathcalligra}{T1}{calligra}{m}{n}
 \DeclareFontShape{T1}{calligra}{m}{n}{<->s*[1.9]callig15}{}

\usepackage{authblk}

\title{\bf Cartan structure equations  and Levi-Civita connection in noncommutative  geometry}
\author{Paolo Aschieri}
\affil{\small{
\centerline{{\sl Dipartimento di Scienze e Innovazione 
Tecnologica, }} {{Universit\`a del Piemonte Orientale}}
\centerline{{\sl  
Viale T. Michel - 15121, Alessandria, Italy}}
\centerline{{\sl INFN, Sezione di Torino, 
  via P. Giuria 1, 10125 Torino, Italy}}
\centerline{\texttt{paolo.aschieri@uniupo.it}}}}
\date{}

\begin{document}

\maketitle
 
\abstract{
We study the differential and Riemannian geometry of algebras $A$ endowed with an action of a triangular
  Hopf algebra $H$ and noncommutativity compatible with the associated
  braiding. 
The modules of one forms and of braided derivations are modules in a
compact closed category of $H$-equivariant $A$-bimodules, whose
internal morphisms correspond to tensor fields. 
Vector fields and forms approaches to curvature and torsion are proven to be
equivalent by extending the Cartan calculus to left (right) $A$-module
(not necessarily $A$-bimodule) connections. The Cartan structure
equations and the Bianchi identities are derived.

Existence and uniqueness of the Levi-Civita connection
  for {\sl arbitrary} pseudo-Riemannian metrics is proven via a Koszul formula. The general
  theory includes Drinfeld twists of commutative geometries and also
  cotriangular Hopf algebras. It is illustrated with the example of the tensor square of Sweedler Hopf algebra which becomes a noncommutative Einstein manifold via a non-central metric}.

\tableofcontents

\section{Introduction}
Noncommutative Riemannian geometry is an active and interdisciplinary
research field.
One line of study follows Connes’ approach, focusing on spectral geometry and the analysis of the Laplacian (see [25] for a recent review). Another line adopts a more algebraic perspective, beginning with a differential calculus on a noncommutative algebra, introducing a notion of Riemannian metric, and addressing the existence and uniqueness of the Levi-Civita connection. Interest in the field is also driven by the possibility that gravity on noncommutative spacetimes may capture aspects of a quantum theory of gravity.

Given a noncommutative algebra $A$ and an associated differential
 calculus $(\dd,\Omega(A))$, there are different approaches to the notion of metric tensor and to that of Levi-Civita connection. On one hand one can consider metrics $\g$ that are compatible with the noncommutative structure, i.e. that are central elements of the bimodule $\Omega(A)\otimes_A\Omega(A)$ ($a\g=\g a$ for all $a\in
A$), or of the dual bimodule. On the other hand one can relax this condition and study arbitrary metrics, which is useful when considering $\g$ as a dynamical field. Similar considerations holds for the study of connections. 
In particular the notion of Levi-Civita connection relies on the possibility of imposing the metric compatibility
condition $\nabla \g=0$. This requires lifting the connection from the module $\Omega(A)$ of one forms to the tensor product module
$\Omega(A)\otimes_A\Omega(A)$. However, lifting connections on tensor products modules is a nontrivial problem. In
the literature it is typically overcome  by constraining the connection to
be compatible with the noncommutative structure in the sense of being a bimodule connection \cite{D-VM96}. For approaches to Riemannian geometry using 
central metrics and bimodule Levi-Civita connections see \cite{Beggs-Majid}, \cite{Beggs-Majid1},
\cite{BGL1}, \cite{BGL2}, \cite{Weber}, \cite{Matassa}, see also the monograph \cite{BM}.

For a selected class of noncommutative algebras on the other hand it is possible to
relax this centrality contraint on  the metric and consider $\g$ just as a (properly defined symmetric)
element in $\Omega(A)\otimes_A\Omega(A)$. As said, this is relevant when
$\g$ is seen as a dynamical field, so that for example $a\g$ is also
a metric for any invertible $a\in A$. 
In case of $\mathbb{R}^n$ with Moyal-Weyl noncommutativity the Levi-Civita
connection of an arbitrary symmetric metric was constructed in
\cite{NCG1} using a noncommutative Koszul formula (see also \cite[\S 3.4, \S 8.5]{book}). A similar result holds  for
the noncommutative torus  \cite{Rosenberg}.
These results and those in \cite{Arnlind} rely on the existence
of (undeformed) derivations of the noncommutative algebra $A$
generating the $A$-module of vector fields (dual to that of
one forms).
A generalization via local charts with Moyal-Weyl noncommutativity is
in \cite{AC} in the deformation quantization context.

In these  examples the noncommutative algebra $A$ is endowed with a representation of a
triangular Hopf algebra $H$ and the noncommutativity is compatible with the braiding $\tau$ given
by the triangular structure. Triangularity of $H$ implies
that $\tau$ is a representation of the permutation group. We shall
call these algebras braided commutative $H$-module algebras,
they are commutative algebras $A$ in the category of $H$-module
algebras with $H$ triangular.
Besides noncommutative tori and  Connes-Landi  spheres \cite{CLandi},
belong to this class the algebras studied in \cite{CDV}, 
the noncommutative algebras obtained from Drinfeld twists of
commutative ones and cotriangular Hopf algebras (see \cite{EG} for
their classification in the finite dimensional case).

In this paper we develop a general theory of Riemannian geometry for
such algebras $A$.
 We study arbitrary pseudo-Riemannian metrics $\g$ (not
only central) and prove existence and uniqueness of 
the Levi-Civita connection (which need not be a bimodule connection). This is achieved by
providing a Koszul formula, see \eqref{Koszul}, that leads to an explicit expression of the Levi-Civita connection, see  \eqref{LEVICIVITA}.

Noncommutative differential geometry is studied primarily in the language of
forms, the complementary vector fields approach being more problematic
since the linear space of derivations of a noncommutative algebra $A$ is not generally
an $A$-bimodule.
For $A$ a braided commutative $H$-module algebra,
vector fields are braided derivations and a braided derivation based
differential calculus and geometry can be constructed. The general theory
of algebras, Lie algebras and differential operators in symmetric
monoidal categories was outlined in the late '80s in \cite[\S 13.5]{Manin}. 
The differential and Cartan calculus was pioneered in \cite{Gurevich},
for a braided derivations approach see \cite{Weber}.  When considering connections one can  contract them with
vector fields to obtain covariant derivatives along vector
fields. These were used in \cite{NCG2} in the context of twist deformation of commutative
algebras, and led  to define curvature and torsion as
left $A$-linear maps on tensor products of vector fields, a quantum
analogue of the standard general relativity definitions. These constructions did not
require lifting covariant derivatives from vector fields to tensor
fields. 
In \cite{AS} such a lifting, within the setting of braided commutative $H$-module algebras, was provided for connections that are not invariant under the $H$-action, these are, in general, not bimodule connections.
Concerning covariant derivatives, the lifting is understood when the vector field is a
derivation, or when the connection is a bimodule connection, or when it is
$H$-equivariant so that the braiding acts trivially.

We study here the covariant derivative and its
lifting to tensor products in full generality. This framework allows us to:  {\sl i)} Extend the Cartan calculus to include covariant derivatives; {\sl ii)} Establish the equivalence between the formulations of curvature and torsion in terms of forms and in terms of vector fields;
 {\sl iii)} Derive the Cartan structure equations and the Bianchi identities; and {\sl iv)} Obtain a
general Koszul formula determining the Levi-Civita connection via the metric compatibility condition $\nabla_{\!u\,}\g=0$ for any vector field (braided derivation) $u$.

Our analysis is based on considering a new perspective on operators on
$H$-equivariant $A$-bimodules: such operators can be expressed as compositions of operators acting from the
left and from the right. This is not merely a matter of
notation but reflects covariance requirements with respect to the
(generally) non-cocommutative Hopf algebra $H$.
In particular, the covariant derivative $\std{}_u$ is a composition of a left
connection $\std$ (acting from the right) and an
inner derivative (acting from the left). Only when
the vector field $u$ is a derivation of the algebra can
$\std{}_u$ be seen as an operator acting form the right. Similarly, its extension $\stDu$ from modules over  $A$ to modules over the exterior algebra $\Omega^\bullet(A)$
is with operators $u$ and $\rm{d}$ acting from the left, while  $\std$ acts from the right. It
is precisely with this composition that we extend the Cartan calculus to include covariant derivatives,
proving the braided Cartan calculus relation
$[\stDu,\ii_v]=\ii_{[u,v]}$. This establishes result {\sl i)} and is the key identity underpinning {\sl ii)}, the equivalence between the vector field and the exterior form
formulations of curvature and torsion.
 Thus, for example, the curvature
tensor defined in \cite{NCG2} as the braided commutator of covariant derivatives
along vector fields
 corresponds to the standard definition given by squaring the
connection $\std$ (cf. Theorem \ref{RD2} and Theorem \ref{Tasd} for
the analogous result on torsion). A similar correspondence for the torsion tensor
was established in \cite{BGL2} within the specific framework of tame differential calculi (where the
$A$-bimodule of 1-forms is generated by central 1-forms, and vector fields
are derivations).
The relation between (left)
connections on vector fields and the dual (right) connections on one
forms then leads to result {\sl iii)}: the Cartan structure equations for curvature and torsion and the associated Bianchi identities. 
We also lift the covariant derivative along vector fields to act on
tensor fields using the results of \cite{AS}, where as mentioned, connections were lifted to tensor product modules without assuming the bimodule connection property.

Finally, upon introducing a pseudo-Riemannian metric, this noncommutative
differential geometry is used to provide {\sl iv)}, a Koszul formula for  metric compatible torsion free connections, yelding both the existence and uniqueness of the Levi-Civita connection, as well as an explicit and global constructive method.
This result is very general, as it assumes
neither $H$-equivariant metrics (as in \cite{Weber}) nor a preferred
set of derivations of the noncommutative algebra. When $H$-equivariant metrics are considered we obtain bimodule Levi-Civita connections and  recover
the results of \cite{Weber}.
It also significantly extends the result of \cite{AC} in the context of formal Drinfeld twist deformation quantization, where the Koszul formula was provided either locally, subordinated to open charts with Moyal–Weyl noncommutativity, or globally, in implicit form  as a formal power series
in $\hbar$ starting with the commutative Koszul formula and determined recursively in $\hbar$.

In the present work, pseudo-Riemannian metrics are just
braided symmetric non-degenerate contravariant tensors. This is the natural context where to formulate a
noncommutative gravity theory where the metric is the dynamical
field. We here present in vacuum Einstein equations
leading to noncommutative Einstein manifolds.
This general theory is illustrated with the example of the
noncommutative torus (and the Moyal-Weyl plane), which is sketched in order to make
contact with previous literature, and the example of the tensor square
of Sweedler Hopf algebra. This is a cotriangular finite dimensional
Hopf algebra with abelian two dimensional Lie algebra of braided
derivations. The Levi-Civita connection, curvature,  Ricci tensor and
scalar curvature are explicitly presented for any metric. We solve the noncommutative Einstein manifold condition
for a non-central non-equivariant metric.
\\

The algebraic structure underlying this study is that of the 
categories of $H$-modules and of $H$-equivariant
$A$-bimodules (compatible $H$-modules $A$-bimodules). Following \cite{AS} and the
sharpened results in \cite{BSS1}, \cite{BSS2}  we recall 
the different structures of modules and module maps we need in
noncommutative Riemannian geometry. This clarifies the constructions
and the different general properties needed in the progress of the paper.
For example, left (right) connections are linear maps but are not morphism in the
category ${}^H\MMM$ of $H$-modules, in categorical terms they are internal morphisms. Left connections have different $H$-action from right connections; they are different internal morphisms in
${}^H\MMM$.  The covariant derivative associated with a left connection, in turn,  can be viewed both as a left $A$-linear map and an internal morphism in the
category $\catm$ of  $H$-equivariant $A$-bimodules (cf. Remark \ref{remcd}). Torsion and
curvature exhibits further structural properties: they are internal morphisms (left $A$-linear maps
transforming under the $H$-adjoint action) in the
subcategory of symmetric $H$-equivariant $A$-bimodules
$\cat$. In particular, we study internal morphisms associated with tensor
products of  finitely generated and projective modules and their duals.
This allows, as in classical differential geometry,  to understand
noncommutative tensors fields in all their different forms, as
elements of $A$-bimodules (sections), or as various left (right)
$A$-module maps. The underlying category in this richest case is the
compact closed category of symmetric $H$-equivariant $A$-bimodules 
finitely generated and projective as $A$-modules (it is a ribbon category
with trivial twist isomorphisms since $H$ is triangular) .
 
This underlying categorical context is presented in Section 2 with
some details, in particular on the biclosed and rigid structures, so
that the paper is self-contained. The section closes with the examples
where $A$ is a cotriangular Hopf algebra and where it is the Drinfeld
twist deformation of the algebra of smooth functions on a manifold $M$. 
In Section 3 the braided derivations formulation of the differential and
Cartan calculus is revisited. 
We provide the example of the calculus on cotriangular Hopf algebras
and show its equivalence with the (categorical) exterior algebra
approach of \cite{Gomez-Majid}.   
In Section 4 we review known properties of right connections
and left connections, since both are relevant for understanding
curvature and torsion. We then introduce the covariant derivative, 
establish the Cartan formula $[\stDu,\ii_v]=\ii_{[u,v]}$ and show that the
curvature tensor, defined as the square of the connection, can be
equivalently defined via the commutator of covariant derivatives along
vector fields. Similarly, two different definitions of torsion are
shown to be equivalent. 
In Section 5 the relation between connections 
on modules and dual modules is recalled and that between curvatures
studied. This leads to relate the vector field formulation of
curvature and torsion to the one involving the dual
connection on forms.  These are the Cartan
structure equations for curvature and torsion since choosing a basis
(a trivialization of the frame bundle) we
obtain those structure equations. The associated
Bianchi identities are also obtained.
In Section \ref{riemgeom} existence and uniqueness of the Levi-Civita
connection is proven for any braided symmetric metric.
The Ricci tensor and the scalar curvature are canonically defined,
this leads to in vacuum
Einstein field equations giving noncommutative Einstein manifolds. 
In Section \ref{expex} the Riemannian geometry of the tensor square
of Sweedler Hopf algebra is presented.

\section{Hopf algebras, braidings and
  representations}\label{section2}
 We work in the category of $\Bbbk$-modules, with $\Bbbk$ a fixed
  field of characteristic zero 
or the ring of formal power series in a variable $\hbar$ over such field;
much of what follows holds for a commutative unital ring. The tensor
product over $\Bbbk$ is denoted $\otimes$. Algebras over $\Bbbk$ are
assumed associative and unital. Hopf algebras are assumed with invertible antipode. 

In Section \ref{SCMC} we study right (left) $\Bbbk$-linear maps
and $A$-linear maps between $H$-modules and between $H$-equivariant
$A$-bimodules, hence introducing biclosed (left and right closed)
monoidal categories and, when the Hopf algebra $H$ is triangular,
symmetric biclosed monoidal categories.
In Section \ref{fgp} we continue the study of right (left) $A$-linear
maps considering, for an arbitrary Hopf algebra $H$, the case of
a finitely generated and projective $A$-module, this is the same as a rigid module. 
In Section \ref{sect:ribbon} we study these rigid $H$-equivariant
$A$-bimodules when $H$ is triangular, the corresponding category is an
example of a compact closed category and hence of a ribbon category.
Most of these results are covered (albeit sometimes implicitly) either in the literature on quantum groups
(see e.g. \cite{Kassel,Ma95}) or in the related one on tensor
categories (see e.g. \cite{EGNO}).
We here set the notation (following \cite{AS,
  BSS1, BSS2}) and present the main results that will be used for the
later sections, spelling out  in particular the properties of 
$\Bbbk$-linear and of right (left) $A$-linear maps transforming under different
left $H$-adjoint actions. 
A last section is devoted to the example of the category of
bicovariant bimodules of a cotriangular Hopf algebra and to
that of noncommutative vector bundles obtained via Drinfeld twist
deformation. 

\subsection{Closed monoidal categories}\label{SCMC}
We start recalling basic Hopf algebra notions, the  category of
$H$-modules for an arbitrary Hopf algebra $H$ and also for 
a triangular Hopf algebra $H$. In this simple context we
introduce  $\Bbbk$-linear
maps that are not invariant under the $H$-action, they come with two
different $H$-actions structuring them as $\Bbbk$-linear maps acting
from the right or from the left.
Their categorical interpretation, as internal morphisms, is also discussed. They structure the monoidal category of $H$-modules
${}^H\MMM$ as a biclosed monoidal category.
If $H$ is a triangular Hopf algebra then $({}^H\MMM, \otimes,
{}_\Bbbk\hom,\hom_\Bbbk)$ is furthermore a symmetric biclosed monoidal category and there is a 
tensor product $\otimes_\RR$ of internal morphisms.
 
Given an $H$-module algebra $A$ we then study $H$-equivariant
$A$-bimodules 
and the associated right $A$-linear maps and
left $A$-linear maps (internal morphisms). This is the biclosed
monoidal category $(\catm, \otimes_A,{}_A\hom, \hom_A)$. If $H$ is triangular and the product in $A$
is compatible with the braiding, restricting to modules where the
$A$-bimodule structure is compatible with the braiding (symmetric $A$-bimodules)
we have the biclosed  monoidal subcategory  $(\cat  \otimes_A,{}_A\hom, \hom_A)$, which is symmetric. 

This chain of results holds as well when we consider a graded algebra
instead of $A$, and graded modules.

\subsubsection{Modules over a Hopf Algebra}\label{subsecMoHA}
Let $H$ be a Hopf algebra $(H,\mu,\eta,\Delta,\varepsilon, S)$ over $\Bbbk$.
We denote by $^{H}\MMM^{}$ the category of left $H$-modules, where objects
 in ${}^{H}\MMM^{}$ are $\Bbbk$-modules $V$ with a left
$H$-action $\ra : H\otimes {V}\to {V}$, while morphisms in ${}^H\MMM$
are $\Bbbk$-module maps $f: V\to W$ that are $H$-equivariant,  i.e.,
\begin{equation}\label{eqn:Hequivariance}
 \hxi\ra f(v) =f(\hxi\ra v) ~,
\end{equation}
for all $\hxi\in H$ and $v\in {V}$; we write $f\in \Hom^{}_{{}^{\,H}\!\MMM}(V,W)$.
In this paper $H$-modules will be always left $H$-modules and will be
simply called $H$-modules.

Since $H$ is a  bialgebra ${}^H\MMM$ is a (strict)  monoidal
category. 
Given two $H$-modules $V$ and $W$
their tensor product $V\otimes W$ 
is an $H$-module  with $H$-action
\begin{equation}
\ra: H \otimes {V} \otimes {W} \longrightarrow  {V} \otimes {W} ~, 
\qquad h \otimes v \otimes w \longmapsto h\ra (v\otimes w):=(h_{(1)} \ra v) \otimes (h_{(2)} \ra w)~,
\end{equation}
where we have used the Sweedler notation $\Delta(h) = h_{(1)}\otimes h_{(2)}$ (with summation understood)
for the coproduct of $H$. 
The tensor product of two morphisms  in ${}^{H}\MMM^{}$, $f: V \to V^\prime,~ g: W\to
W^\prime$ is the morphism in ${}^{H}\MMM^{}$ defined by
$f\otimes g:
V\otimes W\to V'\otimes W', v\otimes w\mapsto f(v)\otimes
g(w)$. The tensor product functor $\otimes$ is associative.
The unit object in ${}^H\MMM$ is $\Bbbk$ with left $H$-action given by the counit of $H$,
$\ra : H\otimes \Bbbk\to \Bbbk\,,~h\otimes \lambda \mapsto \epsilon(h)\,\lambda$.

Since the bialgebra $H$ is  a Hopf algebra,  ${}^H\MMM$ is a closed monoidal category.
For any $V,W$ in ${}^H\MMM$,  we denote by 
$\hom_\Bbbk(V,W)$ in ${}^H\MMM$ the $\Bbbk$-module
${\rm Hom}_{\Bbbk}(V,W)$ of $\Bbbk$-linear maps  $L : V\to W$ 
 equipped with the  adjoint $H$-action
\begin{equation}\label{adjact}
\RA_{\,} : H \otimes^{} {\hom_\Bbbk}^{}\big({V}, {W}\big) \longrightarrow 
 {\hom_\Bbbk}^{}\big({V}, {W}\big)~, ~~
h \otimes L \longmapsto h\RA L:=h_{(1)}\!\!\: \ra\,\circ \,L \,\circ\, S(h_{(2)})\ra ~\!, 
\end{equation}
i.e., $(h\RA L)(v)= h_{(1)} \ra (L(S(h_{(2)})\ra v))$. Given morphisms $f^\op: V\to V', g: W\to W'$ (where $f^\op: V\to V'$ is
just $f: V'\to V$ thought as a morphisms in the opposite category
$({}^H\MMM)^\op$) we have the morphism
\begin{equation}\label{eqn:internalhommorph}
{\hom_\Bbbk}^{}(f^\op,g): {\hom_\Bbbk}^{}\big({V},{W}\big) \longrightarrow 
{\hom_\Bbbk}^{}\big({V}^\prime,{W}^\prime\, \big)~, \qquad L \longmapsto g \circ L \circ f~.
\end{equation}
This way we have defined the so-called internal-hom functor
$
{\hom_\Bbbk}: \big({}^H_{}\MMM^{}\big)^\op \times {}^H_{}\MMM^{} \longrightarrow {}^H_{}\MMM^{}~.
$
This is compatible with the tensor product functor $\otimes$, indeed since any
$H$-equivariant map $f:V\otimes W\to Z$ can be considered as an
$H$-equivariant map $\zeta(f): V\to \hom(W,Z)$ via $\zeta(f)(v):=f(v,
\mbox{-})$, we have that the functor $\mbox{-}\otimes W$ is left adjoint to $\hom_\Bbbk(W,\mbox{-})$,
thus $({}^H\MMM, \otimes, \hom_\Bbbk)$ is a
closed  monoidal category. 
\\

Let $V,W$ be modules in ${}^H\MMM$, we can define another 
$H$-adjoint action
on $\Bbbk$-linear maps
$V\to W$. We denote 
by ${}_\Bbbk\hom(V,W)$ the $\Bbbk$-module ${\rm Hom}_{\Bbbk}(V,W)$
with $H$-action $\RA^\cop$ defined by
\begin{equation}\label{copadjact}
\RA_{\,}^\cop : H \otimes^{} {{\,}_\Bbbk\hom}^{}\big({V}, {W}\big) \longrightarrow 
 {{}_\Bbbk\hom}^{}\big({V}, {W}\big)~\!, ~~
\hxi_{}\otimes \tilde L \longmapsto \hxi\RA^\cop \tilde L:=\hxi_{(2)}
\ra \, \circ \,\tilde L \,\circ\, S^{-1}(\hxi_{(1)}) \ra ~\!,
\end{equation}
$$
~~~~~(h\ra^\cop\tilde L)(v)=\hxi_{(2)}
\ra (\tilde L(S^{-1}(\hxi_{(1)}) \ra v))~\!.
$$
This gives the monoidal structure  $({}^H\MMM, \otimes,
{}_\Bbbk\hom)$, with the functor  $V\otimes \mbox{-}$ that is  left
adjoint to ${}_\Bbbk\hom(V,\mbox{-})$ (via $\zeta(f): W\to {}_\Bbbk\hom(V,Z)$, 
$\zeta(f)(w):=f(\mbox{-},w)$,  for any  $f\in
\Hom^{}_{{}^{\,H}\!\MMM}(V\otimes W,Z)$).

 While $\Bbbk$-linear maps 
$L\in \hom_\Bbbk(V,W)$ naturally act from the left, indeed the $\RA$ adjoint
action satisfies, for all $\hxi \in H$, $v\in V$, $\hxi\ra
(L(v))=(\hxi_{(1)}\RA L)(\hxi_{(2)}\ra v)$, $\Bbbk$-linear maps 
$\tilde L\in {}_\Bbbk\hom(V,W)$ naturally act from the right, indeed the $\RA^\cop$ adjoint
action satisfies, 
\begin{equation}\label{copadjact2}
\hxi\ra
(\tilde L(v))=(\hxi_{(2)}\RA^\cop \tilde L)(\hxi_{(1)}\ra v)~,
\end{equation}
 that, evaluating
$\tilde L$ on $v$ from  the right, reads $\hxi\ra((v)(\tilde L))=(\hxi_{(1)}\ra
v) (\hxi_{(2)}\RA^\cop \tilde L)$.
\\

In summary, associated with the Hopf algebra $H$, we have the
biclosed monoidal category $$({}^H\MMM, \otimes, \hom_\Bbbk,
{}_\Bbbk\hom)~.$$
The submodules $\hom_\Bbbk\!\!\:(V,W)^H\subset \hom_\Bbbk(V,W)$ and ${}_\Bbbk\hom(V,W)^H\subset {}_\Bbbk\hom(V,W)$ of
$H$-invariant elements, i.e., $h\ra L=\varepsilon(h)L$ and $h\ra^\cop
\tilde L=\varepsilon(h)\tilde L$, coincide with
that of $H$-equivariant maps $V\to W$ (compare for example
\eqref{copadjact2} with \eqref{eqn:Hequivariance}) and are hence identified with $\Hom_{{\!\:}^{H\!\!\:}\MMM}(V,W)$. 

\subsubsection{Modules over a triangular Hopf Algebra}
Let now $H$ be a triangular Hopf
algebra with universal 
$\mathcal{R}$-matrix $\mathcal{R}\in H\otimes H$ .
We recall that it satisfies
$$
\Delta^{\cop}(\hxi)
=\mathcal{R}\Delta(\hxi)\mathcal{R}^{-1}
\text{ for all }\hxi\in H,
$$
$$
    (\Delta\otimes\mathrm{id})\mathcal{R}
    =\mathcal{R}_{13}\mathcal{R}_{23}
   ~,~~
    (\mathrm{id}\otimes\Delta)\mathcal{R}
    =\mathcal{R}_{13}\mathcal{R}_{12}
$$
and the triangularity condition $\mathcal{R}_{21}=\mathcal{R}^{-1}$.
Because of the triangular structure the  monoidal category ${}^H\MMM$
is symmetric: 
the braiding $\tau^{}:=\tau(\RR)$ is the natural isomorphism $\tau: \otimes^{} \Rightarrow \otimes^\op$
with components defined by, 
\begin{equation}\label{taubraiding}
\tau^{}_{V,W}: V \otimes^{} W \longrightarrow W \otimes^{} V ~, \qquad v \otimes^{} w 
\longmapsto 
\big(\bar R^{\alpha} \ra w\big) \otimes^{} \big(\bar R_{\alpha} \ra
v\big)
\end{equation}
where we used the notation $\mathcal{R}=R^\alpha\otimes R_\alpha$,
$\mathcal{R}^{-1}=\bar R^\alpha\otimes \bar R_\alpha$.
The  category ${}^H\MMM$ is symmetric because for
all $V, W$, $\tau^{}_{W,V}\circ \tau^{}_{V,W}=\id_{V\otimes W}$, hence
$\tau$ provides a representation of the permutation group.
With slight abuse of notation we  shall frequently omit the indices in the
isomorphisms $\tau^{}_{V,W}$ and  simply write $\tau$.

The functors $\,\mbox{-}\otimes W$ and $W\otimes \mbox{-}$ in this
triangular case are naturally isomorphic via the
braiding; correspondingly, the two internal-hom functors $\hom_\Bbbk$
and  ${}_\Bbbk\hom$ are naturally isomorphic via the family of isomorphisms,
for all $V,W\in {}^H\MMM$,

\begin{equation}\label{LefttoRight} 
{\cal D}_{V,W}:\hom_\Bbbk(V,W) \longrightarrow {}_\Bbbk\hom(V,W)~,~~
L\mapsto {\cal D}^{}_{V,W}(L):=
(\oR^\al\ra L) \circ \oR_\al\ra
\end{equation} 
with inverse given by (cf. \cite[\S 3.2, Rmk. 3.11]{AS} and recall that a
triangular structure $\mathcal{R}\in H\otimes H$ is in particular a
twist or 2-cocycle, cf. \cite[Ex. 2.3.6]{Ma95})
\begin{equation}\label{RighttoLeft} 
~{\cal D}^{-1}_{V,W}: {}_\Bbbk\hom(V,W)\longrightarrow  \hom_\Bbbk(V,W)~, ~~\tilde L\mapsto
{\cal D}^{-1}_{V,W}(\tilde L)=(R^\al\RA^\cop \tilde L)\circ R_{\al\,}\ra~.
\end{equation}   
Here we just prove $H$-equivariance: for all $\hxi\in H$,
$h\ra {\cal D}^{-1}_{V,W}(\tilde L)=~{\cal
  D}^{-1}_{V,W}(h \ra^\cop\tilde L)$. This equality
is equivalent  to
$h\ra {\cal D}^{-1}_{V,W}(\tilde L)\circ h_{(2)}\ra=~{\cal
  D}^{-1}_{V,W}(h_{(1)} \ra^\cop\tilde L)\circ h_{(2)}\ra$
(use $\Bbbk$-linearity and
$\hxi_{(1)}\otimes \hxi_{(2)}S(\hxi_{(3)})=\hxi$). Recalling
\eqref{adjact} the left hand side equals 
$\hxi\ra\;\!\circ\;\!{\cal D}^{-1}_{V,W}(\tilde L)$.
Recalling \eqref{copadjact}, that
$\RA^\cop$ and $\RA$ are actions and using quasi-cocommutativity in the form $\mathcal{R}\Delta(\hxi)
=\Delta^{\cop\!\:}(\hxi)\mathcal{R}$, the right hand side too equals this expression:
${\cal  D}^{-1}_{V,W}(h_{(1)} \ra^\cop\tilde L)\circ h_{(2)}\ra=
(R^\al h_{(1)}\RA^\cop\tilde L)\circ R_\al h_{(2)}\ra=
(h_{(2)}R^\al\RA^\cop\tilde L)\circ h_{(1)}R_\al \ra=
h_{(2)}\!\!\:\RA^\cop\!\!\: (R^\al\RA^\cop\tilde L)
\!\:\circ h_{(1)}R_\al\!\: \ra=
\hxi\ra\;\!\circ\;\!{\cal D}^{-1}_{V,W}(\tilde L)$.\\

In summary, when the Hopf algebra $H$ is triangular the quadruple
 $$({}^H\MMM, \otimes, \hom_\Bbbk, {}_\Bbbk\hom)$$ 
is a symmetric biclosed  monoidal category.

In a braided closed monoidal category we can evaluate, compose and
consider tensor products not just of morphisms but also of internal
morphisms. For quasitriangular Hopf algebras and hence for triangular
Hopf algebras, internal morphisms evaluation and composition are the
usual ones of $\Bbbk$-linear maps on  $\Bbbk$-modules.
The composition of internal morphisms is  easily seen to be an
internal morphism; we give a proof for internal morphisms carrying the $\ra^\cop$
adjoint action.  Let $\tilde L\in {}_\Bbbk\hom(W,Z)$ and
$\tilde L'\in {}_\Bbbk\hom(V,W)$, for all $v\in V$, iterating expression
\eqref{copadjact2} we have,
\begin{equation}\label{eqcompadjcop}
\begin{split}h\ra \big(\tilde L\circ \tilde L'\:(v)\big)=h\ra \big(\tilde L(\tilde
L'(v))\big)&=
(h_{(2)}\ra^\cop \tilde L)\big(h_{(1)}\ra(\tilde L'(v))\big)\\[.2em]&=
(h_{(3)}\ra^\cop \tilde L)\circ (h_{(2)}\ra^\cop \tilde
L')\:(h_{(1)}\ra v)\\[.2em]&=
(h_{(2)}\ra^\cop (\tilde L\circ \tilde L'))(h_{(1)}\ra v)
\end{split}
\end{equation}
where in the last equality we used the definition \eqref{copadjact} of $\ra^\cop$.  
This shows $\tilde L\circ \tilde L'\in {}_\Bbbk\hom(V,Z)$.

The tensor product of morphisms in ${}^H\MMM$ is also the usual one, for all
$f\in \Hom_{{\!\:}^{H\!\!\:}\MMM}(V,W), f'\in 
\Hom_{{\!\:}^{H\!\!\:}\MMM}(V',W')$,  $(f\otimes f')(v\otimes
v'):=f(v)\otimes f'(v')$, for any $v\in V,v'\in V'$.
On the other hand the tensor product of internal morphisms 
differs from that of the category of $\Bbbk$-modules
(cf. \cite[Coroll. 9.3.16]{Ma95}).
Given $\Bbbk$-linear maps $L \in\hom_\Bbbk(V, W)$, $L' \in \hom_\Bbbk(V', W')$ the tensor product $L\otimes_{\mathcal
  R} L' $ is the $\Bbbk$-linear map 
\begin{equation}\label{otimesR}
L\otimes_\RR L' := (L\circ \bar R^\alpha\ra\,)\otimes (\bar R_\alpha\RA 
L')\in \hom_\Bbbk(V\otimes V', W\otimes{W'})~,
\end{equation}
i.e., $(L\otimes_\RR L')(v\otimes v')=L(\oR^\al\ra
v)\otimes (\oR_\al\RA L')(v')$, for any $v\in V,v'\in V'$. It is this associative tensor product that is 
compatible with the $H$-module structure:
$\hxi \ra 
(L\otimes_\RR L')=\hxi_{(1)}\ra L\otimes_\RR  \hxi_{(2)}\ra L'$.
From the definition it follows that
\begin{equation}\label{POQPcircQ}
L\otimes_\RR L'= (L\otimes \id)\circ (
\oR^\alpha\ra\;\otimes \,\oR_\alpha\RA L')
=(L\otimes_\RR \id)\circ (\id\otimes_\RR L')~.
\end{equation}
While
$L\otimes_\RR\id=L\otimes \id$, we have $
\id\otimes_\RR L'=
\oR^\alpha\ra\;\otimes \,\oR_\alpha\RA L'=\tau\circ (L'\otimes \id ) \circ
\tau^{-1}~.$

Similarly, it can be proven that given $\Bbbk$-linear maps $\tilde L \in{}_\Bbbk\hom(V, W)$,
$\tilde L'\in {}_\Bbbk\hom(V', W')$ we have the corresponding tensor
product $\tilde L\,\tilde\otimes_\RR \,\tilde L' $
\begin{equation}\label{tildeotimesR}
\tilde L\,\tilde\otimes_\RR\, \tilde L' :=  (\bar R^\alpha\RA^\cop 
\tilde L)\otimes (\tilde L'\circ \bar R_\alpha\ra\,)\in
{}_\Bbbk\hom(V\otimes  V',W\otimes W')~,
\end{equation}
i.e., for all $v\in V$, $v'\in V'$, $(\tilde L\,\tilde\otimes_\RR \,\tilde
L') (v\otimes v')= (\oR^{\alpha}\RA^\cop\tilde L)(v) \otimes \tilde L'
(R_\alpha\ra v')$. This tensor product is associative and compatible
with the $H$-adjoint action $\RA^\cop$, that is, we have
$\hxi\RA^\cop (\tilde L\,\tilde\otimes_\RR \tilde L') =(\hxi_{(1)}\RA^\cop
\tilde L)\,\tilde\otimes_\RR (\hxi_{(2)}\RA^\cop \tilde L')$. 
\\

\subsubsection{Algebras and bimodules}\label{ALGBIMOD}
Let $H$ be a Hopf algebra. A left $H$-module algebra $A$ is an algebra with a compatible
$H$-module structure,
$$
\hxi\rhd(a b)=(\hxi_{(1)}\rhd a)(\hxi_{(2)}\rhd b) ~,~~
\hxi\rhd 1_{A}=\epsilon(\hxi)1_{A}
$$
for all $\hxi\in H$ and $a,b\in {A}$. 
We denote by $^{H}_{\AA}\MMM^{}_\AA$ the category of
$H$-equivariant $A$-bimodules.
An object $V$ in $^{H}_{\AA}\MMM^{}_\AA$ is an
$A$-bimodule with a compatible $H$-module structure, i.e., $H\otimes V\to
V$, $\hxi\otimes v\mapsto \hxi \ra v\,;~  \hxi \ra av =(\hxi_{(1)}\ra a)(\hxi_{(2)}\ra v)$
and similarly for the right $A$-module structure.
Morphisms in
$^{H}_{\AA}\MMM^{}_\AA$ are $H$-equivariant maps that are also
$A$-bimodule morphisms. 

The category  $^{H}_{\AA}\MMM^{}_\AA$ becomes a monoidal category
with the balanced tensor product $\otimes_A$ (where  by definition $V\otimes_A W$, with $V,W$ in
$\catm$,  is the quotient of $V\otimes W$ in $\catm$  
with the obvious left and right $A$-actions inherited from those of
$V$ and $W$ respectively).

If $V,W$ are modules in $\catm$ also $\hom_\Bbbk(V,W)$ and
${}_\Bbbk\hom(V,W)$ are modules in $\catm$  with $H$-action as in
\eqref{adjact} and \eqref{copadjact} respectively. The $A$-bimodule
structure of  $\hom_\Bbbk(V,W)$ is given via the left $A$-module structure of $V$
and $W$, that of ${}_\Bbbk\hom(V,W)$ via the right $A$-module
structure of 
$V$ and $W$: For all $a\in A, v\in V, L\in
\hom_\Bbbk(V,W), \tilde L\in {}_\Bbbk\hom(V,W)$,
\begin{equation}
\label{aLvLav}
(aL)(v)=a(L(v))~,~~(La)(v)=L(av)~,
\end{equation}
\\[-2.6em]\begin{equation}\label{AhomAA}
(\tilde L a)(v)=\tilde
L(v)\:\!a~,\,\,\,\,~~(a\tilde L)(v)=\tilde L(va)~.
\end{equation} 
Let  $\hom_A(V,W)\subset \hom_\Bbbk(V,W)$ be the  submodule in $\catm$ of right $A$-linear maps: for all $a\in A$, $L(va)=L(v)a$,
and let 
${}_A\hom(V,W)\subset {}_\Bbbk\hom(V,W)$ be the submodule in $\catm$
of left $A$-linear maps: for all $a\in A$, $\tilde L(av)=a\tilde L(v)$.

Associated with $\hom_A(V,W)\subset \hom_\Bbbk(V,W)$  and ${}_A\hom(V,W)\subset {}_\Bbbk\hom(V,W)$  we have the  functors $\hom_A: (\catm)^\op\times \catm\to \catm$ and
${}_A\hom: (\catm)^\op\times \catm\to \catm$, with the action on morphisms 
$(f^\op,g)$ in $(\catm)^\op\times \catm$ as in
\eqref{eqn:internalhommorph}. 
The functor $\,\mbox{-}\otimes_A W$ is left adjoint to
$\hom_A(W,\mbox{-})$ \cite{BSS1},
and similarly,
 $V\otimes_A \mbox{-}$ is left
adjoint to ${}_A\hom(V,\mbox{-})$,
 thus $$(\catm, \otimes_A, \hom_A,{}_A\hom)$$ is a
biclosed  monoidal category. 
\\[-.4em]

When $H$ has a triangular structure $\RR$  we consider $\AA$ to be  braided commutative (also called symmetric  or
quasi-commutative) if, for all $a,b\in \AA$, 
\begin{equation}\label{bsalgA}
  ab=(\bar R^\alpha\ra b)(\bar R_\alpha \ra a)~.
  \end{equation}
Similarly, $V$ in $^{H}_{\AA}\MMM^{}_\AA$  is symmetric if
\begin{equation}\label{bsmodV}
  a\:\!v =
  (\oR^{\al}\ra v )  (\oR_{\al}\ra a)~.
  \end{equation}
We denote by $\cat$ the
full subcategory of symmetric modules in $^{H}_{\AA}\MMM^{}_\AA$. 
Let $V,W$ be modules in $\cat$, then $V\oA W$, $\hom_A(V,W)$
and ${}_A\hom(V,W)$ are also in $\cat$; for example it is easy to see that for all $a\in A$, $L\in \hom_A(V,W)$,  $\tilde L\in
{}_A\hom(V,W)$, $La=(\oR^\al \ra a)(\oR_\al \RA L)$ and $\tilde La=(\oR^\al \ra
a)(\oR_\al \RA \tilde L)$. 
Extending to the biclosed case the results of \cite{BSS1} we have that
$$(\cat, \otimes_\AA, \hom_A, {}_A\hom)~$$ is a full closed monoidal
subcategory of $(\catm, \otimes_A, \hom_A, {}_A\hom)$ which is
symmetric.

The braiding is induced from that in ${}^H\MMM$ and the isomorphisms in \eqref{LefttoRight} restrict to
isomorphisms 
\begin{equation}\label{Dcallr}
{\cal D}_{V,W}: \hom_A(V,W)\to {}_A\hom(V,W)
\end{equation}
 in $\cat$,
thus proving that $\hom_A$ and ${}_A\hom$ are naturally isomorphic
functors (cf. \cite[\S 5.6]{AS},
the $A$-module actions \eqref{AhomAA} there are
denoted $\cdot^\op$, so that e.g. $a\cdot^\op \tilde L\,(v)=\tilde L(va)$).
\\

Since $(\cat, \otimes_\AA, \hom_A)$ is a symmetric closed monoidal
category, besides the usual evaluation and
composition, we have the tensor product of  morphisms, denoted $\oA$. For all
$f\in \Hom_{\cat}(V,W), f'\in 
\Hom_{\cat}(V',W')$,
$(f\oA f')(v\oA
v'):=f(v)\oA f'(v')$, for any $v\in V,v'\in V'$.
We also have the tensor product of internal morphisms
that with slight abuse of notation we
still denote $\otimes_\RR$. Indeed, similarly to the braiding, this can be
seen as induced from the tensor product of internal morphisms  in
$({}^H\MMM, \otimes, \hom_\Bbbk)$. Let $L\in\hom_A(V,W)\subset \hom_\Bbbk(V,W)$ and $L'\in
\hom_A(V',W')\subset  \hom_\Bbbk(V',W')$, then 
$L\otimes_\RR L'\in
\hom_\Bbbk(V\otimes V',W\otimes W')$, as defined in \eqref{otimesR}
is trivially right $A$-linear and
induces a well-defined right $A$-linear map 
in $\hom_A(V\otimes_A V',W\otimes_A W')$ that we still denote
$L\otimes_\RR L'$  (cf. \cite[Thm. 5.16]{AS}).
Associativity
is straightforward. Moreover, for each quadruple
$V,W,V',W'$ of modules in $\cat$ the map
\begin{equation*}{\otimes_\RR}^{}_{{V,W,V',W'}}:
\hom_A(V,W)\otimes_A\hom_A(V',W')\to\hom_A(V\otimes_A
V',W\otimes_A W')~,~~L\otimes_AL'\mapsto L\otimes_\RR L'
\end{equation*}
is a
morphism in the category (cf.  \cite[Prop. 9.3.13]{Ma95}, \cite[\S 5.6]{BSS1}).

Similarly, in $(\cat, \otimes_\AA, {}_A\hom)$ the tensor product of
internal morphisms is denoted $\tilde\otimes_\RR$ and can be seen as induced from \eqref{tildeotimesR}. Here too for each quadruple
$V,W,V',W'$ of modules in $\cat$ we have that
$${\tilde\otimes_\RR}{}^{\phantom{J_i}}_{{V,W,V',W'}}:
{}_A\hom(V,W)\otimes_{\!\!\:A}{}_{\!\:A}\hom(V',W')\to{}_A\hom(V\otimes_A
V',W\otimes_A W')~,~~\tilde L\otimes_A\tilde L'\mapsto \tilde
L_{\,}\tilde\otimes_\RR\tilde L'$$ is a morphism in $\cat$.

\subsubsection{Graded algebras and bimodules}
The results of the previous subsection 
can be extended to the case of $\mathbb{Z}$-graded modules
$V=\bigoplus_{n\in \mathbb{Z}}V^n$.
We consider $H$ to be
$\mathbb{Z}$-graded and nontrivial only in degree zero. Let $\OMM$ be
a graded algebra and an $H$-module algebra, it 
is graded braided commutative if
$$\theta\wedge\theta'=(-1)^{|\theta|_{}|\theta'|}(\bar R^\alpha\ra^\cop \theta')\wedge (\bar R_\alpha
\ra^\cop \theta)~,$$ where $\wedge$ denotes the product in $\OMM$ and
$\theta,\theta'$ are arbitrary elements in $\OMM$ of  homogeneous degree $|\theta|$ and $|\theta'|$.
Correspondingly,  $\CatM$ denotes the category of $\mathbb{Z}$-graded
modules that are relative
$H$-modules $\OMM$-bimodules (with  grade compatible $\OMM$-module
actions)  and that are graded
symmetric:
$V=\bigoplus_{n\in \mathbb{Z}}V^n$ is in  $\CatM$ if 
$$\theta\,v =
(-1)^{|\theta|_{}|v|}(\oR^{\al}\ra v )  (\oR_{\al}\ra^\cop \theta)~,$$ where
$|v|$ is the degree of the homogeneous element $v\in V$.

The category $\CatM$ is monoidal with tensor product $\otimes_\OMM$. It is also
closed, indeed first observe that for each $V,W$ in $\CatM$, we have that ${}_\Bbbk\hom(V,W)$
and $\hom_\Bbbk(V,W)$ are naturally graded $H$-modules. Then define
$$
{}_\OMM\hom(V,W)
$$
to be the graded $H$-submodule of ${}_\Bbbk\hom(V,W)$ spanned by
{\sl{graded}} left $\OMM$-linear maps; these are maps
 $\tilde L\in  {}_\Bbbk\hom(V,W)$ of homogenous degree $|\tilde L|$
 (i.e. $\tilde L:V^n\to
 W^{n+|\tilde L|}$, $n\in \mathbb{Z}$) such that
\begin{equation}\label{gradedOMlinear}
\tilde L(\theta v)=(-1)^{|\tilde L||\theta|}\theta \:\!\tilde
 L(v)~.
\end{equation} 
Similarly, $\hom_\OMM(V,W)$ is the module of right $\OMM$-linear
 maps, $L(v\theta)=L(v)\theta$.
The modules $ {}_\OMM\hom(V,W)$  and $\hom_\OMM(V,W)$ are  in $\CatM$; their $\OMM$-bimodule
structure reads, for all $L\in \hom_\OMM(V,W)$, and for all $\tilde L\in {}_\OMM\hom(V,W)$,  $\theta\in \OMM$, $v\in V$,
respectively of homogenous degree
$|\tilde L|$, $|\theta|$ and $|v|$,
\begin{equation}\label{bimodofOM}
\begin{split}&(\theta L)(v)=\theta (L(v))~,~~~~~~~~~~~~~(L\theta)(v)=L(\theta v)~,\\[.4em]
&(\tilde L \theta)(v)=(-1)^{|\theta||v|}\tilde
L(v)\!\:\theta~,~~~~(\theta\tilde L)(v)=(-1)^{|\theta|(|\tilde L|+|v|)}\tilde L(v\theta)~.
\end{split}
\end{equation}
This defines the functors $\hom_\OMM: (\CatM)^\op\times \CatM\to \CatM$ and
${}_\OMM\hom: (\CatM)^\op\times \CatM\to \CatM$, where their action on
morphisms  
$(f^\op,g)$ in $(\CatM)^\op\times \CatM$ is as in \eqref{eqn:internalhommorph}.

Similarly to the ungraded case we have the symmetric biclosed monoidal
category
$$(\CatM,\otimes_\OMM, \hom_\OMM, {}_\OMM\hom)$$
with tensor product of internal morphisms respectively denoted 
$\otimes_\RR$ and $\tilde\otimes_\RR$.

\subsection{Finitely generated projective modules and their duals}\label{fgp}
In this section $H$ is a Hopf algebra  (no
triangularity structure is assumed).  We study finitely generated and
projective left (right) $A$-modules in $\catm$ and their duals. It is known that finitely
generated and projective $\Bbbk$-modules are rigid modules. Similarly, finitely generated
and projective right (left) $A$-modules in $\catm$ are rigid
modules in $\catm$.  We continue the study of internal morphisms  proving key canonical isomorphisms
for internal morphisms in $\catm$ arising from tensor products with a rigid
module (cf.~{\it 2$\,$} in Theorem \ref{Thmfgp123}).

Among the various equivalent definitions of finitely generated
projective module (see e.g.~the monograph
\cite{LamBook}) we use the convenient characterization in terms of a pair of dual bases 
\begin{lemma}  ({\it Dual Basis Lemma}).
 Let $A$ be an algebra. 
A left $A$-module $\Gamma$ is
finitely generated and projective if and only if there exists a family of elements
$\lbrace s_i\in \Gamma:i=1,\dots,n\rbrace$ and left $A$-linear maps
$\lbrace {}^{*\!}s^i \in\stG
:={_A}\Hom(\Gamma,A):i=1,\dots,n \rbrace$
with $n\in\mathbb{N}$, such that for any $s\in \Gamma$ we have
(sum over repeated indices understood)
\begin{equation}\label{astsisi}
 s= {}^{*\!}s^i(s)\, s_i~.
\end{equation}
A right $A$-module $\Sigma$  is
finitely generated and projective if and only if there exists a family of elements
$\lbrace \sigma^i\in \Sigma :i=1,\dots,n\rbrace$ and right $A$-linear maps
$\lbrace \sigma_i^*\in \Sigma^*:=\Hom_A(\Sigma,A):i=1,\dots,n \rbrace$
with $n\in\mathbb{N}$, such that for any $\sigma\in \Sigma$ we have
\begin{equation}
 \sigma=\sigma^i\, \sigma^*_i(\sigma)~.
\end{equation}
\end{lemma}
The set $\lbrace s_i, {}^{*\!}s^i:i=1,\dots,n \rbrace$ is loosely
referred to a ``pair of dual bases'' for the left $A$-module $\Gamma$,
even though $\lbrace s_i\rbrace$ is just a generating set of $\Gamma$
and not necessarily a basis. Similarly $\lbrace \sigma^i,
{\sigma}{}^*_i:i=1,\dots,n \rbrace$ is a pair of dual bases for the right $A$-module $\Sigma$.
\\

The dual $\stG :={_A}\Hom(\Gamma,A)$ of a finitely generated and projective left $A$-module $\Gamma$ is a finitely
generated and projective right $A$-module, with right $A$-action
 as in \eqref{AhomAA}. Moreover, the dual $({\stG})^{*_{}}
 :=\Hom_A(\stG,A)$  of the dual is a left $A$-module canonically
identified with the original module $\Gamma$. Similarly, we have the
left $A$-module
$\Sigma^*:=\Hom_A(\Sigma,A)$ dual to the right $A$-module $\Sigma$ and the canonical
identification ${}^*(\Sigma^*):= {}_A\Hom(\Sigma^*,A)\simeq
\Sigma$. We state these properties for left $A$-modules (the proof can be easily derived from e.g. ~\cite[\S 2B]{LamBook}).
\begin{proposition}\label{propo:fgpmodprop}
 Let $\Gamma$ be a finitely generated and projective left
 ${A}$-module. Denote by
$\lbrace s_i,\sts^i:i=1,\dots,n \rbrace$ a pair of dual bases. For any $s\in \Gamma$, let
$\iota (s) \in (\stG)^*:= \Hom_A(\stG,A)$ be defined by
$\iota(s)(\sts) := \sts(s)$, for all $\sts \in \stG$. We have
\begin{enumerate}
 \item[{\it 1.}$\phantom i\!\!$] $\lbrace \sts^i,\iota (s_i):i=1,\dots,n\rbrace$ is a pair of dual bases for $\stG$,
 \item[{\it 2.}$\phantom i\!\!$] $\stG$ is a finitely generated and projective right $A$-module,
 \item[{\it 3.}$\phantom i\!\!$] The canonical map $\iota: \Gamma\to ({\stG})^{*_{}}\,,~s\mapsto
  \iota({s})$ is an isomorphism of left $A$-modules.
\end{enumerate}
\end{proposition}
If we consider modules in $\catm$, Proposition
\ref{propo:fgpmodprop} holds in $\catm$. Let  $\Gamma$ be a module in
$\catm$ that is finitely generated and projective as left $A$-module, then
\begin{enumerate} 
\item[{\it  1$'\!.$}] $\stG:={}_A\hom(\Gamma, A)$
is in $\catm$ and  is finitely generated and
projective as right $A$-module,
\item[{\it 2$\:\!'\!.$}] 
$({\stG})^{*_{}}
 :=\hom_A(\stG,A)$
is in $\catm$ and is finitely
generated and projective as  left $A$-module, 
\item[{\it 3$\!\:'\!.$}] The canonical
map $\iota : \Gamma\to ({\stG})^{*_{}}\,,~s\mapsto \iota(s)$ is an isomorphism
in $\catm$. 
\end{enumerate}

We further recall (cf. e.g. \cite[Prop. 6.17]{AS}), 
\begin{proposition}\label{propo:projhomiso}
Let $\Gamma$ be a finitely generated and projective left $A$-module, so
that $\stG$ is a finitely generated and projective right $A$-module. Let $W$ be an $A$-bimodule.
Then there exist  right $A$-module and left $A$-module isomorphisms
(evaluation maps) 
\begin{equation}\label{jmathdef}
\begin{split}
\flat : {}^*\Gamma\otimes_AW\, &\to\, {}_A\Hom(\Gamma,W)\\
\sts\otimes_A w\,\;&\mapsto\, (\sts\otimes_A w)^\flat(s):=
\sts(s)\,w~~~~
\end{split}
\begin{split}
\imath :  W\otimes_A \Gamma\,& \to\, \Hom_A(\stG,W)\\
w\otimes_A s\;\,&\mapsto\,\imath(w\otimes_A s)\,(\sts):= w\:\sts(s)\;.~~
\end{split}
\end{equation}
The inverses are $\tilde L\mapsto
{\tilde L}^{\flat^{-1}\!}= \sts^i\otimes_A\tilde L(s_i)$
and $ L\mapsto
\imath^{-1}(L)= L(\sts^i)\otimes_A s_i$ (sum on $i$ understood).
If in addition $\Gamma$ and $W$ are modules in $\catm$
then 
$$\flat: {}^*\Gamma\otimes_AW \to {}_A\hom(\Gamma,W)~~,~~~\imath :  W\otimes_A \Gamma\, \to\, \hom_A(\stG,W)$$
are module isomorphisms in $\catm$.
\end{proposition}

A pairing between modules $\Gamma$, $\Sigma$ in $\catm$ is a
morphism $\Gamma\otimes_A\Sigma\to A$. We denote by 
\begin{equation}\label{ppair}
\le~,~\re:\Gamma\otimes_A\stG\rightarrow A~,~~s\otimes_A \sts\mapsto \le 
s,\sts\re:=\sts(s)=\iota(s)(\sts) 
\end{equation}
the pairing due to the evaluation of ${}^*\Gamma={}_A\hom(\Gamma,A)$
on $\Gamma$.
It is well-defined on the balanced tensor product $\otimes_A$ because of the left
$A$-module structure of $\stG$. It is easily seen to be left and right $A$-linear,
indeed the notation $\le~,~\re$ conveniently takes into account the $A$-bimodule structures of
$\Gamma$ and
$\stG$, as well as that $\stG$ are left $A$-linear maps while
$\Gamma\simeq({\stG})^*$ are right $A$-linear maps:
For all $a\in A, s\in \Gamma, \sts\in \stG$, $\langle a
s,\sts\rangle=a\langle  s,\sts\rangle$, 
$\langle 
s a,\sts\rangle=\langle  s,a\:\! \sts\rangle$,
$
\langle s,\sts a\rangle=\langle  s,\sts\rangle a
$.
Furthermore, $H$-equivariance
\begin{equation}\label{Hppair}
h\ra\le s,\sts\re=\le h_{(1)} \ra s, h_{(2)}\ra^\cop \sts\re~ 
\end{equation}
is due to the $H$-module structure of $\stG$, cf. \eqref{copadjact2}.

We extend the pairing $\le~,~\re : \Gamma\otimes_A\stG\to
A$ to the morphisms in $\catm$
\begin{equation}\label{convimath}\le ~\,,\,~~\re: \Gamma\otimes_A \stG\otimes_A W\to 
W~,~~s\otimes_A \sts\otimes_A w\mapsto \le s,\sts\otimes_A w\re:=\le 
s,\sts\re w
\end{equation}\\[-2.5em]
\begin{equation}\label{convjmath}\le ~~\,\,,~\re: W\otimes_A \Gamma\otimes_A \stG \to 
W~,~~w\otimes_A s\otimes_A \sts\mapsto \le w\otimes_A s,\sts\re:=w\le 
s,\sts\re
\end{equation}
so that the internal morphisms $(\sts\otimes_A w)^\flat$  and  $\imath(w\otimes_A s)$ of
Proposition \ref{propo:projhomiso} are respectively simply denoted as
$\le~\;\,,\sts\otimes_A w\re$ and
$\le w\otimes_A s\!\:,~\;\re$.

\begin{definition}\label{defdualcat}
Given modules $\Gamma\in \catm$, $\Sigma\in \catm$ we say that $\Gamma$
has right dual $\Sigma$, or equivalently that $\Sigma$ has left dual $\Gamma$, if we have maps 
$\ev:\Gamma\otimes_A\Sigma\rightarrow A$ (evaluation map) and
$\coev:A\rightarrow \Sigma\otimes_A\Gamma$ (coevaluation map) in $\catm$
such that the compositions
\begin{equation}\label{coerecoev}
\begin{split}
\end{split}
\begin{split}\Gamma\simeq \Gamma\otimes_A
A
\xrightarrow{\id_\Gamma\otimes_{\!A\,}\coev\,}\,
\Gamma\otimes_A\Sigma\otimes_A\Gamma\xrightarrow{\ev\otimes_{\!A\,}\id_\Gamma\,}\Gamma\,,\\[.4em]
\Sigma\simeq
A\otimes_A\Sigma
\xrightarrow{\coev\otimes_{\!A\,}\id_{\Sigma}}
\Sigma\otimes_A\Gamma\otimes_A\Sigma\xrightarrow{\id_{\Sigma\,}\otimes_{\!A\,} 
\ev}\Sigma\,,
\end{split}
\end{equation}
are respectively the identity maps 
$\id_\Gamma$ and $\id_{\Sigma}$.
If $\Gamma$ has a right dual we say that it is right rigid.  
If $\Sigma$ has a left dual we similarly  say that it is left rigid.  
\end{definition}
We will also denote the
evaluation and coevaluation maps of the right rigid module $\Gamma$
 by $\ev^{}_\Gamma: \Gamma\otimes_A \Sigma \to A$ and
 $\coev^{}_\Gamma: A\to \Sigma\otimes_A \Gamma$,
and  frequently use the notations $\le~,~\re:\Gamma\otimes_A\Sigma\to
A$ and $\le~,~\re_\Gamma$
for the evaluation map.

Left (right)  duals are unique up to
isomorphisms, in this sense we can simply speak of {\sl the} left (right) dual
of a module, and we say that the module is rigid. 
This implies that if there exists a coevaluation map, it
is uniquely determined by the pairing $\ev: \Gamma\otimes_A \Sigma\to
A$. 
A pairing $\le~,~\re:\Gamma\otimes_A \Sigma\to A$  is
{\it exact} if there exits a map $\coev :A\to
\Sigma\otimes_A\Gamma$ fulfilling the conditions of Definition
\ref{defdualcat}.

The pairing in \eqref{ppair} is exact if and only if $\Gamma$ in
$\catm$ is finitely generated and projective as left $A$-module.

\begin{theorem}\label{Thmfgp123} Let $\Gamma$ in $\catm$ and  $\stG:={}_A\hom(\Gamma,A)$.\begin{enumerate}
\item  $\Gamma$ is finitely generated and projective as left
  $A$-module in $\catm$
if and only if  $\Gamma$ is right rigid.
\item If $\Gamma$ is right rigid, 
 for all $V,W\in \catm$,
 \begin{equation}\label{iso1e2}
{}_A\hom(\Gamma\otimes_A V,W)\simeq {}_A\hom(V,{}^*\Gamma\otimes_A 
  W)~~,~~~~\hom_A(V\otimes_A\stG,W)\simeq 
  \hom_A(V,W\otimes_A\Gamma)~~
\end{equation}
  are isomorphisms in $\catm$. 
\item Let $\Upsilon$ be also right rigid
  in $\catm$ and  ${}^*\Upsilon:={}_A\hom(\Upsilon,A)$, then so is
  $\Gamma\otimes_A\Upsilon$, and we
further have the isomorphism in $\catm$, $\stU\otimes_A \stG\simeq
{}^*(\Gamma\otimes_A \Upsilon)$.
\end{enumerate}
\begin{proof} 
 {\it 1.} Let $W=\Gamma$ in Proposition \ref{propo:projhomiso}. Since $\id_\Gamma\in\Hom_\catm(\Gamma,\Gamma)\subset
 {}_A\hom(\Gamma,\Gamma)$ is central, i.e., $a_{\,}\id_\Gamma=\id_{\Gamma\,}
 a$ (recall the $A$-bimodule structure \eqref{AhomAA} of
 ${}_A\hom(\Gamma,\Gamma)$) so it its image
 $\id_\Gamma^{\flat^{-1}}=\sts^i\otimes_A s_i \in
 \stG\otimes_A\Gamma$. Thus
\begin{equation}\label{expcoev}
\coev: A\rightarrow \stG\otimes_A\Gamma~,~~  a\mapsto a_{\,}\sts^i\otimes_A s_i=\sts^i\otimes_A s_{i\,}a
\end{equation}
is a well defined $H$-equivariant $A$-bimodule map.
The coherence conditions \eqref{coerecoev} follow using \eqref{expcoev}
and recalling that $\lbrace s_i, {}^{*\!}s^i \rbrace$ and $\lbrace  {}^{*\!}s^i, \iota( s_i)\rbrace$ are
a pair of dual bases for $\Gamma$ and $\stG$, respectively
(cf. Proposition \ref{propo:fgpmodprop}). Vice versa, if $\Gamma$ is
rigid then \eqref{expcoev} determines a dual basis.
\\[1em]
{\it 2.}
For all $\tilde L\in {}_A\hom(\Gamma\otimes_A V,W)$, define
\begin{equation}\label{sharpmap}{\tilde
  L}^\sharp:=(\id_{\stG}\otimes_A\tilde L)\circ (\coev\otimes_A\id_V):
V\to \stG\otimes_A W~,~~v\mapsto {\tilde L}^\sharp(v)=\sts^i\otimes_A\tilde
L(s_i\otimes_A v)~.
\end{equation}The map ${\tilde L}^\sharp$ is  left $A$-linear
because so is  $\tilde L$ and because $\coev$ is
$A$-bilinear. 
Hence we have a well-defined map
\begin{equation*}
  \sharp: {}_A\hom(\Gamma\otimes_A V,W)\to
  {}_A\hom(V,{}^*\Gamma\otimes_AW)
\end{equation*}
which is $H$-equivariant since so is
$\coev$. Its inverse is $\flat:
{}_A\hom(V,{}^*\Gamma\otimes_A W)\to  {}_A\hom(\Gamma\otimes_A V,W)$,
$s\otimes_A v\mapsto {\tilde P}^\flat(s\otimes_A v)=\le s, {\tilde
  P}(v)\re$, 
(recall equation \eqref{convimath}).  Indeed we have
${\tilde L}^{\sharp^{\mbox{\scriptsize$\,{}\flat$}}} \!(s\otimes_A v)=\le s, {\tilde 
  L}^\sharp(v)\re=\le s,\sts^i\otimes_A\tilde L(s_i\otimes_Av)\re=\le 
s, \sts^i\re \tilde L(s_i\otimes_A v)=\tilde L(\le s,\sts^i\re 
s_i\otimes_A v)=\tilde L(s\otimes_A v)
$
where in the third equality we used left $A$-linearity of $\tilde L$.
We similarly have
$
{{\tilde P}^{\flat^{\mbox{\scriptsize$\,{}\sharp$}}}} (v)=
\sts^i\otimes_A\tilde P^\flat (s_i\otimes_A v)=\sts^i\otimes_A\le
s_i,\tilde P(v)\re=\tilde P(v)$.
The second isomorphism in \eqref{iso1e2} is similarly proven.
\\[1em]
{\it 3.}
Let 
$\{u_j,\stu^j : j=1,\ldots,m\}$ be a dual basis of $\Upsilon$, then a dual
basis of  $\Gamma\otimes_A\Upsilon$
is given by the elements $t_{ij}=s_i\otimes_A u_j\in
\Gamma\otimes_A\Upsilon$, and the elements $\stt^{ij}\in
{}^*\!(\Gamma\otimes_A\Upsilon)$ defined by, for all $s\in \Gamma,
u\in \Upsilon$,  $\stt^{ij}(s\otimes_A
u) :=\le s\le u,\stu^j\re,\sts^i\re$.

In order to prove the isomorphism $\stU\otimes_A\stG\simeq \stGU$,
define
$\varphi^{}_{\Gamma,\Upsilon}:\stU\otimes_A\stG \to \stGU$,
$\stu\otimes_A\sts\mapsto
\varphi^{}_{\Gamma,\Upsilon}(\stu\otimes_A\sts)=\stt^{ij}\le\le t_{ij},\stu\re,\sts\re
$.
The map $\varphi^{}_{\Gamma,\Upsilon} $ is easily seen to be in
$\catm$; its inverse is $\varphi^{-1}_{\Gamma,\Upsilon}:\stGU\to\stU\otimes_A\stG $,
$\stt\mapsto
\varphi^{-1}_{\Gamma,\Upsilon}(\stt)=\stu^j\otimes_A\sts^i\le s_i\otimes_Au_j,\stt\re$.
\end{proof}\end{theorem}

\begin{remark}\label{newevcoev}
Due to the isomorphism  $\stU\otimes_A \stG\simeq
{}^*(\Gamma\otimes_A \Upsilon)$, in the following the right dual of
$\Gamma\otimes_A\Upsilon$ will be considered to be $\stU\otimes_A
\stG$, with evaluation and coevaluation maps
\begin{equation*}
\begin{split}
\le~~,~~\re &:  \Gamma\otimes_A \Upsilon\otimes_A \stU\otimes_A
\stG\rightarrow A~,~~ s\otimes_A u \otimes_A \stu\otimes_A\sts \mapsto 
\le s\otimes_A u, \stu\otimes_A\sts\re=\le s\le u, \stu\re,\sts\re_{}~,\\[.4em]
\coev&: A\to  \stU\otimes_A\stG\otimes_A \Gamma\otimes_A
\Upsilon~,~~a\mapsto a\,\stu^j\otimes_A\sts^i\otimes_As_i\otimes_A u_j~,
\end{split}
\end{equation*}
where $\{s_i,\sts^i:i=1,\dots,n \}$ and $\{u_j,\stu^j :j=1,\dots,m \}$
are dual basis
of $\Gamma$ and $\Upsilon$, respectively.
\end{remark}

\begin{remark}\label{rem:gmc}
Definition \ref{defdualcat} holds in a generic monoidal category. For
example also when the tensor product is topological. In this case
rigid modules are projective and topologically finitely generated
(i.e., there exists        a finite number of
elements  that span a dense subset of the module). This
is the case of Example \ref{ex:bundleformal}.
\end{remark}

The right dual (or adjoint) of a morphism $f:
\Gamma\to \Upsilon$ of right rigid modules in $\catm$ is
the morphism ${}^*f:\stU\to \stG$ defined by $\le s, {}^{*\!}f(\stu)\re
= \le f(s),\stu\re$, for all $s\in \Gamma, {}^*u\in {}^*\Upsilon$. Explicitly, ${}^{*\!}f(\stu)=\sts^i\,\le
f(s_i),\stu\re$
(using a pair of dual bases for $\Gamma$).
Similarly, the left dual of a morphism
$g:\Sigma \to \Lambda$  of left rigid modules in $\catm$ is
the morphism $g^\ast:\Lambda^\ast\to \Sigma^\ast$ defined by $\le g^\ast(\lambda^\ast), \sigma\re
= \le \lambda^\ast,g(\sigma)\re$  for all $ \lambda^\ast\in
\Lambda^\ast , \sigma\in \Sigma$. 
Duals of internal morphisms
will be studied in Section \ref{dualinternalmorph}.

\subsection{Compact closed categories\label{sect:ribbon}}
If the Hopf algebra $H$ has a triangular structure $\RR$ and $A$ is
braided commutative (cf. \eqref{bsalgA}) 
we can consider finitely generated (left or right) $A$-modules in the
symmetric category $\cat$ of symmetric modules; we
recall that 
by definition these satisfy $av= (\oR^\al\ra v)(\oR_\al\ra a)$ for all $a\in A, v\in
V$, cf. \eqref{bsmodV}.

Let $\Gamma$ in $\cat$ be finitely generated and
projective as a left $A$-module, hence it is right rigid in $\cat$, 
with right dual ${}^*\Gamma={}_A\hom(\Gamma,A)$ which is finitely
generated and projective as a right $A$-module  in $\cat$. 
All the results in section \ref{fgp} concerning  modules in $\catm$
holds true in the full subcategory $\cat$ of symmetric modules.

In a symmetric category a right rigid module $\Gamma$ is also
left rigid, and vice versa a left rigid module is also right rigid.
We give an explicit proof for our category of interest $\cat$, where
we recall that the braiding $\tau$ is defined in \eqref{taubraiding}.
\begin{proposition}\label{propep'}
Let $\Gamma \in \cat$ be right rigid, with right dual ${}^\ast \Gamma$,
evaluation $\le~,~\re:\Gamma\oA {}^\ast\Gamma\to A$,
and coevaluation $\coev:A\rightarrow {}^\ast\Gamma\oA\Gamma$.
Then $\Gamma$ is left rigid with left dual  ${}^\ast\Gamma$ and evaluation
and coevaluation maps
\begin{equation}\label{tevcoev}
\begin{split}
\le~,~\re':=\langle~,~\rangle\circ \tau^{}_{\stG,\Gamma}&:
\stG\otimes_A\Gamma\to A~,\\[.4em]
\coev':=
\tau^{}_{\stG,\Gamma}\circ\coev&: A\to \Gamma\otimes_A\stG~.
\end{split}
\end{equation}
\end{proposition}\begin{proof}
We have to show that the compositions
\begin{equation}\begin{split}
    \Gamma\simeq
A\otimes_A\Gamma
\xrightarrow{\coev'\otimes_{\!A\,}\id_{\Gamma}}
\Gamma\otimes_A{}^\ast\Gamma\otimes_A\Gamma\xrightarrow{\id_{\Gamma\,}\otimes_{\!A\,} 
\le~,~\re'}\Gamma\,,
~~\,\\[.2em]
{}^\ast\Gamma\simeq {}^\ast\Gamma\otimes_A
A
\xrightarrow{\id_{{}^\ast{}\!\!\:\Gamma}\otimes_{\!A\,}\coev'\,}\,
{}^\ast\Gamma\otimes_A\Gamma\oA{}^\ast\Gamma\xrightarrow{\le~,~\re'\otimes_{\!A\,}\id_{{}^\ast\!\!\:\Gamma\,}}{}^\ast\Gamma\,,
\end{split}
\end{equation}
equal $\id_\Gamma$ and $\id_{{}^\ast\!\!\:\Gamma}$,
respectively. Using a pair of dual bases we have
$
\coev(1_A)=\sts^i\otimes_A s_{i\,}\in  \stG\otimes_A\Gamma$,
$\coev'(1_A)=\oR^\al\ra s_i\otimes_A \oR_\al\ra \sts^i\in
\Gamma\oA\stG$ so that these conditions read, for all $s\in\Gamma$, $\sts\in\stG$,
\begin{equation}\label{sstsslr}\oR^\al \ra s_i\le\oR^\beta\ra s,
\oR_\beta\oR_\al\ra\sts^i\re=s~~,~~~\le\oR^\be\oR^\al \ra
s_i,\oR_\beta\ra\sts\re\oR_\al\ra\sts^i=\sts~.
\end{equation}
We prove the first equation 
\begin{equation}\label{proofdualitylr}
\begin{split}
  \oR^\al \ra s_i\le\oR^\beta\ra s,\oR_\beta\oR_\al\ra\sts^i\re&=
  \oR^\gamma\ra\le\oR^\beta\ra s,\oR_\beta\oR_\al\ra\sts^i\re\;\oR_\gamma \oR^\al
  \ra s_i\\[.2em]
  &=
\le \oR^\gamma\oR^\beta\ra
s,\oR^\delta\oR_\beta\oR_\al\ra\sts^i\re\oR_\delta\oR_\gamma \oR^\al \ra
s_i\\[.2em]
&=
\le \oR^\gamma\oR^\beta\ra
s,\oR^\delta\oR_\al\oR_\gamma\ra\sts^i\re\oR_\delta\oR^\al \oR_\beta \ra
s_i\\[.2em]
&=\le \oR^\beta\ra
s,{\oR_\beta}{}^{}_{\!\!\:(1)}\ra\sts^i\re{\oR_\beta}{}^{}_{\!\!\:(2)}\ra s_i\\[.2em]
&=\le
s,\sts^i\re s_i\\[.2em]
&=s
\end{split}
\end{equation}
where in the second line we used $(\Delta \otimes
\id)\RR^{-1}=\RR^{-1}_{23}\RR^{-1}_{13}$, in the  third line  the Yang--Baxter equation
$\RR^{-1}_{12}\RR^{-1}_{13}\RR^{-1}_{23}=\RR^{-1}_{23}\RR^{-1}_{13}\RR^{-1}_{12}$,
in the fourth that $\oR^\delta\oR_\al\otimes\oR_\delta\oR^\al=1\otimes
1$ due to  triangularity of $\RR$ and then $(\id\otimes
\Delta)\RR^{-1}=\RR^{-1}_{12}\RR^{-1}_{13}$. In the fifth we used that $\coev$
is $H$-equivariant so that for all $h\in H$, $h\ra (\sts^i\otimes
s_i)=h\ra (\coev(1_A))=\coev(h\ra 1_A)=\varepsilon(h) \sts^i\otimes
s_i$ and in the last line the hypothesis that  ${}^\ast \Gamma$ is right dual to $\Gamma$.
The proof of the second equation in \eqref{sstsslr} is similarly obtained.
    \end{proof}
Vice versa, if $\Gamma$ is left rigid  with left dual $\stG$, evaluation and
coevaluation maps $\le~,~\re'$ and $\coev'$, then $\stG$ is
also right dual to $\Gamma$ with evaluation and coevaluation maps
$\le~,~\re$ and $\coev$ implicitly defined by \eqref{tevcoev}, i.e.,
$
\le~,~\re:=\langle~,~\rangle'\circ \tau^{}_{\Gamma,\stG}$,
$\coev:=
\tau^{}_{\Gamma,\stG}\circ\coev'
$; (just read \eqref{proofdualitylr} from the fifth line to the first
and then use \eqref{sstsslr}).
We therefore speak of finitely generated and projective $A$-modules in
$\cat$, there is no need to specify if they are left or right $A$-modules.

  We denote by $\catfp$ the full subcategory in $\cat$ of finitely generated and
  projective modules. 
  Every module has a left and a right dual and therefore the  category
  $\catfp$ is a rigid category. Rigid symmetric monoidal categories are called compact closed
  categories. We have shown that the  symmetric monoidal
  category $\catfp$ is a compact closed category.

\begin{remark}\label{ep'}The exact pairing $\le~,~\re':
\stG\oA\Gamma\to A$ induces the isomorphisms in $\cat$ (cf. Proposition \ref{propo:projhomiso})
\begin{equation}\label{jpip}
\begin{split}
\le~\,,\,~~\re'&: \Gamma\oA W \longrightarrow {}_A\hom(\stG, W)~,~~
s\otimes_A w\mapsto \le ~\,, s\otimes_A w\re':=\le ~\,, s\re' \otimes_A w~, 
\\[.4em]
\le~~~,~\re'&: W\oA \Gamma \longrightarrow \hom_A(\Gamma^*, W)~,~~w\otimes_A 
s\mapsto \le w\otimes_A s, ~\, \re':=w\oA\le s,\, ~\re' ~.
\end{split}
\end{equation}
\end{remark}

For later use we notice that the
dual of the braiding isomorphism equals the braiding on the dual
modules: In the
notations of Remark \ref{newevcoev} (and recalling the definition of
right dual morphism given after Remark \ref{rem:gmc}),
\begin{equation}\label{taudual}
{}^{*\!}(\tau^{}_{\Gamma,\Upsilon})=\tau^{}_{\stG,\stU}~,
\end{equation} 
that is, for all $s\in \Gamma,
u\in\Upsilon, \sts\in \stG, \stu\in\stU$,
$\le s\oA u, \tau^{}_{\stG,\stU}(\sts\oA\stu)\re
=\le\tau^{}_{\Gamma,\Upsilon}(s\oA u),\sts\oA\stu\re$.
This is equivalent to
$\le\tau^{-1}_{\Gamma,\Upsilon}(u\otimes_A s),
\tau^{}_{\stG,\stU}(\sts\otimes_A\stu)\re=\le u\oA s, \sts\oA\stu\re$
and  is easily proven,
\begin{equation*}
  \begin{split}
\le\tau^{-1}_{\Gamma,\Upsilon}(u\otimes_A s),\tau^{}_{\stG,\stU}(\sts\otimes_A\stu)\re
&=\le R_\al\ra s\,,\le R^\al \ra u,\oR^\beta \ra \stu\re
\oR_\beta\ra\sts \re\\[.2em]
&=\oR^\gamma\ra\le R^\al \ra u,\oR^\beta \ra \stu\re\, \le \oR_\gamma R_\al\ra s,
\oR_\beta\ra\sts \re\\[.2em]
&=\le\oR^\gamma R^\al \ra u,\oR^\delta\oR^\beta \ra \stu\re \le\oR_\delta \oR_\gamma R_\al\ra s,
\oR_\beta\ra\sts \re\\[.2em]
&=\le u,\oR^\delta\oR^\beta \ra \stu\re \le\oR_\delta \ra s,
\oR_\beta\ra\sts \re\\[.2em]
&=\le u,\oR^\beta \ra\stu\re\,\oR_\beta\ra\le s,\sts \re\\[.2em]
&=\le u\otimes_A s,\sts\otimes_A\stu\re~.
\end{split}
\end{equation*}

\subsection{Examples}
Let $H$ be a triangular Hopf algebra and $A$ a braided commutative $H$
module algebra. We study examples of compact closed categories $\catfp$ of symmetric relative $H$-equivariant $A$-bimodules finitely generated and projective as $A$-modules.
The first example arises when $A=K$ is a cotriangular Hopf algebra and
$H=U^\op\otimes U$ is the triangular Hopf algebra obtained from the
triangular Hopf algebra $U$ dual to $K$. 
Another example is that of equivariant vector
bundles on a manifold and a further one (adapting the treatment in \cite[\S 6]{BSS1}
to the compact closed category context)  is obtained via 
noncommutative Drinfeld twist deformation of equivariant vector bundles.
\begin{example}\label{cotriangularK}
 {\it{Cotriangular Hopf Algebra.}} Let $A=K$ be a
   finite dimensional
  Hopf algebra
  over a field  $\Bbbk$
  and  let $U$
  be the
dual Hopf algebra. Right (left) $K$-coactions correspond to left
(right) $U$-actions on modules (using Sweedler like notation, given a
right $K$-coaction $V\to V\otimes K, v\mapsto v_0\otimes v_1 $ we have the left
$U$-action $\ra: U\otimes V\to V, \xi\ra v=v_0\!\: \xi(v_1)$, while 
given a left $K$-coaction $V\to K\otimes V, v\mapsto v_{-1}\otimes v_0 $ we have the right
$U$-action $\triangleleft: V\otimes U\to V, v\triangleleft \xi=
\xi(v_{-1})\;\! v_0$). Vice versa, since $K$ is finite dimensional over the field
$\Bbbk$, $\Hom_\Bbbk(U,V)\simeq V\otimes K$  (cf. \eqref{jmathdef}); this
implies that given a left
$U$-action,  the map $V\to\Hom_\Bbbk(U,V)$, $v\mapsto \Delta_v$; $\Delta_v(\xi)=\xi\ra
v$ defines a right $K$-coaction $\Delta_\Rs: V\to V\otimes K$; similarly,
right $U$-actions define left $K$-coactions. Moreover, right
(left) $K$-comodule algebras are equivalently left (right)
$U$-module algebras.
In particular, since  $K$ is a $K$-bicomodule
algebra via the coproduct, then it is a $U$-bimodule algebra. 
Recall that a $K$-bicovariant bimodule \cite[Def. 2.3]{Wor} is a
$K$-bimodule with compatible and commuting left and right $K$-coactions.
The duality between
$K$-coactions and $U$-actions
implies that this category is equivalent to that of relative $U$-bimodules
$K$-bimodules (these are relative left $U$-modules $K$-bimodules and
relative right $U$-modules $K$-bimodules with commuting left and right
$U$-actions). Now $U$-bimodules ($U$-bimodule algebras) are equivalently left
${U}^\op\otimes U\,$-modules (left  ${U}^\op\otimes U\,$-module algebras), where $U^\op$ is the Hopf algebra with opposite
product and ${U}^\op\otimes U$ is the tensor product Hopf algebra; similarly,
relative $U$-bimodules $K$-bimodules are equivalently relative left
${U}^\op\otimes U$-modules $K$-bimodules (the
${U}^\op\!\otimes \!U\,$-action on a $K$-bicovariant bimodule $V$ reads
$(\zeta\otimes \xi)\ra v=   \zeta(v_{-1})v_0\xi(v_1)$ and on
$K$ itself $(\zeta\otimes \xi)\ra k=    \zeta(k_{1})k_2\xi(k_3)$).
Hence the monoidal category of $K$-bicovariant bimodules 
is equivalent to the monoidal category 
${}^{U^\op\!\:\otimes U}_{}{}^{}_{\!\!\!\!\!\!\!K\!}\MMM^{}_K$ of relative
${U^\op\otimes U}\,$-modules $K$-bimodules. These are free $K$-modules
since $K$-bicovariant bimodules are free $K$-modules
(cf. \cite[Thm. 2.1]{Wor}).

Let now $K$ be cotriangular, this is the case if and only if its dual $U$ is
triangular. 
Let 
$\mathscr{R}$ be the triangular structure of
$U$, i.e., the cotriangular structure of $K$; in particular we have the quasi-commutativity property, for all 
$k,k'\in K$,
\begin{equation}\label{cotriangprop}
  k'k=\mathscr{R}(k^{}_{(1)}\otimes k'_{(1)}) k^{}_{(2)}k'_{(2)}
  \mathscr{R}^{-1}(k^{}_{(3)}\otimes k'_{(3)})~.
  \end{equation}
Furthermore $U^\op$ is triangular with  $\RR$-matrix $\mathscr{R}^{-1}$ and
$U^\op\otimes U$ is also triangular with $\RR$-matrix
${\RR}=(\id\otimes \rm{flip}\otimes\id)(\mathscr{R}^{-1}\otimes
\mathscr{R})$. We show that the quasi-commutativity property
\eqref{cotriangprop} of $K$ is just
the braided commutativity property
of $K$ with respect to the $\RR$-matrix of  $U^\op\otimes U$, for all
$k,k'\in K$,
\begin{equation}\label{cotquasi}
  \begin{split}
    k'k=\oR^\al\ra k\,\oR_\al\ra k'&=k^{}_{(2)}k'_{(2)}\,
    \oR^\al(k^{}_{(1)}\otimes k^{}_{(3)})\oR_\al(k'_{(1)}\otimes k'_{(3)})\\[.2em]
 &=
   k^{}_{(2)}k'_{(2)}{\mathscr{R}}^\al(k^{}_{(1)})\overline{\mathscr{R}}^\beta(k^{}_{(3)})
{\mathscr{R}}_\al(k'_{(1)})    \overline{\mathscr{R}}_\beta(k'_{(3)})\\[.2em]
    &=      {\mathscr{R}}(k^{}_{(1)}\otimes k'_{(1)}) k^{}_{(2)}k'_{(2)}
  \mathscr{R}^{-1}(k^{}_{(3)}\otimes k'_{(3)})~,
    \end{split}
  \end{equation}
where we used the notation $\RR^{-1}=\oR^\al\otimes\oR_\al=(\id\otimes \rm{flip}\otimes\id)(\mathscr{R}\otimes
\mathscr{R}^{-1})={\mathscr{R}}^\al\otimes\overline{\mathscr{R}}^\beta
\otimes {\mathscr{R}}_\al\otimes    \overline{\mathscr{R}}_\beta$.

In conclusion,  if $K$ is a finite dimensional
cotriangular Hopf algebra over a field then  $K$ is a braided commutative
algebra with respect to the triangular Hopf algebra $U^\op\otimes U$
and, recalling also the previous section,
the category ${}^{{U^\op\!\:\otimes
    U}}_{}{}^{}_{\!\!\!\!\!\!\!K\!}\MMM^{\sym,{\rm{f}}}_K$
of finitely generated
 symmetric relative ${U^\op\otimes
    U}$-modules $K$-modules is a
rigid symmetric category (compact closed category) 
of free  $K$-modules.

Similarly, we can consider $K$ a cotriangular quantum group, for
example of the $A,B,C,D$ series given via multiparametric $\RR$
matrices, with dually paired triangular topological Hopf algebra
$U$ over $\mathbb{C}[[\hbar]]$ \cite{res}. 
As before, $K$ is braided commutative and the category
${}^{{U^\op\!\:\otimes U}}_{}{}^{}_{\!\!\!\!\!\!\!K\!}\MMM^{\sym,\rm{fp}}_K$ of
finitely generated and projective symmetric  ${U^\op\otimes
    U}$-equivariant $K$-modules is a
compact closed category. 
\endex\end{example}

Compact closed categories are easily obtained via Drinfeld twists
of compact closed categories. Recall that a Drinfeld twist is an
invertible element
$\FF\in H\otimes H$ satisfying the cocycle and normalization properties:
\begin{equation}\label{twist}
(\FF\otimes
\id)(\Delta\otimes\id)\FF=(\id\otimes\FF)(\id\otimes\Delta)\FF~,~~(\varepsilon\otimes\id)
\FF=(\id\otimes\varepsilon)\FF=1_H~.
\end{equation}
Given a twist we deform the triangular Hopf algebra
$H$ with universal $\RR$-matrix $\RR$ in the triangular Hopf algebra $H^\FF$
that as algebra is the same as $H$, while it has coproduct $\Delta^\FF(\hxi)=\FF\Delta(\hxi)\FF^{-1}$ and
universal $\RR$-matrix $\RR^\FF=\FF_{21}\RR\FF^{-1}$. We deform the
braided commutative algebra
$A$ in the braided commutative algebra $A^\FF$ that equals $A$ as
$\Bbbk$-module, while it has
multiplication,  for all $a,b\in A$,
$a\cdot^{{_{^\FF}}} b:=(\of^\al\ra a) (\of_\al\ra b)$, where $\FF^{-1}=\of^\al\otimes\of_\al$. The
normalization conditions imply that the unit of
$A$ is also the unit of $A^\FF$.
We further associate to any module $\Gamma$ in $\catfp$ a new module 
$\Gamma^\FF$ in $\catFF$ that has the same $\Bbbk$-module structure.
The $H^\FF$-module structure is given by the $H$-action (indeed as 
algebras $H^\FF$ equals $H$) and the left and right $A^\FF$-actions are 
given by, for all $a\in A^\FF$, $s\in \Gamma$,
$a\dotF s:=(\of^\al\ra a)(\of_\al \ra s)
$
and $s
\dotF a=(\of^\al\ra s)(\of_\al \ra a)$.
Due to these assignments the compact closed category $\catfp$ is 
equivalent to the compact closed category $\catfpFF$. (They are
equivalent as symmetric monoidal categories and hence, since any monoidal
functor preserves rigidity of objects, as compact closed categories). 

 \begin{example}\label{ex:bundle}{\it $U\gg$-equivariant
  $C^\infty(M)$-bimodules from $G$-equivariant vector bundles on $M$.}
Let $G$ be a Lie group and $M$ a $G$-manifold with right $G$-action
$r: M\times G\to M\,,~(m,g)\mapsto m\:\!g$.  We assume manifolds to be
finite dimensional and second countable. In this case every finite
rank vector bundle has a finite open covering of the base trivializing the bundle.

Recall that a $G$-equivariant vector bundle $E\to M$ is a vector
bundle $E\to M$ where $E$ and $M$ are right $G$-manifolds, the
$G$-actions on $E$ and $M$ are compatible with the projection $E\to M$ and $G$ acts via linear isomorphisms on the fibers of $E$.

The right $G$-action on $M$ pulls back to a left $G$-action 
$G\times C^\infty(M)\to C^\infty(M)$ on the algebra $C^\infty(M)$ 
of smooth complex valued functions on $M$.
Using the $G$-action on $E$  we also have a left $G$-action on the
$C^\infty(M)$-module of sections $\Gamma(E)$ (for all $m\in M,  s\in
\Gamma(E)$, $(g\ra s)(m)=s(mg)g^{-1}$). With this action $\Gamma(E)$
is a relative $G$-module $C^\infty(M)$-module (or $G$-equivariant
$C^\infty(M)$-module), i.e., for all $a\in
C^\infty(M), s\in  \Gamma(E) , g\in G, g\ra as=(g\ra a) (g\ra s)$. 
Since $E\to M$ is trivialized by a finite open covering of $M$ then $\Gamma(E)$ is finitely generated and projective as a $C^\infty(M)$-module.
This construction leads to a functor $\Gamma: G$-{\sl VecBun}$_M\to
{}^G_{C^\infty(M)}\MMM^{\sym,\fp}_{C^\infty(M)} $ from the category of
$G$-equivariant vector bundles over $M$ to that of 
$G$-equivariant $C^\infty(M)$-modules that are
finitely generated and
projective over $C^\infty(M)$. (If $G$ is the trivial group the
functor is an equivalence proving Serre-Swan theorem in the smooth
context).
Both categories are rigid symmetric monoidal, hence are compact closed
categories and the functor is symmetric monoidal.

Associated with the $G$-action on $M$ we have a (smooth) action of the
Lie algebra $\gg=$Lie$\,G$ on $A=C^\infty(M)$ via derivations.
It canonically lifts to an action of the universal enveloping
algebra $U\gg$ on $A$ that therefore  becomes a
$U\gg$-module algebra. Similarly, $\Gamma(E)$ is a relative
$U\gg$-module $A$-bimodule.  (Technically, $C^\infty(M)$ and the
finitely generated projective module $\Gamma(E)$ are nuclear
Fr{\'e}chet spaces with respect to the usual
$C^\infty$-topology. Furthermore, $C^\infty(M)$ is a unital Fr{\'e}chet algebra
with the usual pointwise product $\mu:= \mathrm{diag}_M^\ast :
C^\infty(M)\,\widehat\otimes\, C^\infty(M) \to C^\infty(M) $ given by pull-back of
the diagonal map $\mathrm{diag}_M : M\to M\times M$.
Here $C^\infty(M)\,\widehat\otimes\, C^\infty(M)\simeq C^\infty(M\times M)$
denotes the completed tensor product. Therefore $C^\infty(M)$ is an algebra object in the category of $U\gg$-Fr{\'e}chet spaces).

It follows that we have a symmetric monoidal functor {\sl ULie}$\,\circ \Gamma$ from the
category $G$-{\sl VecBun}$_M$
to the category of $U\gg$-equivariant $A$-bimodules. The functor
is valued in the sub-category
 ${}^{U\gg}_A\MMM_{A\,\geom}^{\sym,\fp}$, where objects are
the modules of sections $\Gamma(E)$ of $G$-equivariant vector bundles
$E\to M$ and morphisms are $U\gg$-equivariant $A$-bimodule morphisms.  
This is a compact closed category.
\endex\end{example}

\begin{example}\label{ex:bundleformal} {\it Formal braided equivariant vector bundles.}
Associated with the ring $\mathbb{C}[[\hbar]]$ of formal power
series in $\hbar$ with coefficients in $\mathbb{C}$ we have the
 formal power series extension
of the $\mathbb{C}$-modules $U\gg$, $A$ and $\Gamma(E)$, denoted as
usual $U\gg[[\hbar]]$, $A[[\hbar]]$ and $\Gamma(E)[[\hbar]]$.
The tensor product in the definition of the
${\mathbb{C}}[[\hbar]]$-Hopf algebra  $U\gg[[\hbar]]$ is the
completed one in the $\hbar$-adic topology (cf. e.g. \cite[\S XVI]{Kassel}), so that 
 $U\gg[[\hbar]]\,\widehat\otimes\, U\gg[[\hbar]]\simeq (U\gg\otimes
U\gg)[[\hbar]]$ (and the tensor product in the definition of the
$\mathbb{C}[[\hbar]]$-algebra $A[[\hbar]]$ is the 
completed one in the Fr\'echet and $\hbar$-adic topologies, so that 
$A[[\hbar]]\,\widehat \otimes\, A[[\hbar]]\simeq C^\infty(M\times M)[[\hbar]]$). 
The $A[[\hbar]]$-bimodule $\Gamma(E)[[\hbar]]$ is topologically  finitely
generated and projective and we denote by
${}^{U\gg[[\hbar]]}_{A[[\hbar]]}\MMM_{A[[\hbar]]\,\geom}^{\sym,\fp}$ the associated
compact closed category.
It has been obtained via the change of base ring $\mathbb{C}\to\mathbb{C}[[\hbar]]$
that induces the braided monoidal functor $\mathbb{C}[[\hbar]]\widehat\otimes:
{}^{U\gg}_A\MMM_{A\,\geom}^{\sym,\fp}
\to{}^{U\gg[[\hbar]]}_{A[[\hbar]]}\MMM_{A[[\hbar]]\,\geom}^{\sym,
  \fp}$ that to
$\Gamma(E)$ associates the extension
$\mathbb{C}[[\hbar]]\!\;\widehat\otimes\!\;\Gamma(E)\simeq \Gamma(E)[[\hbar]]$
and extends morphisms by $\mathbb{C}[[\hbar]]$-linearity.

By considering a twist $\FF$ of $U\gg[[\hbar]]$  we obtain the rigid braided commutative 
monoidal category ${}^{U\gg[[\hbar]]^\FF}_{A[[\hbar]]^\FF}\!\MMM_{A[[\hbar]]^\FF\,\geom}^{\sym,
      \fp}$. This provides a deformation quantization of $G$-equivariant
    bundles on $A=C^\infty(M)$.
\endex\end{example}

A further example, as shown in the next section, is provided by the
symmetric monoidal category of the modules of
covariant and contravariant tensor fields canonically associated with
the differential calculus on a braided commutative $H$-module algebra $A$.

\section{Differential and Cartan Calculus}\label{sectCC}
We give a self contained  thorough exposition
of the differential and Cartan calculus on a braided commutative
$H$-module algebra $\AA$ (cf. \eqref{bsalgA}), where
$H$ is a triangular Hopf algebra. 
For later use, details on the duality between covariant and contravariant vector
fields are provided. The presentation is constructive and when 
$A$ is a cotriangular Hopf algebra we obtain
a bicovariant differential calculus \`a  la Woronowicz \cite{Wor}. We then show its equivalence with that in \cite{Gomez-Majid}.
The notations and conventions adopted stem
from Section 2 and will be used throughout the paper.
They differ from \cite{Gurevich}, where the construction of a K\"ahler
differential and Cartan calculus of an algebra in a symmetric
rigid category was outlined, and from \cite{Weber} where, as here, a braided
derivations approach is pursued. Apart from the
example and the notions of Lie algebra and bimodule of
braided derivations (found in the beginning)  the expert reader can
skip this section and go back to relevant
formulae and theorems when they will be referred to in  later sections.

\subsection{Braided derivations and differential calculus}
Let $\Der_\RR(\AA)$ be the $\Bbbk$-module of
braided derivations. These are $\Bbbk$-linear maps $u\in
\hom_\Bbbk(\AA,\AA)$ that satisfy the braided Leibniz rule, for all $a,b\in A$,
\begin{equation}\label{braidedleibniz}
u(ab)=u(a)b\,+\,\oR^\alpha\ra a \,(\oR_\alpha \RA u)(b)~.
\end{equation}
$\Der_\RR(\AA)$  is an $H$-submodule of
$\hom_\Bbbk(\AA,\AA)$. Indeed, for all $h\in H, u\in \Der_\RR(\AA)$ and
$a,b\in A$
\begin{align*}
  (h\ra u)  (ab)&=h_{(1)}\ra \big(u(S(h_{(3)}) \ra a\,S(h_{(2)})\ra b)\big)\\
&=h_{(1)}\ra \big(u(S(h_{(3)}) \ra a)\,S(h_{(2)})\ra b\big)+h_{(1)}\ra
                                                                        \big((\oR^\al
                                                                        S(h_{(3)})
                                                                        \ra
                                                                        a)(\oR_\al\ra u)(S(h_{(2)})\ra
                                                                        b)\big)\\
 &=h_{(1)}\ra (u(S(h_{(2)}) \ra a))\, b+\oR^\al h_{(2)}S(h_{(3)})\ra 
                                                                        a\,(\oR_\al
    h_{(1)}\ra u)(b)\\
 &=(h\ra u)(a)\, b + \oR^\al \ra 
                                                                        a\,(\oR_\al\ra
    (h\ra u))(b)
\end{align*}
where in the third line we used \eqref{adjact} in the form $h\ra
L(v)=(h_{(1)}\ra L)(h_{(2)}\ra v)$ (with $L$ given by $\oR_\al\ra u$,
$h$ given by $h_{(1)}$ and
$v$ given by $S(h_{(2)})\ra b$) and then that 
$\Delta(h)\RR^{-1}=\RR^{-1}\Delta^{cop}(h)$.

As for derivations on a commutative algebra, it is not difficult to see
(cf. \cite[Lemma 3.1]{Weber}) that the braided commutator
$$[~,~]: \Der_\RR(\AA)\otimes \Der_\RR(\AA)\to \Der_\RR(\AA)~,~u\otimes v\mapsto
~[u,v]:=u v-{\oR}^\al\ra v\,{\oR}_\al\ra u~,$$
(where composition of operators is understood) closes in
$\Der_\RR(\AA)$, is an $H$-equivariant map (for all $h
\in H, u,v, \in  \Der_\RR(\AA)$, $h\ra [u,v]=[h_{(1)}\ra u,h_{(2)}\ra v]$)
and structures the $H$-module $\Der_\RR(\AA)$ as an $(H, \RR)$-Lie
algebra, i.e., we have the braided antisymmetry property and the  Jacobi identity, for all $u,v,z\in \Der_\RR(\AA)$,
\begin{equation}\label{brLie}
[u,v]=-[{\oR}^\al \ra v,{\oR}_\al \ra u]~,~~
[u,[v,z]]=[[u,v],z]+[{\oR}^\al \ra v,[{\oR}_\al\ra u,z]]~.
\end{equation}
Braided derivations are furthermore a module in ${}^{H}_\AA\MMM_\AA^\sym$
by defining, for all $a,b\in A$, 
\begin{equation}\label{aub}
(au)(b)=a\,u(b)~~,~~~ua=(\oR^\al\ra a)\oR_\al\ra u~;
\end{equation}
$a u$ is again a braided derivation because of the braided commutative
property \eqref{bsalgA} of $\AA$, for all $a, b, c \in A, u\in \Der_\RR(A)$,
\begin{align*}
(au)(bc)&=a\,u(bc)\\&=a\,u(b)c\,+a\,\oR^\alpha\ra b \,(\oR_\alpha \RA
u)(c)\\ &=
(au)(b)c\,+(\oR^\beta\oR^\alpha\ra b) (\oR_\beta\ra a)\,(\oR_\alpha
          \RA u)(c)\\
&=  (au)(b)c\,+(\oR^\alpha\ra b) (\oR_\alpha\ra (a u))(c)~,
\end{align*}
where in the last passage  we used the quasitriangularity property
$\RR^{-1}_{12}\RR^{-1}_{13}=(\id\otimes \Delta)\RR^{-1}$. One
 similarly proves the compatibility $[u,av]=u(a)v + {\oR}^\al \ra
a[{\oR}_\al \ra u, v]$.

From now on we set $$\Vect(A):=\Der_\RR(A)$$ and refer to it as the bimodule
of (braided) vector fields.
\\

We denote by $\underline{\Omega}(A):={}^*\Vect(A)={}_\AA\hom(\Vect(A),
A)\subset {}_\Bbbk\hom(\Vect(A),
A)$ the right dual module of left $A$-linear maps $\Vect(A)\to  A$
with $H$-action $\ra^\cop$ defined in \eqref{copadjact} and
$A$-bimodule structure defined in \eqref{AhomAA}. It  is a module in $\cat$
and as in \eqref{ppair} we denote the evaluation of an element in
$\underline{\Omega}(A)$ on a vector field via the bracket
\begin{equation}
\label{evaluation<>}
    \langle~,~\rangle :\Vect(A) \otimes_\AA \underline{\Omega}(\AA)\to
    \AA~,~~u\otimes_A\omega\mapsto \le u,\omega\re
  \end{equation}
that is $H$-equivariant (cf.\eqref{Hppair}).

We define the map  $\dd: A\to \underline{\Omega}(A)$  by 
\begin{equation}\label{defofd}
  \langle u, \dd a\rangle=u(a)~,
  \end{equation}
 for all $u\in \Vect(A)$. This definition is well-defined since  both  $\le\,~,\dd a\re:\Vect(A)\to A$
and $\hat a: \Vect(A)\to A$, $u
\mapsto u(a)$ are left $A$-linear maps. The map
$\dd$ is $H$-equivariant, indeed, for all
$\hxi\in H$ the identities $$\hxi\ra\le u,\dd a\re= \le\hxi_{(1)}\ra u,\hxi_{(2)}\ra^\cop \dd
a\re~,~~\hxi\ra (u(a))=(\hxi_{(1)}\ra
u)(\hxi_{(2)}\ra a)
$$
imply $\hxi\ra^\cop (\dd\:\!a)=\dd(\hxi\ra a)$.
Next we prove the undeformed Leibniz rule $\dd(ab)=(\dd a)b+a\dd b$:
\begin{align*}\le u,\dd (ab)\re &=u(ab)=u(a)b+\oR^\al\ra a\,(\oR_\al \ra
  u)(b)=\le u,(\dd a)b\re+(\oR^\al\ra a)\le\oR_\al\ra u,\dd b\re
  \\
                                &=\le u,(\dd a)b\re+\le u , a\dd b\re~,
\end{align*}
where we used \eqref{aub} and
that the pairing $\le~,~\re$ is well-defined on the balanced tensor
product $\VV\otimes_A\underline{\Omega}(A)$. 

The module of one-forms $\Omega(A)$ is the submodule of
$\underline{\Omega}(A)$ in $\cat$ defined by
\begin{equation}\label{defofOmega}
  \Omega(A):=A\dd A=\{\omega\in\underline{\Omega}(A)~;~\omega=a^i \dd a_i\}
\end{equation}
for all $a^i,a_i\in A$, with finite sum over the index $i$ understood
(the right $A$-action closes in $\Omega(A)$ due to the Leibniz rule).  
A {\it first order differential
calculus on an algebra $A$} is an $A$-bimodule $\Omega(A)$ and a map $d:A\to \Omega(A)$ 
that satisfies the Leibniz rule and such that  every element of
$\Omega(A)$ is of the form $a^i\dd a_i$. We thus see that
$(\Omega(A),\dd)$ as defined in \eqref{defofd} and \eqref{defofOmega}
is a first order differential calculus on $A$. It is an
$H$-equivariant differential calculus since $A$ and $\Omega(A)$ are
in $\cat$ and $\dd: A\to \Omega(A)$ is $H$-equivariant.

 The contraction operator is the morphism in $\cat$ defined by
\begin{equation}\label{contractionoperator}
\ii: \Vect(A)\to \Omega(A)^\ast:=\hom_A(\Omega(A),A) \:,~u\mapsto \ii_u~,~~\ii_u( (\dd a^i)a_i):=\langle
u, (\dd a^i)a_i\rangle=u(a^i)a_i~.
\end{equation}

\begin{proposition}\label{Om*V} The contraction $\ii: \Vect(A)\to
  \Omega(A)^\ast:=\hom_A(\Omega(A),A)$ is an isomorphism in $\cat$.
\end{proposition}
\begin{proof}
The $\Bbbk$-linear map $\vv:\hom_A(\Omega(A), A)\to \hom_\Bbbk(A, A)$,
$\psi\mapsto \vv_\psi$, given by, for all $a\in A$,
$\vv_\psi(a):=\psi(\dd a)$ is $H$-equivariant, for all $h\in H$, $a\in A$,
$$
h\ra(\vv_\psi(a))=h\ra(\psi(\dd a))=(h_{(1)}\ra\psi)(h_{(2)}\ra^\cop\dd
a)
=(h_{(1)}\ra\psi)(\dd (h_{(2)}\ra a))=\vv_{h_{(1)}\ra\psi}(h_{(2)}\ra a)~.
$$
We show
$\vv(\hom_A(\Omega(A), A))\subseteq\Vect(A)$, indeed,  for all
$\psi\in \hom_A(\Omega(A), A)$, $\vv_\psi$ is a braided
derivation: for all $a,b\in
A$,
\begin{equation}
  \begin{split}
\vv_\psi(ab)=\psi((\dd a)b+a\dd b)&=\psi(\dd a)b+
\psi((\oR^\al\ra^\cop\dd b)\oR_\al\ra
a)\\[.2em]
&=\vv_\psi(a)b+\psi(\oR^\al\ra^\cop\dd b)\oR_\al\ra
a\\[.2em]
&=\vv_\psi(a)b+(\oR^\beta\oR_\al\ra
a)\:\!\oR_\beta\ra(\psi(\oR^\al\ra^\cop\dd b))\\[.2em]
&=
\vv_\psi(a)b+(\oR^\beta\oR^\gamma\oR_\al\ra
a)\!\:(\oR_\beta\ra\psi)(\oR_\gamma\oR^\al\ra^\cop\dd b)\\[.2em]
&=\vv_\psi(a)b+(\oR^\beta\ra
a)(\oR_\beta\ra\psi)(\dd b)\\[.2em]
&=
\vv_\psi(a)b+(\oR^\beta\ra
a)\vv_{\oR_\beta\ra\psi}(b)\\[.2em]
&=\vv_\psi(a)b+(\oR^\beta\ra
a)(\oR_\beta\ra\vv_{\psi})(b)~,
\end{split}
\end{equation}
where  in the fourth line we used 
$
    (\mathrm{id}\otimes\Delta)\mathcal{R}^{-1}
    =\mathcal{R}^{-1}_{12}\mathcal{R}^{-1}_{13}
    $ and in the fifth triangularity of $\RR$.
    
The induced map $\vv: \hom_A(\Omega(A), A))\to\Vect(A)$ is the inverse
of the contraction operator: The equality $\vv\circ
\ii=\id_{\Vect(A)}$ is easily proven, for all $u\in \Vect(A)$, $a\in
A$, $\vv_{\ii_u}a=\ii_u(\dd a)=u(a)$. We are left to prove $\ii\circ
\vv=\id_{\Omega^\ast(A)}$; for all $\psi\in {\Omega^\ast(A)}$, $(\dd
a^i)a_i\in \Omega(A)$, $\ii_{\vv_\psi}((\dd
  a^i)a_i)=\vv_\psi(a^i)a_i=\psi(\dd a^i)a_i=\psi((\dd a^i)a_i)$.
\end{proof}
\begin{corollary}\label{pnondeg}
The pairing $\le~,~\re: \Vect(A)\oA\Omega(A)\to A$ is non-degenerate:
i) $\le u,\omega\re=0$ for all $\omega\in \Omega(A)$ implies  $u=0$; ii)  $\le u,\omega\re=0$  for all $u
\in \Vect(A)$ implies $\omega=0$.
\end{corollary}
\begin{proof}
i) for all $a\in A$, $0=\le u, \dd a\re =u(a)$ implies
$u=0$. ii) Follows from $\Vect(A)\simeq \Omega(A)^*$.
\end{proof}
\begin{corollary}\label{omega=uomega}
If $\Omega(A)$ is finitely generated and projective over $A$ then
$\underline{\Omega}(A)= \Omega(A)$.
\end{corollary}
\begin{proof}
If $\Omega(A)$ is in $\catfp$
so
is $\Vect(A)\simeq \Omega(A)^*$ and 
the canonical isomorphism $^*(\Omega(A)^*)\simeq \Omega(A)$
(cf. paragraph before Proposition  \ref{propo:fgpmodprop}) implies
$\underline{\Omega}(A):={}^*\Vect(A)\simeq \Omega(A)$; henceforth, since  $\Omega(A)\subseteq \underline{\Omega}(A)$, 
$\underline{\Omega}(A)= \Omega(A)$.
\end{proof}

Associated with the modules $\Omega(\AA)$ and $\Vect(A)$ in
${}^{H}_\AA\MMM_\AA^{\sym}$ we have the modules
$\Tau^{p,0}=\Omega(A)^{\otimes_A p}$
and $\Tau^{0,q}=\Vect(A)^{\otimes_A q}$, $p,q\in
\mathbb{N}$, with $\Tau^{0,0}=A$,  and the
graded $H$-module algebras 
of contravariant tensor fields $\Tau^{\bullet,0}=\bigoplus_{p\in\mathbb{N}} \Tau^{p,0}$ and of covariant tensor
fields  $\Tau^{0,\bullet}=\bigoplus_{q\in \mathbb{N}} \Tau^{0,q}$. We also have
the graded $H$-module algebra $\Tau^{\bullet,\bullet}=\bigoplus_{p,q\in\mathbb{N}}
\Tau^{p,q}$ with product that on elements of homogeneous degree is defined by
\begin{equation}\label{minT} 
\otimes_A : \Tau^{p,q}\otimes \Tau^{p',q'}\to
\Tau^{p+p',q+q'}~,~~\theta\otimes_A\nu\otimes\theta'\otimes_A\nu'\mapsto
\theta\otimes_A\oR^\al \ra^\cop \theta'\otimes_A
\oR_\al\ra\nu\otimes_A\nu'~,
\end{equation}
where $\theta\in \Tau^{p,0}$, $\nu\in \Tau^{0,q}$,  $\theta'\in
\Tau^{p',0}$, $\nu'\in \Tau^{0,q'}$.
\\

We define the module $\Omega^{ r}(A)$ of $r$-forms as the submodule of
$\TT^{r,0}$ of completely braided antisymmetric tensors. This is obtained
by writing every permutation $\wp$ of $r$ elements as composition of elementary
nearest neighbour transpositions $\wp=t_{k_1}\circ t_{k_2}\ldots \circ t_{k_n-1}$ (where
$t_k:(1,2,\ldots k,k+1,\ldots n)\to(1,2,\ldots,k+1,k\ldots r)$,
$k=1,2,\ldots r-1$) and by associating to every permutation $\wp$
the morphism in $\catfp$ given by
$$\Pi_\wp=\tau_{k_1}\circ\tau_{k_2}\ldots\circ\tau_{k_r}: \TT^{r,0}\to \TT^{r,0}~,$$
where $\tau_k=\id\otimes_A\ldots\otimes_A\id\otimes_A
\tau\otimes_A\id\otimes_A\ldots\otimes_A\id$ is
a product of $r-1$ factors with the
braiding $\tau=\tau(\RR)$ occurring in the $k^{\rm th}$ factor. Since
the universal $\RR$ matrix is triangular we have $\tau_k^2=\id$
besides the braid relations
$\tau_{k+1}\circ\tau_k\circ\tau_{k+1}=\tau_{k}\circ\tau_{k+1}\circ\tau_{k}$ and
$\tau_k\circ\tau_\ell=\tau_\ell\circ \tau_k$ for $|k-\ell|\geq 2$. Thus the map
$\wp\to \Pi_\wp$ is a well-defined representation of the permutation
group.

The module $\Omega^{ r}(A)$ of $r$-forms is then the image of the projector
\begin{equation}\label{PAproj}
  A_r:=\frac{1}{r!}\sum_{\wp} {\rm{sign}}(\wp)\Pi_\wp:\TT^{r,0}\to \TT^{r,0}~,
\end{equation}
where the sum is over all the $r!$ permutations and  ${\rm{sign}}(\wp)$ is $1$
or $-1$ depending on the even or
odd number of elementary transpositions occurring in $\wp$.

The
module of exterior forms is  $\OM=\bigoplus_{r\in
  \mathbb{N}}\Omega^{ r}(A)$ with $\Omega^0(A)=A$,
$\Omega^1(A)=\Omega(A)$. 
From $\Pi_\wp\circ \tau_k=\Pi_{\wp\;\!\circ t_k}$ it follows that
$A_r\circ \tau_k=-A_r$ and hence $A_r\circ\Pi_\wp={\rm{sign}}(\wp)A_r$ and 
$A_r=A_r\circ (A_\ell\otimes_A \id^{\otimes_A (r-\ell)})=A_r\circ
(\id^{\otimes_A\ell} \otimes_A A_{r-\ell})$ for all $\ell=0,1,\ldots r$.
As in the classical case (see e.g. \cite[\S 6.1]{AMR}) these
properties imply that the wedge product
$\wedge: \OM\otimes\OM\to \OM$ defined on exterior forms of homogenous degree as 
$$\wedge: \Omega^r(A)\otimes \Omega^{r'}(A)\longrightarrow
\Omega^{r+r'}(A)~,~~\theta\otimes\theta'\mapsto
\theta\wedge\theta':=\frac{(r+r')!}{r!\,r'!}A_{r+r'}(\theta
\otimes_A\theta')~$$
is associative and graded braided commutative: $\theta\wedge\theta'= (-1)^{rr'}\oR^\al\ra^\cop\theta'\wedge
\oR_\al\ra^\cop\theta$. 
 The module of
exterior forms with the wedge product becomes the graded braided
commutative $H$-module algebra $(\OM,\wedge)$.
It is generated by $\Omega(A)^0=A$ and 
$\Omega^1(A)=\Omega(A)$ since so is the tensor algebra
$\Tau^{\bullet,0}$, explicitly, let $\theta_1,\theta_2,\ldots\theta_r\in \Omega(A)$,
 then
$\theta_1\wedge\theta_2\ldots\wedge\theta_r=r!\,A_r(\theta_1\otimes\theta_2\ldots\otimes\theta_r)$.\

For later use we also observe that
 the wedge product of exterior forms 
 $\theta\in \Omega^r(A)$ and $\theta'\in
\Omega^{r'}(A)$ also reads
\begin{equation}\label{twt'Arr'}
\theta\wedge\theta'=A_{r,r'}(\theta
\otimes_A\theta')~.
\end{equation}
Here
$A_{r,r'}:=\sum
{{\rm{sign}}(\mathcalligra{s}_{\,r,r'})}\Pi_{\mathcalligra{s}_{\,r,r'}}$
with the sum running over all 
$(r,r')$-shuffles $\mathcalligra{s}_{\,r,r'}$, i.e., permutations
$(1,\ldots r,r+1,\ldots r+r')\to
(\wp(1),\ldots \wp(r),\wp(r+1),,\ldots \wp(r+r'))
$ with $\wp(1)<\wp(2)<\ldots\wp(r)$ and $\wp(r+1)<\wp(r+2)<\ldots
\wp(r+r')$. (This follows form the decomposition of any permutation 
 of $r+r'$ elements as
  ${\mathcalligra{s}}_{\,r,r'}\circ \wp_{r}\circ \wp_{r'}$
where $\wp_{r}$ permutes the
first $r$ elements, $\wp_{r'}$ the last $r'$).
\\

The first order  differential calculus  $(\Omega(A),\dd)$ defined in
\eqref{defofd} and \eqref{defofOmega} uniquely extends to the graded
algebra of exterior forms $\OM$ by defining,
for all
$a, a_1,a_2,\ldots a_r\in A$,
$\dd(\dd  a)=0$ and $\dd(a\dd a_1\wedge\dd a_2\wedge \ldots\dd
a_r)=
\dd a\wedge\dd a_1\wedge\dd a_2\wedge \ldots\dd a_r$. One has to
prove however that the definition is well-defined, i.e.,  if
 $\sum_{i_1,...i_r}a_{i_1,...i_r}\dd a_{i_1}\wedge\dd a_{i_2}\wedge \ldots\dd
a_{i_r}=0$ then
$\sum_{i_1,...i_r}\dd a_{i_1,...i_r}\wedge\dd a_{i_1}\wedge \dd a_{i_2}\wedge \ldots\dd
a_{i_r}=0$ (finite sum understood).

\begin{theorem}
Let $H$ be a  Hopf algebra with triangular structure $\RR$ and $A$ a
braided commutative $H$-module algebra. Let $\Omega(A)$ be the
module in $\cat$ of one-forms as defined in \eqref{defofOmega} and $(\OM,\wedge)$
the graded braided commutative algebra of exterior forms.
There exists one and only one $H$-equivariant map of degree one 
$\dd:
\OM\to \Omega^{\bullet+1}(A)$ such that:
\begin{enumerate}
\item $\dd$ restricted to the degree zero subalgebra $A$  is the original
  derivative $\dd: A\to \Omega(A)$.
  \item It satisfies the graded Leibniz rule
\begin{equation}\label{gLr}
\dd
(\theta\wedge\theta')=\dd\theta\wedge\theta'+(-1)^k\theta\wedge\dd\theta'
\end{equation}
for any $\theta\in\Omega(A)^k$ and $\theta'\in\OM$.
\item  It is nihilpotent, $\dd\circ\dd=0$.
\end{enumerate}
\end{theorem}
  \begin{proof}
 We follow the extended bimodule method of \cite[Thm
 4.1]{Wor}. Let $AX$ be the free left $A$-module with one generator
 $X$. The direct sum $\tilde\Omega(A):=AX\oplus\Omega(A)$ is a left
 $A$-module. Any element in $\tilde\Omega(A)$ has a unique
 decomposition as $bX+\theta$ with $b\in A, \theta\in \Omega(A)$. The
 left $A$-module $\tilde\Omega(A)$ becomes an $A$-bimodule by defining, for any
 $bX+\omega\in \tilde\Omega(A)$, 
 $$(bX+\omega)a=baX+b\dd a+\omega a~.$$ In particular we have $\dd
 a=Xa-aX$. The bimodule $\tilde\Omega(A)$  is canonically a bimodule in $\cat$ by defining, for all $h\in H$, $h\ra (bX+\omega)=(h\ra b)X+h\ra^\cop\omega$, so that the
 generator $X$ is $H$-equivariant, for all $h\in H$, $h\ra X=\epsilon(h)X$. 

 Let $\tilde{\Tau}^{r,0}=\tilde{\Omega}(A)^{\otimes_A r}$ and $\tOM$ be the exterior algebra associated with
 $A$ and $\tilde{\Omega}(A)$; it is defined via the projector
 $\tilde{A}_r:\tilde{\Tau}^{r,0}\to \tilde{\Tau}^{r,0}$ based on the
 braiding $\tau:\tilde{\Omega}(A)\otimes_A
 \tilde{\Omega}(A)\to \tilde{\Omega}(A)\otimes_A
 \tilde{\Omega}(A)$. We show that $\OM$ is a graded subalgebra
 of $\tOM$. The degree zero subalgebras of $\OM$ and of $\tOM$ equal $A$.  The module $\Omega(A)$ is a submodule in $\cat$ of
 $\tilde{\Omega}(A)$, in particular the braiding of
 $\tilde\Omega(A)\oA \tilde\Omega(A)$
 restricts to the braiding of
 $\Omega(A)\oA \Omega(A)$. This implies that
 $\tilde{A}_r|_{\Omega^n(A)}=A_r$ and therefore that the wedge
product of $\tOM$ restricts to the  wedge
product of $\OM$, i.e., $\tilde\wedge|_{\OM\otimes_A\OM}=\wedge$.

For any $\theta\in \OM$ of homogeneous degree $k$ we define
$$
\dd \theta:=X\:\!\tilde\wedge\:\!\theta-(-1)^k\theta\:\!\tilde\wedge\:\! X~.
$$
On $A$ we recover the derivative $\dd$ of the first order differential
calculus. The graded Leibniz rule \eqref{gLr} is easily verified. The property
$\dd\circ\dd=0$ follows from $\dd(\dd \theta)=0$ for any $\theta\in
\Omega^k(A)$ and is due to the $H$-equivariance of
$X$,  $X\:\!\tilde\wedge\:\! X=X\otimes X-\oR_\alpha\ra
X\otimes \oR^\al\ra X=0$. From the graded Leibniz rule and $\dd\circ \dd=0$ we obtain, for all
$a_1,a_2,\ldots a_r\in A$,
$$\dd(a_1\dd a_2\wedge\dd a_3\wedge \ldots\dd
a_r)=
\dd a_1\wedge\dd a_2\wedge\dd a_3\wedge \ldots\dd a_r~.$$ This
 shows that $\dd:\OM\to \OMPu\subset \tOMPu$; it also shows uniqueness of $\dd:\OM\to \OMPu$.
\end{proof}
The triple $(\OM, \wedge, \dd)$ with $\OM=\bigoplus_{r\in
  \mathbb{N}}\Omega^r(A)$, $\Omega^0(A)=A$,  is an $H$-equivariant
braided commutative differential graded algebra because  $(\OM, \wedge)$ is a
graded braided commutative $H$-module algebra and $\dd:
\OM\to \Omega^{\bullet+1}(A)$ is an $H$-equivariant map of degree one
that satisfies $\dd\circ \dd=0$ and the graded Leibniz rule \eqref{gLr}.
It is an {\it $H$-equivariant differential calculus on $A$} because it is generated
in degree zero and one.

\subsection{Cartan Calculus}
In this section and in the rest of the paper we assume $\Omega(A)$ to be finitely generated and
projective, hence  $\Omega(A)=\underline{\Omega}(A)$ is in $\catfp$
with left dual $\Vect(A)$ (cf. Lemma \ref{Om*V} and Corollary
\ref{omega=uomega}) and the contraction operator 
becomes  just the canonical isomorphism $\Vect(A)\simeq ({}^*\Vect(A))^*$ of Proposition
\ref{propo:fgpmodprop}. Besides the evaluation map
   $ \langle~,~\rangle:\Vect(A) \otimes_\AA {\Omega}(\AA)\to
    \AA$ 
we  have the coevaluation map
$$\coev : A\to \Omega(\AA)\otimes_A \Vect(A)~,~~a\mapsto
a\omega^i\otimes_A e_i
$$
 where $\{ e_i\in \Vect(A), \omega^i\in \Omega(A) :
i=1,\ldots, n\}$ is a pair of dual bases for $\Vect(A)$.

We study the contraction operator on tensors and the inner derivative
on forms, the
Lie derivative and the Cartan calculus. While the
existence of these operators and of the Cartan calculus is independent
from assuming $\Omega(A)$ in $\catfp$,
if this is the case, all the relevant modules will be in $\catfp$.

 \begin{proposition}\label{onionprop}
The module $\Tau^{0,r}=\Vect(A)^{\otimes_A r}$ in $\catfp$  is left dual to
$\Tau^{r,0}=\Omega(A)^{\otimes_A r}$, with
evaluation and coevaluation maps
\begin{equation*}
\begin{split}
\le~,~\re &
: \TT^{0,r} \otimes_\AA\TT^{r,0}\;\to\;
A~,~~
\le v_r\otimes_A\ldots v_2\otimes_A
v_1,\omega_1\otimes_A\omega_2\otimes_A\ldots \omega_r\re=\\[.2em]
&~~~~~~~~~~~~~~~~~~~~~~~~~~~~~~~~~~~~~~~~~~~~~~~~~~~~~~~~~~=
\le v_r,\ldots\le v_2,\le v_1,\omega_1\re\omega_2\re\ldots\omega_r\re \\[.2em]
\coev &: A \;\to\;\TT^{r,0}\otimes_\AA  \TT^{0,r} ~,~~\coev(a)=a_{\,}
\omega^{i_1}\otimes_A\omega^{i_2}\ldots\otimes_A\omega^{i_r}\otimes_A
e_{i_r}\otimes_A \ldots e_{i_2}\otimes_A e_{i_1}~.
\end{split}
\end{equation*}
\begin{proof}
For $r=1$ this is rigidity of $\Omega(A)$. The proposition follows from $\Omega(A)={}^*\Vect(A)$ and, by iteration,
from Remark \ref{newevcoev}.  
\end{proof}
\end{proposition}
Recalling the algebra structure \eqref{minT} of $\Tau^{\bullet,
  \bullet}$, we have that  $\TT^{p,q}=\TT^{r,0}\otimes_A \TT^{p-r,q}$ (if
$p\geq r$) and we
immediately extend the pairing $\le~,~\!\!\:\re$ to the evaluation of $\TT^{0,r}$ on
$\TT^{p,q}$. This is the morphism in $\catfp$ defined to be trivial if
$r>p$ and otherwise given by
\begin{equation}\label{evTT}
\le~,~\;\re 
: \TT^{0,r} \otimes_\AA\TT^{p,q}\;\to\;
\TT^{p-r,q}~,~~\le\nu,\theta\otimes_A\eta\re:=\le\nu,\theta\re\eta
\end{equation}
for all $\nu\in \Tau^{0,r}$, $\theta\in \Tau^{r,0}$, $\eta\in
\Tau^{p-r,q}$.
In particular, for $r=1$  we obtain the extension of the contraction
operator to $\ii:\Vect(A)\to \hom_A(\TT^{p,q},\TT^{p-1,q})$. Hence the
evaluation \eqref{evTT} is just the iteration of the contraction operator 
$r$-times:
$\le v_r\otimes_A\ldots v_2\otimes_A
v_1,\eta\re=
\ii_{v_r}\circ \ldots \ii_{v_2}\circ\ii_{v_1}(\eta)$. 
\\

The exterior algebra construction applied to $\Vect(A)$ rather than
$\Omega$ gives the graded braided
antisymmetric $H$-module algebra $\Vect^\bullet(A):=\bigoplus_{r\in
 \mathbb{N\,}}\Vect^{r}(A)$ of polyvector fields.  
The duality of Proposition \ref{onionprop} then implies the duality
between the graded modules $\OM$ and $\Vect^\bullet(A)$.
\begin{proposition}\label{onionpropwedge}
The module $\Vect(A)^{r}\subset\Vect(A)^{\otimes_A r}=\Tau^{0,r}$ in
$\catfp$ is left dual, for any $r\in \mathbb{N\,}$,  to
$\Omega^r(A)\subset\Omega(A)^{\otimes_A r}=\Tau^{r,0}$. The
evaluation map 
$\le~,~\re^{}_{\Vect(A)^{r}}:\Vect(A)^{r}\otimes_A\Omega^r(A)\to A$ is the restriction to $\Vect(A)^{r}\otimes_A\Omega^r(A)$ of the
evaluation map in  Proposition \ref{onionprop}, the coevaluation map is
\begin{equation*}
\coev^{}_{\Vect(A)^{r}} : A \;\to\;\Omega^r(A)\otimes_\AA \Vect(A)^{r} ~,~~\coev^{}_{\Vect(A)^{r}} (a)=\frac{1}{(r!)^2}a_{\,}
\omega^{i_1}\wedge\omega^{i_2}\ldots\wedge\omega^{i_r}\otimes_A
e_{i_r}\wedge \ldots e_{i_2}\wedge e_{i_1}~.
\end{equation*}
\begin{proof}
We  begin by observing that
  \begin{align}
\omega^{s_1}&\wedge\,\omega^{s_2}\ldots\wedge\,\omega^{s_r}\otimes_A
                                                            e_{s_r}\otimes_A \ldots e_{s_2}\otimes_A e_{s_1}=\nn\\
    &\,=
\omega^{i_1}\otimes_A\omega^{i_2}\ldots\otimes_A\omega^{i_r}\otimes_A
\le e_{i_r}\otimes_A \ldots e_{i_2}\otimes_A e_{i_1}\,,\,
\omega^{s_1}\wedge\omega^{s_2}\ldots\otimes_A\omega^{s_r}\re\,
      e_{s_r}\otimes_A \ldots e_{s_2}\otimes_A e_{s_1}\nn\\
     &\,=
\omega^{i_1}\otimes_A\omega^{i_2}\ldots\otimes_A\omega^{i_r}\otimes_A
\le e_{i_r}\wedge\,\ldots e_{i_2}\wedge\, e_{i_1}\,,\,
\omega^{s_1}\otimes_A\omega^{s_2}\ldots\otimes_A\omega^{s_r}\re\,
       e_{s_r}\otimes_A \ldots e_{s_2}\otimes_A e_{s_1}\nn\\
    &=\omega^{s_1}\otimes_A\omega^{s_2}\ldots\otimes_A\omega^{s_r}\otimes_A
                                                            e_{s_r}\wedge\,\ldots e_{s_2}\wedge\, e_{s_1}
\label{prelobsome}
  \end{align}
where we first rewrote
$\omega^{s_1}\wedge\,\omega^{s_2}\ldots\wedge\,\omega^{s_r}$ using
that $(\id_{\Tau^{r,0}}\otimes_A\le~,~\re)\circ
(\coev\otimes_A\id_{\Tau^{r,0}})=\id_{\Tau^{r,0}}$
(cf. \eqref{coerecoev}); then we used that
the adjoint of the projector
$A_r:\TT^{r,0}\to \TT^{r,0}$ is the projector
$  A_r:\TT^{0,r}\to \TT^{0,r}$, cf. \eqref{taudual}, and finally that
 $(\le~,~\re\otimes_A \id_{\Tau^{0,r}})\circ
 (\id_{\Tau^{0,r}}\otimes_A\coev)=\id_{\Tau^{0,r}}$.
 
 We now prove that 
 $\le~,~\re^{}_{\Vect(A)^{r}}$
 and
 $\coev^{}_{\Vect(A)^{r}}$
 are evaluation and coevaluation maps. We have
 $(\le~,~\re^{}_{\Vect(A)^{r}}\otimes_A \id_{\Vect(A)^{r}})\circ
(\id_{\Vect(A)^{r}}\otimes_A \coev^{}_{\Vect(A)^{r}})=\id_{\Vect(A)^{r}}$,
indeed, for all
 $u_r\wedge \ldots u_2\wedge u_1\in {\Vect(A)^{r}}$,
   \begin{align}
\frac{1}{(r!)^2}\le    u_r\wedge  \ldots u_2\wedge u_1,
\omega^{s_1}&\wedge\,\omega^{s_2}\ldots\wedge\,\omega^{s_r}\re^{}_{_{ \!\!\;\Vect(A)^{r}}}\,
                                                            e_{s_r}\wedge \ldots e_{s_2}\wedge e_{s_1}=\nn\\
&=\frac{1}{r!}\le     u_r\wedge  \ldots u_2\wedge u_1,
A_r(\omega^{s_1}\oA\,\omega^{s_2}\ldots\oA\omega^{s_r})\re^{}_{_{ \!\!\;\Vect(A)^{r}}}\,
                                                            e_{s_r}\wedge
                                                                                                            \ldots e_{s_2}\wedge e_{s_1}\nn\\
     &=\frac{1}{r!}\le     A_r(u_r\wedge  \ldots u_2\wedge u_1),
\omega^{s_1}\oA\,\omega^{s_2}\ldots\oA\omega^{s_r}\re\,
                                                            e_{s_r}\wedge
       \ldots e_{s_2}\wedge e_{s_1}\nn\\
          &=\frac{1}{r!}\le     u_r\wedge  \ldots u_2\wedge u_1,
\omega^{s_1}\wedge\omega^{s_2}\ldots\wedge\omega^{s_r}\re\,
                                                            e_{s_r}\oA
            \ldots e_{s_2}\oA e_{s_1}\nn\\
          &=\le     A_r(u_r\wedge  \ldots u_2\wedge u_1),
\omega^{s_1}\oA\,\omega^{s_2}\ldots\oA\omega^{s_r}\re\,
                                                            e_{s_r}\oA
       \ldots e_{s_2}\oA e_{s_1}\nn\\
          &=\le     u_r\wedge  \ldots u_2\wedge u_1,
\omega^{s_1}\oA\omega^{s_2}\ldots\oA\omega^{s_r}\re\,
                                                            e_{s_r}\oA
            \ldots e_{s_2}\oA e_{s_1}\nn\\
     &=u_r\wedge \ldots u_2\wedge u_1~,\nn
   \end{align}
   where in the third line we used that the adjoint of the projector
   $A_n:\Tau^{n,0}\to \Tau^{n,0}$ is the projector $A_n:\Tau^{0,n}\to
   \Tau^{0,n}$, cf. \eqref{taudual}, in the fourth, equation
   \eqref{prelobsome}, and
   in the fifth we considered again the adjoint of the projector
   $A_n:\Tau^{n,0}\to \Tau^{n,0}$.

Similarly one proves that $(\id_{\Omega^r(A)}\otimes_A\le~,~\re^{}_{\Vect(A)^{r}})\circ
(\coev^{}_{\Vect(A)^{r}}\otimes_A\id_{\Omega^r(A)})=\id_{\Omega^r(A)}$.
\end{proof}
\end{proposition}
As in the commutative case the contraction operator restricted to
forms becomes the inner derivative.
\begin{proposition}
The contraction operator
$\ii:
\VV\otimes_A\Omega^{\otimes_A r}\to\Omega^{\otimes_A r-1}$ restricts to
\begin{equation}\label{innerder}
  \ii:\VV\otimes_A\Omega^{r}\to\Omega^{r-1}~,
\end{equation}
giving, for all $u\in
\Vect(A)$,
the (graded) braided inner derivative
$\ii_u:\OM\to\OMmu(A)$, 
\begin{equation}\label{lieinnerder}\ii_u(\theta\wedge\theta')=\ii_u(\theta)\wedge\theta'+(-1)^{|\theta|}(\oR^\al\ra^\cop
\theta)\wedge
\ii_{\oR_\al\ra u}\theta'~,
\end{equation}
where $|\theta|\in \mathbb{N}$ is the degree of the homogeneous form
$\theta$, while $\theta'\in \OM$. 
Furthermore, for all $u, v\in \Vect(A)$, on $\OM$ we have
\begin{equation}\label{ii}
\ii_u\circ \ii_v+\ii_{\oR^\al\ra v}\circ \ii_{\oR_\al \ra u}=0~.
\end{equation}
\end{proposition}
\begin{proof}
Use that $\Omega(A)$ is a braided
commutative $A$-bimodule, $K$-equivariance of $\,\ii$, 
the  property    $ (\mathrm{id}\otimes\Delta)\mathcal{R}
=\mathcal{R}_{13}\mathcal{R}_{12}$
and triangularity of the universal $\RR$-matrix
to prove the first of the following equalities, for all $u\in \Vect(A)$, $\omega\in \Omega(A)$,
$\eta\in \Tau^{0,\bullet}$ and 
integers $k\geq 2$,
\begin{equation}\label{iut1}
  \ii_u\circ \tau_1(\omega\otimes_A \eta)
  =\oR^\alpha\ra^\cop \omega\otimes_A \ii_{\oR_\alpha \ra u} \eta~\;;~~~~ i_u\circ \tau_k=\tau_{k-1}\circ i_u
  \end{equation}
where $\tau_1$ is the braiding in the first two factors of the tensor
product, while $\tau_k$ is that on the $k^{\rm th}$ and $k^{\rm th}+1$
factor in $\Tau^{0,r}$, $r\geq k+1$. 
The second equality is the commutativity of the contraction with the braiding,
which holds whenever the braiding does not involve the first factor of the tensor product.

 We now prove \eqref{innerder} by induction on the degree of exterior forms.
Equation \eqref{innerder} is
 trivially true for $r=0,1$. We assume it holds for $r$ and show it
 holds for $r+1$.
 Let $\omega\in \Omega(A)$, $\theta\in \Omega^r(A)$; from \eqref{twt'Arr'}
 we have
 $\omega\wedge\theta=A_{1,r}(\omega\oA \theta)$ where
$A_{1,r}$ represents the alternating sum of all $(1,r)$-shuffles,
 \begin{align*}
 A_{1,r}&=\id^{\otimes_{\! A} r}-\tau_1+\tau_2\circ\tau_1-\tau_3\circ\tau_2\circ\tau_1+\,\ldots\,
            (-1)^{r-1}\tau_{r}\circ\tau_{r-1}\,...\circ\tau_1\\
            &=\id^{\otimes_{{\! A}}
   r}-(\id\otimes_A A_{1,r-1})\circ \tau_1~.
 \end{align*}
 From the second equality in \eqref{iut1}, $\ii_u\circ
 A_{1,r}=\ii_u- \ii_u\circ (\id\oA A_{1,r-1})\circ \tau_1=\ii_u- A_{1,r-1}\circ \ii_u\circ \tau_1$, hence, for all
 $u\in \Vect(A)$,
 \begin{align}
   \ii_u(\omega\wedge\theta)&=\ii_u\circ A_{1,r}(\omega\oA\theta)\nn\\
&=\ii_u(\omega)\:\!\theta-A_{1,r-1}(\oR^\alpha\ra^\cop \omega\:\!\oA \ii_{\oR_\alpha
                                                                      \ra
                                                                      u}\theta)\nn\\
   &=\ii_u(\omega)\:\!\theta-(\oR^\alpha\ra^\cop \omega)\wedge \ii_{\oR_\alpha    \ra  u}\theta
     \label{iuomtheta}
 \end{align}
where in the second line we used the first equality in
\eqref{iut1}. 
The inductive hypothesis implies that this  expression is in
$\OM$.
\\

The braided Leibniz rule \eqref{lieinnerder} is also easily proven by
induction. It holds for $|\theta|=0,1$ (cf. \eqref{iuomtheta}). We
assume that it holds for $|\theta|\leq r$  and show that it holds also
for forms of homogenous degree $r+1$. This follows from, for all $\omega\in
\Omega(A)$, $\theta'\in\OM$, $u\in \Vect(A)$,
\begin{align*}
  \ii_u(\omega\!\!\:\wedge\!\!\:\theta_{\!\!\:}\wedge&\theta')
 =
 \ii_u(\omega)\:\!\theta\wedge\theta'-(\oR^\alpha\ra^\cop \omega)\wedge \ii_{\oR_\alpha    \ra  u}(\theta\wedge\theta')\\
&=
 \ii_u(\omega)\:\!\theta\wedge\theta'-(\oR^\alpha\ra^\cop\omega)\wedge\ii_{\oR_\alpha
                                                                                                          \ra
                                                                                                          u}(\theta)\wedge\theta'+(-1)^{k+1}\!\!\:(\oR^\alpha\ra^\cop\omega)\!\!\:\wedge\!\!\:(\oR^\beta\ra^\cop\theta)\!\!\:\wedge\!\!\:\ii_{\!\!\:\oR_\beta\!\:\!\oR_\alpha\ra u}\theta'\\                                   
&=\ii_u(\omega\wedge\theta)\wedge\theta'+(-1)^{k+1}(\oR^\alpha\ra^\cop(\omega\wedge\theta))\wedge\ii_{\oR_\alpha\ra u}\theta' 
\end{align*}
where in the first equality we used \eqref{iuomtheta}, in the second
 the inductive hypothesis, in the third again \eqref{iuomtheta}
and where we have rewritten the last addend using $(\Delta\otimes
\id)\RR^{-1}=\RR^{-1}_{23}\RR^{-1}_{13}$.
\\

Finally, equality \eqref{ii} trivially holds on forms of homogenous
degree $0$ and $1$, 
and for higher forms $\theta\in
\Omega^r(A)$, $r\geq 2$ it follows from 
$$(\ii_u\circ \ii_v+\ii_{\oR^\al\ra v}\circ \ii_{\oR_\al \ra u})\theta=\langle (\id+\tau)(u\otimes_A
v),\theta\rangle=
\langle u\otimes_A
v,(\id+\tau_1)\theta\rangle=0$$
where in the second passage we used \eqref{taudual} and in the last
that $\tau_1\circ A_r=  \frac{1}{r!}\sum_{\wp} {\rm{sign}}(\wp)\:\!\tau_1\circ
\Pi_\wp=\frac{1}{r!}\sum_{\wp} {\rm{sign}}(\wp)\Pi_{t_1\circ \wp}=-A_r$. 
\end{proof}

We next study the action of the $(H,\RR)$-Lie algebra of vector fields
$\Vect(A)$ on tensor fields, i.e.,  the Lie derivative. On  $A$ and $\Vect(A)$ we define $$\LLL:\VV\otimes A\to A~, ~~\LLL_u(a):=u(a) 
~;~~~\LLL:\VV\otimes \VV\to \VV~,~~\LLL_u(v):=[u,v]~.$$ 
Since for any $h\in H$, $\hxi\ra (\LLL_u(a))=\LLL_{\hxi_{(1)}\ra u}(\hxi_{(2)}\ra a)$ and
$\hxi\ra (\LLL_u(v))=\LLL_{\hxi_{(1)}\ra u}(\hxi_{(2)}\ra v)$, $\LLL$ is $H$-equivariant.
It is not difficult to check its compatibility with the $A$-bimodule structure
\eqref{aub} of $\Vect(A)$,
$\LLL_u(av)=\LLL_u(a)v+\oR^\al\ra a\,\LLL_{\oR_\al\ra u}(v)$,
$\LLL_u(va)=\LLL_u(v)a+\oR^\al\ra v\,\LLL_{\oR_\al\ra u}(a)$.
We extend  $\LLL$  to an $H$-equivariant
action of $\Vect(A)$ on the tensor product
$\Vect(A)\otimes\Vect(A)$ by defining,  
$\LLL_u(v\otimes z)=
\LLL_u(v)\otimes z+
\oR^\alpha\ra v\otimes \LLL_{\oR_\al\ra u}z$
that is, the action of $u\in \Vect(A)$ is given by
$\LLL_u\otimes_\RR\id_{\Vect(A)}+\id_{\Vect(A)}\otimes_\RR\LLL_u$;
this is a well defined internal morphism in the symmetric category
${}^H\MMM$ in  accordance with \eqref{otimesR}, in particular it
commutes with the braiding $\tau$.
The compatibility of $\LLL$ with the $A$-bimodule structure of
$\Vect(A)$ implies that also the following action on the balanced
tensor product over $A$ is well defined (that is, independent from the representatives):
$\LLL_u(v\otimes_A z)=
\LLL_u(v)\otimes_A z+
\oR^\alpha\ra v\otimes_A \LLL_{\oR_\al\ra u}z$.
It is extended to the whole tensor algebra $\Tau^{0,\bullet}$
by iterating this braided derivation rule,
 i.e, 
for all $\nu,\nu'\in  \Tau^{0,\bullet}$,
$$\LLL_u(\nu\otimes_A\nu')=\LLL_u(\nu)\otimes_A
\nu'+\oR^\al\ra \nu\otimes_A\LLL_{\oR_\al\ra u}(\nu')~.$$
The Lie derivative on contravariant tensor fields
is canonically defined by duality, for all $\nu\in \Tau^{0,r}$ and
$\theta\in  \Tau^{r,0}$,
\begin{equation}\label{LonOm} \LLL_u\langle \nu,\theta\rangle=\langle\LLL_u\nu,\theta\rangle
+\langle \oR^\al\ra\nu,\LLL_{\oR_\al\ra u}\theta\rangle
\end{equation} 
i.e., $\langle \nu,\LLL_{u}\theta\rangle:=\LLL_{\oR^\al \ra u}\langle \oR_\al\ra
\nu,\theta\rangle-\langle\LLL_{\oR^\al\ra u}\oR_\al \ra\nu,\theta\rangle
$.
It follows that vector fields act  on the tensor algebra
$\Tau^{\bullet,\bullet}$ as braided derivations. 
On tensor fields  $\Tau^{\bullet,\bullet}$ we have, for all $u, v\in \VV$, 
\begin{equation}\label{LLL}
\LLL_u\circ \LLL_v-\LLL_{\oR^\al \ra v}\circ \LLL_{\oR_\al \ra
  u}=\LLL_{[u,v]}~.
\end{equation}
On $A$ this is the definition of the braided commutator $[u,v]$; on
$\Vect(A)$ it is just the  Jacobi identity in \eqref{brLie}, then
\eqref{LonOm} implies that it  holds on $\Omega(A)$. Since
both sides of \eqref{LLL} are braided derivations this equality
extends to all  $\Tau^{\bullet,\bullet}$.
Equation \eqref{LLL} shows
that the Lie derivative $\LLL:\VV\otimes\Tau^{\bullet,\bullet}\to  \Tau^{\bullet,\bullet}$ 
is an action of the $(H,R)$-Lie algebra of derivations $\VV$ on  $\Tau^{\bullet,\bullet}$.

The commutativity with the braiding, $\LLL_u\circ\tau=\tau\circ\LLL_u$, and the wedge
product expression \eqref{twt'Arr'} imply that vector fields  also act as braided derivations on the exterior algebras
$\OM$ and 
$\Vect^{\bullet\!}(A):=\bigoplus_{n\in \mathbb{N}}\Vect^n(A)$ (thar
are both $H$-submodules of  $\Tau^{\bullet,\bullet}$).
From \eqref{LonOm} it is immediate to compute, for all $u,v\in
\Vect(A)$, $\omega\in \Omega(A)$, 
\begin{equation}\label{Lii}
(\LLL_u\circ \ii_v-\ii_{\oR^\al \ra v}\circ \LLL_{\oR_\al \ra
  u})\omega=\ii_{[u,v]}\omega~;
\end{equation}
since both left hand side and right hand side are braided derivations
on $\OM$, this relation extends to arbitrary exterior forms.
\\

The Lie derivative commutes with the exterior derivative on $A$, for
all $a\in A$, $v\in \Vect(A)$, 
$\LLL_v\dd a=\dd \LLL_v a$, indeed, for all $u \in \Vect(A)$,
\begin{equation*}
  \begin{split}\le u,\LLL_v\dd a\re&=\LL_{\oR^\al\ra v}\le\oR_\al\ra u,\dd
  a\re-\le[\oR^\al\ra v,\oR_\al\ra u],\dd
  a\re=\LL_{\oR^\al \ra v}\LL_{\oR_{\al\ra u}} a-\LL_{[\oR^\al\ra v,\oR_\al\ra
    u]}a\\[.2em]
  &=\LL_u\LL_va
  =\le u,\dd\LL_v a\re~.
\end{split}
\end{equation*}
Using induction on the form degree we have that $\dd
\LLL_v\theta=\LLL_v\dd\theta$ for any $\theta\in \OM$.

Similarly, for all $v\in \Vect(A)$,
$$\LLL_v=\ii_v\circ \dd+\dd\circ \ii_v$$ trivially holds on $A$ and  by
induction on the form degree it holds on
$\OM$ since both the right hand side and the left hand side are
braided derivations of
$\OM$.
\\

The equations $\LLL_z\circ \dd =\dd\circ \LLL_z$, $\LLL_v=\ii_v\circ \dd+\dd\circ \ii_v $ and  \eqref{ii}, \eqref{LLL} (restricted to exterior forms),
\eqref{Lii}, $\dd^2=0$,  constitute the Cartan calculus of the
exterior, Lie and inner derivatives \cite{Gurevich, Weber}. Notice that the derivation of these
equations holds true also if $\Omega(A)$ is not finitely generated and
projective over $A$ (indeed we never used coevaluation maps, just
nondegeneracy of the pairing, cf. Corollary \ref{pnondeg}). We summarize
these equations in the following theorem,

\begin{theorem}[Braided Cartan calculus]
Let ${A}$ be a braided commutative left $H$-module algebra
and consider the associated braided differential algebra
$(\OM,\wedge,\mathrm{d})$. 
The exterior derivative, the Lie derivative and inner derivative along vector fields
$u,v\in \Vect(A)$ are graded braided derivations of
$\OM$ (respectively of degree $1,0,
-1$) that 
 satisfy
\begin{equation*}
\begin{split}
        [\LLL_u,\LLL_v]
        =&\,\LLL_{[u,v]},~~\\
        [\LLL_u,\ii_v]
        =&\,\ii_{[u,v]},\\
        [\LLL_u,\mathrm{d}]=&\,0\,,
\end{split}
\begin{split}
        [\ii_u,\ii_v]=&\,0\,,\\
        [\ii_u,\mathrm{d}]
        =&\,\LLL_u,\\
        [\mathrm{d},\mathrm{d}]=&\,0\,,
\end{split}
\end{equation*}
where 
$[L,L']=L\circ L'-(-1)^{|L||L'|}\oR^\al\ra L'\circ \oR_\al\ra L$
is the graded braided commutator of $\Bbbk$-linear maps $L,L'$  on
$\OM$ of
degree $|L|$ and $|L'|$.
\end{theorem}

\subsection{Examples}\label{exgeneral}
\begin{example}\label{excotriang}{\it Braided derivations of a cotriangular Hopf algebra $K$
  define a bicovariant differential calculus \`a la Woronowicz.} 
Let $A=K$ be a finite dimensional cotriangular Hopf algebra over a
field $\Bbbk$. Let $U$ be the dual triangular Hopf algebra with $\RR$-matrix
$\mathscr{R}=\mathscr{R}^\al\otimes{\mathscr{R}}_\al$, inverse
${\mathscr{R}}^{-1}=\overline{\mathscr{R}}^\al\otimes\overline{\mathscr{R}}_\al$
and antipode $S$.
Recall from Example \ref{cotriangularK} that $K$ is a $U^\op\otimes
U$-module algebra and is braided commutative with
${\RR}=R^\alpha\otimes R_\alpha=(\id\otimes
\rm{flip}\otimes\id)(\mathscr{R}^{-1}\otimes\mathscr{R})$.
The braided derivations
\begin{equation}\label{brderK}
  \Der_\RR(K)=\{u\in \hom_\Bbbk(K,K)\,|\,
u(k\ell)=u(k)\ell\,+\,\oR^\alpha\ra k \,(\oR_\alpha \RA u)(\ell) \,,\,\mbox{ for
  all }k,\ell\in K\}
\end{equation}are then a relative $U^\op\otimes U$-module
$K$-bimodule. We set $\Vect(K):=\Der_\RR(K)$ and call it the
module of (braided) vector fields. As in Example \ref{cotriangularK}, since $K$ is finite
dimensional over the field $\Bbbk$, $\Vect(K)$ is equivalently a
$K$-bicovariant bimodule, it is therefore free over $K$.
In particular, the
adjoint left
$U^\op$-action 
$\ra: U^\op\otimes  \Vect(K)\to  \Vect(K)$, $\zeta\otimes
u\mapsto\zeta\ra u$, where $(\zeta\ra u)(k)=\zeta_{(1)}\ra (u(S^{-1}(\zeta_{(2)})\ra k))$
for all $k\in K$, with $\zeta\ra k=\zeta(k_{(1)})k_{(2)}$ and $S^{-1}$
the antipode in $U^\op$, is dual to the left $K$-coaction
\begin{equation}\label{deltaL}\Delta_\Ls: \Vect(K)\to
  K\otimes\Vect(K)~,~~ u\mapsto  \Delta_\Ls(u)=u_{-1}\otimes u_0
  \end{equation}
defined by
$u_{-1}\otimes u_0(k):=u(k_{(2)})^{}_{(1)}S^{-1}(k_{(1)})\otimes u(k_{(2)})^{}_{(2)}$ for
all $k\in K$. (In order to prove that
$\zeta\ra u=\zeta(u_{-1})u_0$ for all $\zeta\in U^\op$, evaluate
both members on $k\in K$).

Let
${}_{\rm{inv}}\Vect(K)\subset \Vect(K)$ be the $\Bbbk$-submodule of
left-invariant vector fields, i.e., of vector fields that under the
adjoint left
$U^\op$-action 
transform in the trivial representation: $\zeta\ra
u=\varepsilon(\zeta)u$, for all $\zeta\in U^\op$; equivalently, of
vector fields invariant under the coaction
$\Delta_\Ls$:  $\Delta_\Ls(u)=1_K\otimes u$.
The $K$-bicovariant bimodule structure of $\Vect(K)$ implies
\begin{equation}\label{XXinv}
  \Vect(K)=K\otimes{}_{\rm{inv}}\Vect(K)~.
  \end{equation}The dual $K$-module of one-forms is the 
$K$-bicovariant bimodule $\Omega(K):={}^\ast \Vect(K)$. The differential
$\dd: K \to  \Omega(K)$ is defined as in 
\eqref{defofd}; left $K$-linearity of the pairing $\le~,~\re:
\Vect(K)\otimes_K \Omega(K)\to K$ implies that $\dd$ is
determined by the left-invariant vector fields
${}_{\rm{inv}}\Vect(K)$.
Let $\{\uXu_j\}$, $j=1,2,\ldots n$, be a basis of ${}_{\rm{inv}}\Vect(K)$ and $\{\omega^j\}$
the dual basis of left-invariant one-forms, $\le \uXu_j,\omega^i\re=\delta_j^i$;
from $\omega^j\le u_j,~\,\re=\id_{\Omega(K)}$  (sum over $j$ understood) we have
\begin{equation}\label{defofbicdd}
  \dd k=\omega^j\;\!\uXu_j(k)~.
\end{equation}

We study the module of left-invariant vector fields ${}_{\rm{inv}}\Vect(K)$ and prove
that $\dd: K \to  \Omega(K)$ defines a bicovariant
differential calculus \`a la Woronowicz (in particular this implies
the surjectivity property $\Omega(K)=K\dd K$,
cf. \eqref{defofOmega}). 

Due to left-invariance, the $U^\op\otimes U$-action on ${}_{\rm{inv}}\Vect(K)$ reduces to the $U$-action, so that these vector fields
satisfy the braided derivation property, for all $k,\ell\in K$, $u\in {}_{\rm{inv}}\Vect(K)$,
\begin{equation}\label{beqXRr}
\uXu(k\ell)=\uXu(k)\ell\,+\,\oR^\alpha\ra k \,(\oR_\alpha \RA \uXu)(\ell)=\uXu(k)\ell\,+\,\overline{\mathscr{R}}^\al\ra k \,(\overline{\mathscr{R}}_\al \RA \uXu)(\ell)~.
\end{equation}
This shows that the $U$-module ${}_{\rm{inv}}\Vect(K)$ is a
$(U,\mathscr{R})$-Lie algebra: that of left-invariant vector
fields. 
Let $\gg$ be the image of
${}_{\rm{inv}}\Vect(K)$ under the $\Bbbk$-linear map ${}_{\rm{inv}}\Vect(K)\to U$,
$\uXu\mapsto\chi^{}_\uXu:=\varepsilon\circ \uXu$.
The identity $u(k)=k_1 \chi^{}_\uXu(k_2)$ for all $k\in K$ provides
the inverse map, hence the isomorphism
\mbox{${}_{\rm{inv}}\Vect(K)\simeq\gg$}. This identity follows applying
$\id\otimes \varepsilon$ to
$\Delta(u(k))=k_{(1)}\otimes u(k_{(2)})$. This latter identity in turn
is
equivalent to the left-invariance condition  $\Delta_\Ls(u)=1_K\otimes u$.
Indeed, for all $k\in K$, $u_{-1}\otimes u_0(k)=1_K\otimes
u(k)$ is equivalent to   $u_{-1}k_{(1)}\otimes u_0(k_{(2)})=k_{(1)}\otimes
u(k_{(2)})$  which,  recalling the
definition of $\Delta_\Ls$, reads $\Delta(u(k))=k_{(1)}\otimes u(k_{(2)})$. 

Furthermore, ${}_{\rm{inv}}\Vect(K)\simeq\gg$
is a $U$-module isomorphism, that is,  for all
$\xi\in U, u\in {}_{\rm{inv}}\Vect(K)$,  we have the
$U$-equivariance $\chi^{}_{\xi\ra \uXu}=\xi\ra\chi^{}_\uXu$, where $\ra: U\otimes U\to U$,
$\xi\otimes\zeta\mapsto\xi\ra\zeta=\xi_1\zeta S(\xi_2)$ is the adjoint $U$-action.
Indeed, for all $k\in K$ (and using the standard  pairing notation
$\le~,~\re : U\otimes K\to \Bbbk$),
\begin{equation*}\begin{split}
  \chi^{}_{\xi\ra \uXu}(k)&=\varepsilon((\xi\ra \uXu)(k))\\
  &=\varepsilon (\xi_{(1)}\ra
  (\uXu(S(\xi_{(2)})\ra k)))\\
  &=\varepsilon(\uXu(k_{(1)})^{}_{(1)})\le \xi_{(1)},\uXu(k_{(1)})^{}_{(2)}\re\le
  S(\xi_{(2)}), k_{(2)}\re\\
  &=\le\xi_{(1)},\uXu(k_{(1)})\re\le S(\xi_{(2)}),k_{(2)}\re
\end{split}
\end{equation*}
equals
\begin{equation*}\begin{split}
  (\xi\ra \chi^{}_{\uXu})(k)&=(\xi_{(1)}\chi^{}_\uXu S(\xi_{(2)}))(k)\\
  &=\le \xi_{(1)}, k_{(1)}\re \varepsilon(\uXu(k_{(2)}))\le
  S(\xi_{(2)}), k_{(3)}\re\\
  &=\le \xi_{(1)}, \uXu(k_{(1)})^{}_{(1)}\re \varepsilon(\uXu(k_{(1)})^{}_{(2)})\le
  S(\xi_{(2)}), k_{(2)}\re\\
      &=\le\xi_{(1)},\uXu(k_{(1)})\re\le S(\xi_{(2)}),k_{(2)}\re
\end{split}
\end{equation*}
where in the third line we used that $\uXu\in {}_{\rm{inv}}\Vect(K)$
is left-invariant, i.e.,
$\Delta(u(k))=k_{(1)}\otimes u(k_{(2)})$.

$U$-equivariance of the  isomorphism ${}_{\rm{inv}}\Vect(K)\simeq\gg$ and the braided derivation property
\eqref{beqXRr} imply the braided derivation property $\chi^{}_\uXu(k\ell)=
 \chi^{}_\uXu(k)\varepsilon(\ell)\,+\,\overline{\mathscr{R}}^\al(k)
 \,(\overline{\mathscr{R}}_\al \RA \chi^{}_\uXu)(\ell)$ for all $k,
 \ell\in K$, hence $\gg$
 is the $U$-module and $(U,\mathscr{R})$-Lie algebra
\begin{equation}\label{gforbic}
  \gg=\{\chi\in U; \Delta^{}_U\chi=\chi\otimes
  1^{}_U+\overline{\mathscr{R}}^\al\otimes\overline{\mathscr{R}}_\al\ra\chi\}~.
\end{equation}
From $(
  \varepsilon^{}_U\otimes\id)\Delta^{}_U\chi=\chi$ and \eqref{gforbic} we see that $\gg\subset
  \ker\,\varepsilon^{}_U$.  In terms of the unit and counit
of $K$ we have that the elements of the $U$-module $\gg$ satisfy, for all $\chi\in\gg$, $\chi(1_K)=0$,
$\Delta\chi-\chi\otimes \varepsilon_K\in U\otimes \gg$. Equivalently, $S(\gg)$ 
is a right $U$-module under the adjoint action
$\chi'\triangleleft\xi:=S(\xi_{(1)})\chi'\xi_{(2)}$ and its elements satisfy: 
for all $\chi'\in S(\gg), \chi'(1_K)=0,
\Delta\chi'- \varepsilon_K\otimes\chi'\in S(\gg)\otimes U$.
This proves that $S(\gg)$ is a quantum Lie algebra  associated with a bicovariant differential calculus
\`a la Woronowicz on the Hopf algebra $K$ \cite[\S 14.2.3, Coroll. 10]{KS}. The differential calculus is the one defined by
$\gg\simeq {}_{\rm{inv}}\Vect(K)$ as in \eqref{defofbicdd}. This follows from \cite[Thm.
5.2]{Wor} where it is proven that 
$\dd k=(\chi'_j\ra k)\omega^j$
with $\{\chi'_j\}$ a basis of $S(\gg)$. Indeed choosing
$\chi'_j=-S(\chi^{}_{\uXu_j})$ we have $\dd k=(\chi'_j\ra
k)\omega^j=\omega^j\;\!\uXu_j(k)$ as in \eqref{defofbicdd}; for the last equality  see for example
\cite[eq. (2.17)]{AschieriSchupp}.
\endex\end{example}
\begin{remark} 
In Section \ref{H4DC} the example of the 
 Sweedler Hopf algebra is explicitly presented.
\end{remark}

\begin{remark} 
The bicovariant calculus on the cotriangular Hopf algebra
$(K,\mathscr{R})$ that we have constructed and which is determined by
$\gg$ is equivalent to the one
defined in \cite[\S 4.3]{Gomez-Majid}.
Indeed we have, 
\begin{align}
\gg&=\{\chi\in U; \chi(k\ell)=\chi(k)\ell+
(\overline{\mathscr{R}}^\al\ra k)
(\overline{\mathscr{R}}_\al\ra\chi)(\ell), \mbox{ for all } k,\ell\in K\}\nonumber\\
&=\{\chi\in U; \chi(k\ell)=\chi(k)\ell+
{\mathscr{R}}^{-1}(k\otimes \ell_{(1)}S(\ell_{(3)}))\,\chi(\ell_{(2)}), \mbox{ for all } k,\ell\in K\}\label{kerepsi2}\\
&=\{\chi\in U; \chi(1_K)=0 \mbox{ and } \chi(\ker{\varepsilon}\underline\cdot\!\;\ker{\varepsilon})=0\} \label{kerepsi2p}
\end{align}
where $ \ker{\varepsilon}\underline\cdot\!\;\ker{\varepsilon}$ is the
left ideal in $K$ considered in \cite[Prop. 4.8]{Gomez-Majid},
with, for all $\ell,k\in K$,
$\ell\underline{\cdot}k:=\ell_{(1)}k_{(2)}{\mathscr{R}}(k_{(1)}S(k_{(3)})\otimes
S(\ell_{(2)})
=k_{(2)}\ell_{(2)}{\mathscr{R}}(k_{(1)}\otimes \ell_{(1)}S(\ell_{(3)}))
$ 
(see also Proposition 4.7 for this last equality, wherein 
the notation used is
$\mbox{\bf{r}}_2(\ell)={\mathscr{R}}_\alpha(\ell){\mathscr{R}}^\alpha$).
The equality \eqref{kerepsi2} follows from duality, for all
$\zeta\in U$,
$$(\zeta\ra
\chi)(\ell)=\le\zeta_{(1)}\chi S(\zeta_{(2)})\,,\,\ell\re=\le\zeta_{(1)}\otimes
S(\zeta_{(2)})\otimes\chi\,,\,\ell_{(1)}\otimes\ell_{(3)}\otimes\ell_{(2)}\re=\le\zeta\otimes\chi,\ell_{(1)}S(\ell_{(3)})\otimes
\ell_{(2)}\re.
$$
We prove the equality in \eqref{kerepsi2p} by proving the inclusions
$\subseteq$ and $\supseteq$. If $\chi\in \gg$ then
we have already seen that $\chi(1)=0$; we have $\chi(\ker{\varepsilon}\underline\cdot\!\;\ker{\varepsilon})=0$ because 
$\ell\mapsto \ell_{(1)}S(\ell_{(3)})\otimes\ell_{(2)}=:\ell_{-1}\otimes\ell_0$ is
the adjoint left coaction of $K$ on $K$,
so that:
\begin{equation*}
  \begin{split}
  \chi(\ell\underline\cdot\!\;k)&=
\chi(k_{(2)}\ell_{0}){\mathscr{R}}(k_{(1)}\otimes
\ell_{-1})\\ &=\chi(k_{(2)})\varepsilon(\ell_{(0)}){\mathscr{R}}(k_{(1)}\otimes
\ell_{-1})+{\mathscr{R}^{-1}}(k_{(2)}\otimes \ell_{-1})\chi(\ell_{0}){\mathscr{R}}(k_{(1)}\otimes
\ell_{-2})\\
&=
\chi(k_{(2)}){\mathscr{R}}(k_{(1)}\otimes
\varepsilon(\ell)1)+\chi(\ell_{0})\,\mathscr{R}\mathscr{R}^{-1}(k\otimes \ell_{-1})\\
&=\chi(k)\varepsilon(\ell)+\chi(\ell)\varepsilon(k)~,
\end{split}
\end{equation*}
where we  used \eqref{kerepsi2} in the second line. This shows the
inclusion $\subseteq$ in line \eqref{kerepsi2p}. For the other inclusion, if 
$\chi(1)=0$ and
$\chi(\ker{\varepsilon}\underline\cdot\!\;\ker{\varepsilon})=0$ then
from
$\chi\big((\ell-\varepsilon(\ell)1)\,\underline{\cdot}\,(k-\varepsilon(k)1)\big)=0$
we have ${\mathscr{R}}(k_{(1)}\otimes
\ell_{-1})\chi(k_{(2)}\ell_0)=\chi(\ell \,\underline{\cdot}\,
k)=\chi(k)\varepsilon(\ell)+\varepsilon(k)\chi(\ell)$ and equivalently
$$k_{(1)}\otimes \ell_{-2}\,{\mathscr{R}}(k_{(2)}\otimes
\ell_{-1})\chi(k_{(3)}\ell_0)=k_{(1)}\otimes \ell_{-1}\big(\chi(k_{(2)})\varepsilon(\ell_0)+\varepsilon(k_{(2)})\chi(\ell_0)\big)~.$$
Applying ${\mathscr{R}}^{-1}$ to this expression we obtain
\begin{equation*}
  \begin{split}
    \chi(k\ell)\,&=\,{\mathscr{R}}^{-1}(k_{(1)}\otimes\ell_{-1})\chi(k_{(2)})\varepsilon(\ell_0)\,+\,{\mathscr{R}}^{-1}(k\otimes\ell_{(1)}S(\ell_{(3)}))\chi(\ell_{(2)}) \\
&=\,{\mathscr{R}}^{-1}(k_{(1)}\otimes\varepsilon(\ell)1)\chi(k_{(2)})
\,+\,{\overline{\mathscr{R}}}^{\al}(k){\overline{\mathscr{R}}}_{\al}
(\ell_{(1)}S(\ell_{(3)}))\chi(\ell_{(2)}) \\
&=\chi(k)\varepsilon(\ell)\,+ \,{\overline{\mathscr{R}}}^{\al}(k)
\le {{\overline{\mathscr{R}}}_{\al}}^{}_{(1)}\otimes
S({{\overline{\mathscr{R}}}_{\al}}^{}_{(2)})\otimes\chi\,,\,\ell_{(1)}\otimes\ell_{(3)}\otimes\ell_{(2)}\re\\
  &
=\chi(k)\varepsilon(\ell)\,+ \,{\overline{\mathscr{R}}}^{\al}(k)
({{\overline{\mathscr{R}}}_{\al}}\ra \chi)(\ell)~,
  \end{split}
  \end{equation*}
  that shows $\chi \in \gg$.
\end{remark}

\begin{remark} 
A slight variation of the construction in Example \ref{excotriang}
allows to consider $A=K$ an infinite dimensional Hopf algebra over
a field $\Bbbk$ with cotriangular structure $\mathscr{R}:K\otimes K\to
\Bbbk$. This is obtained by rewriting the $U^\op$ and $U$-actions in
terms of $K$-coactions, cf. Example \ref{cotriangularK}. 
In this infinite dimensional case we consider the rational morphisms
  ${\rm{HOM}}_\Bbbk(K,K)$, see e.g. \cite{Ulbrich}, i.e.,  those $\Bbbk$-linear maps
  $L:K\to K$ such that the left and right coadjoint actions
$\Delta_\Ls:{\rm{HOM}}_\Bbbk(K,K)\to K\otimes {\rm{HOM}}_\Bbbk(K,K)$ and
$\Delta_\Rs:{\rm{HOM}}_\Bbbk(K,K)\to {\rm{HOM}}_\Bbbk(K,K)\otimes K$ are well defined
by, for all
$k\in K$,
$L_{-1}\otimes L_0(k):=L(k_{(2)})^{}_{(1)}S^{-1}(k_{(1)})\otimes
L(k_{(2)})^{}_{(2)}$ and
$L_{0}(k)\otimes L_{1}:=L(k_{(1)})^{}_{(1)}\otimes
L(k_{(1)})^{}_{(2)}S(k_{(2)})$, cf. \eqref{deltaL}. 
Then, setting $U:=K^\circ$ (the Hopf dual of
$K$) and recalling the definition of $U$-adjoint action $\ra$ in
\eqref{adjact} and the associated $U^\op$ one just before
\eqref{deltaL},  we have, for all  $\zeta\otimes\xi\in U^\op\otimes U$,
$(\zeta\otimes\xi)\ra L=
\zeta(L_{-1})L_0\xi(L_1)$.
 Recall that $K$ is a right $K^\cop\otimes K$-comodule where $K^\cop\otimes K$ has
cotriangular structure
 ${\RR}=(\mathscr{R}^{-1}\otimes\mathscr{R})\circ (\id\otimes 
\rm{flip}\otimes\id)$. Braided
derivations $u\in \Der_\RR(K)$ are then defined as in  \eqref{brderK}
but with $u\in {\rm{HOM}}_\Bbbk(K,K)$ and
where, for all $k\in K$, $\oR^\alpha\ra k\, \oR_\alpha \RA
u=\RR^{-1}\ra (k\otimes u)$
stands now for
$k_{(2)}u_0\!\:\RR^{-1}(k_{(1)}\otimes k_{(3)}\,\otimes\, u_{-1}\otimes
u_1)$, (cf. \eqref{cotquasi}).
They span a $K$-bicovariant bimodule so that they are a free left
$K$-module as in \eqref{XXinv}. The right $K$-comodule of left invariant vector fields
${}_{\rm{inv}}\Vect(K)$ is isomorphic to the right $K$-comodule $\gg$ as defined
in \eqref{gforbic} where
$\overline{\mathscr{R}}^\al\otimes\overline{\mathscr{R}}_\al\ra\chi$
is the element in $U\otimes U$ defined by ${\mathscr{R}}^{-1}(\;\mbox{-}\,\otimes
\chi_1)\otimes\chi_0$. (That, for all $k\in K$, ${\mathscr{R}}^{-1}(\,\mbox{-}\,\otimes k)\in
U$ follows from  cotriangularity of
$\mathscr{R}$).
We then proceed as after equation \eqref{gforbic} and, if $U$
separates the elements of $K$, we conclude that $S(\gg)$ defines a
quantum Lie algebra. The bicovariant differential calculus is as in \eqref{defofbicdd}. 
\end{remark}

\begin{example}
Let $G$ be a Lie group and $M$ a $G$-manifold as in Example
\ref{ex:bundle}. Recall the Drinfeld twist deformation of
vector bundles on $M$ of Example \ref{ex:bundle}.
Drinfeld twist deformations of the differential and Cartan
calculus on $M$  have
been considered in \cite{NCG2}  and in more generality in
\cite{Weber}. They exemplify the constructions presented in
this section. 
\endex\end{example}

\section{Right connections, left connections, curvature and torsion}
\label{sec:connections}

We review connections on modules $\Gamma\in \cat$ considered as right 
$A$-modules or as left $A$-modules and study their 
extensions to $\Cat$, where
 as usual $H$ is a triangular Hopf algebra, $A$ a braided
  commutative $H$-module algebra and
  $\bigl(\OM,\wedge,\dd\bigr)$ the associated
  differential calculus constructed in Section \ref{sectCC}. This first
  explanatory part can be found e.g. in \cite[\S 4]{BM}, the
  braided commutativity aspects being in \cite{AS}.
  Covariant derivatives along vector fields are then introduced, their
braided bracket with inner derivatives is an inner derivative
generalizing the braided Cartan identity
$[\LLL_u,\ii_v]=\ii_{[u,v]}$ to $[\stDu,\ii_v]=\ii_{[u,v]}$. The latter implies that the
curvature tensor, defined as the square of the connection, can be
equivalently defined via the commutator of covariant derivatives along
vector fields. This extends to the braided commutative geometry
setting the usual two equivalent definitions of curvature. Similarly for
the torsion tensor. 

Right (left) connections on modules in $\cat$ are
 more general connections than the bimodule connections considered in \cite{D-VM96}. As
 shown in the last subsection, using the braided commutativity
 property of $A$ and of $\Gamma$ they can be summed to give connections on
tensor product modules.

\subsection{Connections and Cartan formula}
  Let $H$ be a triangular Hopf algebra, $A$ be a braided
  commutative $H$-module algebra and
  $\bigl(\OM,\wedge,\dd\bigr)$ the associated
 $H$-equivariant differential calculus constructed in Section \ref{sectCC}.  
\begin{definition}\label{defi:connection} 
A right connection on a module $\Gamma$ in $\cat$ 
is a $\Bbbk$-linear map 
\begin{equation} \dst : \Gamma \rightarrow
\Gamma \otimes_A \Omega(A)
\end{equation} 
in $\hom_\Bbbk(\Gamma,\Gamma\oA \Omega(A))$,
which satisfies the Leibniz rule, for all
$s\in\Gamma$, $a\in A$, 
\begin{equation} 
\label{prop_Delta} 
\dst(s a) =
\dst(s) \:\!a+ s\otimes_A\dd a ~.
\end{equation}  
A left connection on $\Gamma$ is a $\Bbbk$-linear map
\begin{equation} \std : \Gamma \rightarrow
\Omega(A)\otimes_A\Gamma
\end{equation} 
in ${}_\Bbbk\hom(\Gamma,\Omega(A)\oA\Gamma)$,
which satisfies the Leibniz rule,
\begin{equation}  \label{stdasleib}
\std(a s) =
\dd a\otimes_A s+ a\std(s)~.
\end{equation}
\end{definition}

We denote by $\Con_A(\Gamma)$ and ${}_A\Con(\Gamma)$ the set of all
right, respectively left connections. Notice that in the definition of
$\Con_A(\Gamma)$
(respectively ${}_A\Con(\Gamma)$),  no compatibility condition with the left
(right) $A$-module structure of $\Gamma$ is required.

Given any connection $\dst \in \Con_A(\Gamma)$ and any right $A$-linear map
$L\in\hom_A(\Gamma,\Gamma\otimes_A\Omega(A))$, the sum $\dst+ L$ is a
connection  in $\Con_A(\Gamma)$.
The  action $\dst\mapsto \dst +L$ is free and transitive and
hence $\Con_A(\Gamma)$ is an affine space over the module $\hom_A(\Gamma,\Gamma\otimes_A\Omega(A))$ in $\cat$.

Similarly, the difference $\std-\std'$ of two left
connections is a left $A$-linear map, and ${}_A\Con(\Gamma)$ is an
affine space over the module ${}_A\hom(\Gamma,\Gamma\otimes_A\Omega(A))$ in $\cat$.
Notice that left connections, as left $A$-linear maps, acts from the
right, cf. \eqref{copadjact} and \eqref{copadjact2}. Their properties
become more intuitive by evaluating them on
the right of elements $s\in \Gamma$, hence writing  $(s)\std$ rather than $\std(s)$.
\\

We can act with the 
$H$-adjoint action $\RA$ defined in \eqref{adjact} on $\dst\in \Con_A(\Gamma)\subset \hom_\Bbbk(\Gamma,\Gamma\otimes_A\Omega(A))$, for all $\hxi\in H$,
\begin{equation}
\label{eqn:conadjoint}
 \hxi\RA\dst := \hxi_1\ra\,\circ \dst \circ S(\hxi_2)\ra\,~.
\end{equation}
This $\Bbbk$-linear map $\hxi\RA\dst\in
\hom_\Bbbk(\Gamma,\Gamma\otimes_A\Omega(A))$ is easily seen to satisfy
(cf. \cite[\S 6.2]{AS}),
for all $s\in \Gamma$ and $a\in A$,
\begin{equation}\label{hdsa} (\hxi\RA\dst\,) (sa)=
(\hxi\RA \dst\,) (s)\;\! a\, +\,s\otimes_A  \varepsilon(\hxi)\dif a~.
\end{equation}
In particular we see that if $\varepsilon(\hxi)=0$ then
$\hxi\RA \dst \in\hom_A(\Gamma,\Gamma\otimes_A \Omega(A))$, while if
$\varepsilon(\hxi)=1$ then $\hxi\RA\dst\in\Con_A(\Gamma)$.
Similarly for $\std\in{}_A\Con(\Gamma)$. 
Using this action and the braided commutativity of the $A$-bimodule $\Gamma$,  a right connection
$\dst$ on $\Gamma$ is shown to be also a braided left connection,
cf. \cite[Prop. 6.8]{AS}, and similarly 
a left connection $\std$ on $\Gamma$ is also a braided right
connection, for all $a\in A, s\in \Gamma$, 
\begin{equation}\label{braidedconn}
\begin{split}
~~~~~~~~~~~~~
\end{split}
\begin{split}
\dst(a s)&=(\oR^\al\ra a) (\oR_\al\ra 
\dst)(s)+\oR^\al\ra s \otimes_A \oR_\al\ra \dif a~,\\
\std(s a)&= (\oR^\al\RA^\cop
\std)(s) (\oR_\al\ra a)+\oR^\al\ra \dif a\otimes_A\oR_\al\ra s ~.
\end{split}
\end{equation}
If $\dst$ is $H$-equivariant we have $\dst(a s)= a
\dst(s)+\oR^\al\ra s\otimes_A \oR_\al\ra\dif a
$, and thus recover the notion of bimodule connection
studied in \cite{D-VM96}.
\\

The $H$-module algebra homomorphism (injection) $A\rightarrow \OM$
allows to associate to modules in $\cat$ modules in $\Cat$ via the
change of base ring
$\Gamma\to \OM\otimes_A\Gamma$.
Indeed, given $\Gamma$ in $\cat$, the bimodule in $\cat$ 
\begin{equation}\label{defOMM}
\OM\otimes_A \Gamma
\end{equation}
is naturally a left $\OM$-module
by defining, for all $\theta, \theta' \in \OM$, $s\in \Gamma$, $\theta' \wedge (\theta \otimes_\AA s):=
(\theta' \wedge \theta) \otimes_\AA s$. 
Since $\OM$ is a graded braided commutative
$H$-module algebra, $\OM\otimes_A \Gamma=\bigoplus_{n\in\mathbb{N}\,}
\Omega^{n\!}(A)\otimes_A\Gamma$  becomes a graded braided
commutative $\OM$-bimodule by defining 
$(\theta \otimes_\AA s)\wedge\theta':=
\theta\wedge \oR^\al\ra^\cop \theta' \,\otimes_\AA \oR_\al \ra s$. Hence 
$\OM\otimes_A \Gamma$
is a module in $\Cat$. We notice that, for any $\Gamma,\Gamma'$ in
$\cat$, the canonical isomorphism in $\Cat$, 
$$(\OM\otimes_A \Gamma)\otimes_\OM (\OM\otimes_A \Gamma')\,\simeq\,\OM\otimes_A
(\Gamma\otimes_A\Gamma')~,$$
given by $(\theta\oA s)\otimes_\OM(\theta'\oA s')=(\theta\oA
s)\otimes_\OM\theta'\wedge(1_\OM\oA
s')\longmapsto
\theta\wedge\oR^\al\ra^{\cop\,}\theta'\oA (\oR_\al\ra s\oA s')$
 implies that the 
association $\Gamma\to \OM\otimes_A\Gamma$
defines a strict monoidal functor from $\cat$ to $\Cat$, and
hence from $\catfp$ to $\Catfp$.

More generally,  given $W$ in $\Cat$  (and hence in $\cat$)
and $\Gamma$ in $\cat$,  we have
$W\otimes_A \Gamma$ in $\cat$ and we structure it as a module in $\Cat$ with the obvious left 
$\OM$-module structure and with  right $\OM$-module structure
determined by 
requiring  $W\otimes_A\Gamma$ to be graded symmetric. Notice that
if $\Sigma \in \cat$, the associativity $(W\otimes_A\Gamma)\otimes_A\Sigma =
W\otimes_A(\Gamma\otimes_A\Sigma)$ in $\cat$ lifts to $\Cat$.

Analogously, we define the module $ \Gamma\otimes_A \OM$ in
$\Cat$. The  $\OM$-bimodule actions read, for
all $s\in \Gamma$, $\theta, \theta'\in \OM$, 
\begin{equation}\label{sotwsot}
(s\otimes_\AA \theta)\wedge \theta':=s\otimes_A (\theta\wedge
\theta')~,~~\theta'\wedge (s\otimes_A\theta ):=(\oR^\al \ra s)
\otimes_A(\oR_\al\ra^\cop\theta')\wedge \theta~.
\end{equation}

Every left $A$-linear map $\tilde L\in {}_A\hom(\Gamma,\Omega^{n}(A)\otimes_A
\Gamma)$, uniquely extends to a graded left $\OM$-linear  map 
in  ${}_\OM\hom(\OM\otimes_A\Gamma,\Omega^{\bullet +n}(A)\otimes_A
\Gamma)$, that we still denote $\tilde L$ and with degree $|
\tilde L|=n$; it is given by
$\tilde L(\theta\otimes_A s)=(-1)^{|\tilde
  L||\theta|}\theta\wedge\tilde L(s)$ for all $\theta\in \OM$ of
homogeneous degree $|\theta|$ and $s\in
\Gamma$. This provides an isomorphism 
${}_A\hom(\Gamma,\Omega(A)\otimes_A
\Gamma)\simeq {}_\OM\hom(\OM\otimes_A\Gamma,\OMPu(A)\otimes_A
\Gamma$) of modules in $\cat$. Similarly, 
$\hom_A(\Gamma,\Gamma\otimes_A
\Omega(A))\simeq \hom_\OM(\Gamma\otimes_\OM\OM,\Gamma\otimes_A
\OMPu(A))$ as modules in $\cat$.

Correspondingly, right and left connections on $\Gamma$ uniquely extend as $\Bbbk$-linear maps of
degree one on the modules $\Gamma\otimes_\AA \OM$ and $\OM\otimes_\AA
\Gamma$.
\begin{lemma} 
The connections $\dst\in\Con_A(\Gamma)$ and $\std\in{}_A\Con(\Gamma)$ extend to the graded maps
$\Dst\in \hom_\Bbbk(\Gamma\otimes_\AA \OM, \Gamma\otimes_\AA
\Omega^{\bullet +1}(A))$ and $\stD\in {}_\Bbbk\hom(\OM\otimes_\AA \Gamma,
\Omega^{\bullet +1}(A)\otimes_\AA\Gamma)$ well-defined by
\begin{equation}
    \begin{split}\label{dOMdef}
\Dst: \Gamma\otimes_\AA \OM &\longrightarrow  \Gamma\otimes_\AA \Omega^{\bullet+1\!}(A)~~,\\
s\otimes_\AA 
\theta~ &\longmapsto \Dst(s\otimes_\AA 
\theta):=\dst(s)\wedge\theta+s\otimes_\AA \dd\theta~,~~~~~
\end{split}
\end{equation}
and, for all $k\in \mathbb{N}$,\begin{equation}
    \begin{split}\label{defofdnabla}
\stD:  \Omega^{k}(A)\otimes_\AA \Gamma &\longrightarrow  \Omega^{k+1\!}(A)
\otimes_\AA\Gamma~,\\
\theta\otimes_\AA 
s~ &\longmapsto \stD(\theta\otimes_\AA 
s):=\dd\theta\otimes_\AA s+(-1)^k\theta \wedge {\std}( s)~.
\end{split}
\end{equation}
 More generally, for any $h\in H$, the maps $h\RA \dst$ and
$h\RA^{\cop}\std$ extend to the well-defined internal morphisms in ${}^H\!\MMM$,
\begin{equation} \label{hDD}
\begin{split}
\end{split}
\begin{split}\dd^{\phantom{J_J}}_{_{}h_{}\RA}\ddst(s\otimes_\AA
\theta)&:=(h\RA\dst)(s)\wedge\theta+s\otimes_\AA
\varepsilon(h)\dd\theta~,~~\\[.2em]
\dd^{\phantom{J_J}}_{_{}h_{}\RA^{\!\:\cop\!\:}}\stdd(\theta\otimes_\AA 
s)&:=\varepsilon(h)\dd\theta\otimes_\AA s+(-1)^k\theta \wedge
{h\RA^\cop \std}( s)~. 
\end{split}
\end{equation}
The $H$-action on $\Dst$ and $\stD$ reads, for all $h'\in H$, 
$ h'\RA \Dst=\dd^{\phantom{J_i}}_{h'\RA\!}\ddst,$
$h'\RA^\cop \stD=\dd^{\phantom{J_J}}_{_{}h'_{}\RA^{\!\:\cop}}\stdd$,
and we have the Leibniz rule, for all $\varsigma\in
\Gamma\otimes_\AA\Omega^k(A)$, $\vartheta\in\OM$
and for all $\sigma\in
\OM\otimes_\AA\Gamma$ and
$\theta\in \Omega^k(A)$, 
\begin{equation}\label{dleibOM}
\Dst(\varsigma\wedge \vartheta)= \Dst_{}\varsigma\wedge 
\vartheta+(-1)^k\varsigma\wedge \dd\vartheta~~,~~~
\stD(\theta\wedge \sigma)= \dd\theta\wedge 
\sigma+(-1)^k\theta\wedge \stD\sigma~.
\end{equation}
\begin{proof}
We give the proof for left connections, the other case is similar. The definitions in \eqref{defofdnabla} and \eqref{hDD} are well-defined because  independent from the
representative chosen for the balanced tensor product over $A$.
We now show that
$\stD\in {}_\Bbbk\hom(\OM\otimes_\AA\Gamma, \Omega^{\bullet
  +1}(A)\otimes_\AA\Gamma)$.
Let  $V,W\in \catm$ and $\pi_{V,W}: V\otimes W\to V\oA W$ be the projection to the
balanced tensor product; consider the $\Bbbk$-linear map
\begin{equation}\label{ddRididd} 
\begin{split}
\phantom{(\id_\Gamma\otimes_{\!\!\:A}\wedge)\!\circ\! \pi^{}_{\Gamma\oA\Omega,\OM}\!\circ\! (\dst\otimes_\RR\id_\OM)+\pi^{}_{\Gamma,\Omega^{\bullet +1}(A)}\!\circ\!(\id_\Gamma\otimes\dd)} &\phantom{: \Gamma\otimes
\OM\longrightarrow  \Gamma\otimes_\AA
\Omega^{\bullet +1}(A)}\\[-1.6em]
\pi^{}_{\Omega^{k +1}\!,\Gamma}\!\circ\!(\dd\!\:\!\otimes\!\:\! \id_\Gamma)+(-1)^{\!\!\:k}\!\;\!(\wedge\otimes_{\!\!\:A}\id_\Gamma)\!\circ\!\pi^{}_{\Omega^{k}\!,\Omega(A)\otimes_{\!\!\:A}\Gamma}\!\circ\!(\id_{\Omega^{k}\!\!\:(A)}\tilde\otimes_\RR \std) \!&:\!\:\! \Omega^{k}(A)\otimes\Gamma\longrightarrow \Omega^{k +1}\!\!\:(A)\otimes_{\!\!\:\AA}\Gamma\,
\end{split}
\end{equation}
where 
$\wedge\oA\id_\Gamma$ is a tensor product
of morphisms in $\catm$ while
 $\dd\otimes \id_\Gamma$ is a  tensor product
of morphisms in ${}^{H}\MMM$
and 
 $\id_{\Omega^{k}\!\!\:(A)}\tilde\otimes_\RR
\std$ is a tensor product
of internal morphisms in ${}^{H}\MMM$.
The maps in \eqref{ddRididd} are sums of compositions of $\Bbbk$-linear maps carrying
the $\ra^\cop$ adjoint action  (recall that a morphism in
${}^{H}\MMM$ or  $\catm$ can be seen as an internal morphism carrying
trivial $\ra^\cop$ adjoint action, cf.~the last paragraph of
Section \ref{subsecMoHA}). Since the composition of internal morphisms is an
internal morphism (cf. equation \eqref{eqcompadjcop}), the map in
\eqref{ddRididd} is therefore  in
${}_\Bbbk\hom(\OM\otimes\Gamma, \Omega^{\bullet
  +1}(A)\otimes_\AA\Gamma)$.
By evaluating these internal morphisms on elements  $\theta\otimes s\in  \Omega^{k}\otimes\Gamma$
we see that they equal
$\stD\circ\pi_{\Omega^{k}(A),\Gamma}$, thus showing that
$\stD\in {}_\Bbbk\hom(\OM\otimes_\AA\Gamma, \Omega^{\bullet
  +1}(A)\otimes_\AA\Gamma)$.  Similarly we prove that for all $h\in H$,
$\dd^{\phantom{J_J}}_{_{}h_{}\RA^{\!\:\cop}}\stdd\in {}_\Bbbk\hom(\OM\otimes_\AA\Gamma, \Omega^{\bullet
  +1}(A)\otimes_\AA\Gamma)$.

The $H$-action on the connection follows observing that the only map in \eqref{ddRididd} that
carries a nontrivial $H$-action is  $\std$.

The Leibniz rule is easily proven, indeed because of $\Bbbk$-linearity it is enough to consider
$\sigma=\eta\otimes_\AA s$  with $\eta\in \Omega^k(A)$. 
\end{proof}
\end{lemma}

Vice versa,  internal morphism
$\Dst \in \hom_\Bbbk(\Gamma\otimes_\AA \OM, \Gamma\otimes_\AA
\Omega^{\bullet +1}(A))$,
$\stD\in {}_\Bbbk\hom(\OM\otimes_\AA\Gamma, \Omega^{\bullet
  +1}(A)\otimes_\AA\Gamma)$ that satisfy the
Leibniz rule \eqref{dleibOM} are uniquely defined by connections
$\dst$, $\std$, respectively. In view of this bijection we also call  $\Dst$ and
$\stD$ connections.
\\

We canonically extend the inner derivative $\ii:\Vect(A)\to
\hom_A(\OM,\OMmu(A))$ to $\ii:\Vect(A)\to
\hom_A(\OM\otimes_A \Gamma,\OMmu\oA \Gamma)$ by defining, for all $u\in
\Vect(A)$, 
\begin{equation}\label{defcontraction}
\ii_u: \OM \otimes_\AA\Gamma\to
\OMmu(A)\otimes_\AA\Gamma~,~~\theta\otimes_A s\mapsto\ii_u(\theta\otimes_\AA s)=\ii_u(\theta)\otimes_\AA s~.
\end{equation}
The {\it covariant derivative 
 along a vector field} $u\in
\Vect(A)$ of a left connection is the
$\Bbbk$-linear operator of zero degree $\stD{}_{\mbox{$\!{}^{}_u$}}: \OM\otimes_\AA\Gamma\to \OM\otimes_\AA\Gamma$ defined by 
\begin{equation}\label{defCD}
\stD{}_{\mbox{$\!{}^{}_u$}}:=\ii_u\circ \stD+\stD
\circ \ii_u~,
\end{equation}
in particular, on $\Gamma$ we have $\stD{}_{\mbox{$\!{}^{}_u$}}=\ii_u\circ \std$
that as usual denote by $\std{}_u$,
hence $\std{}_u:=\ii_u\circ \std$.
Notice however that $\std{}_u$ is the composition of a
$\Bbbk$-linear map $\std$ acting from the right and a $\Bbbk$-linear
map $\ii_u$ acting from the left. From \eqref{defofdnabla} we have
that the covariant derivative satisfies, for all $u,v\in\VV,~a\in A, 
\;s,t \in\Gamma$, 
\begin{equation}\label{PROPDu}
\begin{split}
\std{}_{u+v} s&=\std{}_u s+\std{}_v s~,\\ \std{}_{au} s&=a\std{}_{u}s~,\\
 \std{}_u (s+t)&=\std{}_u (s)+\std{}_u (t)~,
\end{split}
\end{equation}
and 
\begin{equation}\label{LEIBDu}
{}~~~~~~~~~~~~~~~~~~
\std{}_u(a s) 
=\LLL_u(a) s+ (\oR^\al\ra a)\std{}_{\oR^\al\ra u} s~.
\end{equation}

\begin{remark} \label{remcd}
We also term covariant derivative a map $\std^\cd:\Vect(A)\times \Gamma\to \Gamma$ that satisfies
\eqref{PROPDu} and \eqref{LEIBDu}; equivalently, since \eqref{PROPDu} and
\eqref{LEIBDu} imply $\Bbbk$-bilinearity, we term covariant derivative
a left $A$-linear map 
$\std^{\cd}:\Vect(A)\otimes\Gamma\to \Gamma$ that satisfies the Leibniz
rule  $\std^\cd(u\otimes a s) 
=\LLL_u(a) s+ (\oR^\al\ra a)\std^\cd(\oR^\al\ra u\otimes
s)$. This latter more elegantly reads 
\begin{equation}\label{leibnablacd}
\std^\cd(u\otimes a s) 
=\LLL_u(a) s+ \std^\cd(ua\otimes s)~.
\end{equation} 
Since $\std\in
{}_A\Con(\Gamma)\subset {}_\Bbbk\hom(\Gamma,\Omega(A)\oA\Gamma)$, we also require $\std^\cd\in{}_\Bbbk\hom(\VV\otimes
\Gamma,\Gamma)$. A {\it{covariant derivative}} is therefore a map
 $\std^\cd\in{}_A\hom(\VV\otimes
\Gamma,\Gamma)$ satisfying the Leibniz rule \eqref{leibnablacd}.
Notice that ${}_A\hom(\VV\otimes
\Gamma,\Gamma)$ is a module in $\catm$, not in $\cat$, indeed $\VV\otimes
\Gamma$ is  in $\catm$ not in $\cat$.

Let  $\pi: \VV\otimes \Omega(A)\to \VV\otimes_A \Omega(A)$ be the
canonical projection to the balanced tensor product. To any connection  $\std\in
{}_A\Con(\Gamma)$ we associate a covariant derivative  via the map
 $$\std\longmapsto
\std^\flat:=(\le~,~\re\circ \pi\otimes_A \id_\Gamma)\circ
(\id_{\VV}\otimes\std)~,~~ \std^\flat(u\otimes s)=\ii_u\circ\std(s)~.$$
If the module $\VV$ is finitely generated and projective, this map is a bijection with inverse
$\std^\cd\mapsto (\std^\cd)^\sharp:=(\id_\Omega\otimes_A\std^\cd)\circ
(\coev\circ \eta \otimes \id_\Gamma)$,
$(\std^\cd)^\sharp(s)=\omega^i\otimes_A\std^\cd(e_i\otimes
s)$; here as usual $\Bbbk\otimes \Gamma\simeq\Gamma$,
$\eta: \Bbbk\to A$ is the unit $1_\Bbbk\mapsto
\eta(1_\Bbbk)=1_A$, 
$\{e_i,\omega^i : i=1\ldots, n\}$
 is a dual basis of $\VV$ and $\coev(1_A)=\omega^i\otimes_A e_i$.

It is easy to see that covariant derivatives form an affine
space over the module ${}_A\hom(\Vect(A)\oA\Gamma,\Gamma)$ in $\cat$. 
The bijection $\flat=\sharp^{-1}$ is compatible with
the affine structures of the space of connections ${}_A\Con(\Gamma)$
and the affine space of covariant derivatives; thus it lifts to an 
isomorphism of affine spaces over the isomorphic modules
${}_A\hom(\Gamma,\Omega\oA\Gamma)$ and
${}_A\hom(\Vect(A)\oA\Gamma,\Gamma)$ in $\cat$ (cf. Theorem \ref{Thmfgp123}).
\end{remark}

\begin{lemma}
The covariant derivative $\stDu$ along a vector field $u\in \VV$ of a left connection satisfies the braided Leibniz rule, for all
$\theta\in \OM$, $\sigma\in\OM\otimes_\AA\Gamma$
\begin{equation}
\stDu (\theta\wedge\sigma)=\LLL_u(\theta)\wedge\sigma+(\oR^\al\ra^\cop\theta)\wedge \stDalu(\sigma)~.
\end{equation}
\begin{proof} Because of $\Bbbk$-linearity it is enough to consider a form of
  homogeneous degree $\theta\in \Omega^k(A)$. We apply the
  definition and use the Leibniz rule for $\stD$ and
  the braided one for the inner derivative,
$\ii_u(\theta\wedge\sigma)=\ii_u(\theta)\wedge\sigma+(-1)^k(\oR^\al\ra^\cop\theta)\wedge
\ii_{\oR_\al\ra u}\sigma$,
that immediately follows from \eqref{lieinnerder},
thus obtaining
\begin{eqnarray*}
&&\!\!\!\!\!\!\!\!\!\!\!\!\!\!\! \!\! \stDu\!(\theta\wedge\,\sigma)=\nn\\[.2em]
&=&\!\!\!(\ii_u\circ \stD+\stD
\circ \ii_u)(\theta\wedge\sigma)\\[.2em]
&=&\!\!\!\ii_u\big(\dd\theta\wedge\sigma+(-1)^k\theta\wedge
\stD\sigma\big)+
\stD\big(\ii_u(\theta)\wedge\sigma+(-1)^k(\oR^\al\ra^\cop\theta)\wedge
\ii_{\oR_\al\ra u}\sigma\big)\\[.2em]
&=&\!\!\!\ii_u(\dd\theta)\wedge\sigma-(-1)^k
(\oR^\al\ra\dd\theta)\wedge \ii_{\oR_\al\ra u}\sigma+(-1)^{k\,}\ii_u(\theta)\wedge
 \stD\sigma+(\oR^\al\ra^\cop\theta)\wedge \ii_{\oR_\al\ra u} \stD\sigma
\\[.2em]
&&\!\!\!+\,\dd(\ii_u\theta)\wedge\sigma-(-1)^k\:\!\ii^{}_u(\theta)\wedge
\stD\sigma+(-1)^k\dd(\oR^\al\ra^\cop\theta)\wedge
\ii_{\oR_\al\ra u }\sigma+(\oR^\al\ra^\cop\theta)\wedge
\stD\ii_{\oR_\al\ra u}\sigma
\\[.2em]
&=&\!\!\!
\LLL_u(\theta)\wedge\sigma+(\oR^\al\ra^\cop\theta)\wedge   
\std{}_{\oR_\al\ra u}(\sigma)
\end{eqnarray*}
where we  used the identity
$\dd(\oR^\al\ra^\cop\theta)\otimes\oR_\al=\oR^\al\ra
(\dd\theta)\otimes\oR_\al$ due to  $H$-equivariance of the exterior derivative. 
\end{proof}
\end{lemma}
\begin{theorem}\label{THM44}
For all vector fields $u, v\in \Vect(A)$ the covariant derivative 
$\stDu: \OM\otimes_A\Gamma\to \OM\otimes_A\Gamma$ and the inner derivative
$\ii_v: \OM\otimes_A\Gamma\to\OMmu \otimes_A\Gamma$, satisfy the
braided Cartan relation
\begin{equation}
\stDu\circ \ii_v -\ii_{\oR^\al\ra v} \circ \stDalu=\ii_{[u,v]}~,\label{cartanCalculuscovder}
\end{equation}
equivalently, 
\begin{equation}\label{CRID}
 \ii_u \circ \stDv -\stDalpv\circ \ii_{\oR_\al\ra u} =\ii_{[u,v]}~.
\end{equation}
 \begin{proof}
Because of $\Bbbk$-linearity, in order to  prove the first relation it is enough to
evaluate it on  $\theta\otimes_A s$, with 
$\theta\in \OM$ and $s\in \Gamma$.
We then compute
\begin{eqnarray*}
\big(\stDu\circ \ii_v -\ii^{}_{\oR^\al\ra v}
  &&\!\!\!\!\!\!\!\!\!\!\!\!\!\!\!\!\!\circ\, \stDalu{} \,\big)
(\theta \otimes_\AA s)\,=\\[.4em]
&=&\!\!
\stDu
\big(\ii_v\theta\otimes_\AA s\big) 
-\ii_{\oR^\al\ra v}\big(\LLL_{\oR_\al\ra u}\theta\otimes_\AA s+
(\oR^\beta\ra^\cop\theta)\otimes_\AA \std{}_{\oR_\be\oR_\al\ra u}  s\big)
\\[.4em]
&=&\!\!
\LLL_u\circ  \ii_v(\theta)\otimes_\AA s
+\oR^\al\ra (\ii_v\theta)\otimes_\AA  \std{}_{\!\!^{\phantom{.}}_{\oR_\al\ra u}}  s\\[.4em]
&&\!\!-\ii_{\oR^\al\ra v}\circ \LLL_{\oR_\al\ra u}(\theta)\otimes_\AA s
-\ii_{\oR^\al\ra v}(\oR^\be\ra^\cop \theta)\otimes_\AA
\std{}_{\oR_\be\oR_\al\ra u} s
\\[.4em]
&=&\!\!
\ii_{[u,v]}(\theta \otimes_\AA s)~;
\end{eqnarray*}
where in the last passage the second and fourth term cancel because 
$\oR^\al\ra(\ii_v\theta) \otimes\oR_\al
=\oR^\al\ra \le v,\theta\re \otimes \oR_\al=
\le \oR^\al_{{(1)}} \ra v,\oR^\al_{{(2)}}\ra^\cop\theta\re
\otimes\oR_\al 
=\le \oR^\al\ra
          v , \oR^\be\ra^\cop\theta\re\otimes\oR_\be\oR_\al=
          \ii_{\oR^\al\ra v}(\oR^\be\ra^\cop\theta)\otimes
\oR_\be\oR_\al~.
$

The equivalent Cartan relation \eqref{CRID} is obtained by observing that
both the left hand side and the right hand side in \eqref{cartanCalculuscovder} are $\Bbbk$-linear expressions
in $u$ and in $v$,  by considering the
substitution $u\otimes v \to \oR^\al\ra v \otimes \oR_\al\ra u$,  and by recalling the braided antisymmetry of the
braided commutator (cf. the first identity in \eqref{brLie}). 
\end{proof}
\end{theorem}

The braided Cartan relation \eqref{cartanCalculuscovder} equivalently
reads $[\stDu,\ii_v]=\ii_{[u,v]}$ where the braided bracket, despite the
connection $\stD$ is not $H$-equivariant, braids nontrivially only
the vector fields $u$ and $v$, as in the Cartan relation $[\LLL_u,\ii_v]=\ii_{[u,v]}$.

\subsection{Curvature}

The curvatures of the connections $\dst\in \Con_A(\Gamma)$ and $\std\in
{}_A\Con(\Gamma)$ are respectively defined by
\begin{equation}~~~~
\Dst^{\!\!\!\!\!2}=\Dst\circ \Dst~\,,\,~~
{\stD}^{\!\!\!\!\!2}=\stD\circ \stD~.~~~
\end{equation}
These are respectively right and left $\OM$-linear maps,  
$$
\Dst^{\!\!\!\!\!2}\;\in\hom_{\OM}(\Gamma\otimes_\AA \OM, \Gamma\otimes_\AA \OMP2\,)~,~~ 
{\stD}^{\!\!\!\!\!2}\;
\in
{}_{\OM}\hom(\OM\otimes_A\Gamma, \OMP2\otimes_A\Gamma)~.$$
Recalling that $\Omega(A)={}_A\hom(\VV,A)$ we have a morphism 
\begin{equation}\label{oVVG}
{}_A\hom(\Gamma,\Omega^2(A)\otimes_A\Gamma)\longrightarrow {}_A\hom(\Vect^2(A)\otimes_A\Gamma,\Gamma)
\end{equation}
(that becomes an isomorphism 
if the $A$-module $\Vect(A)$ is  finitely generated and projective in
$\cat$, cf. Theorem \ref{Thmfgp123} and Remark \ref{newevcoev}).
In $ {}_A\hom(\Vect^2(A)\otimes_A\Gamma,\Gamma)$ we have a second definition of curvature of a
left connection $\std\in  {}_A\Con(\Gamma)$. As in \cite{NCG2}
we define the curvature $R\stdd$ to be the  $\Bbbk$-linear map
$R\stdd:\Vect(A)\otimes\Vect(A)\otimes\Gamma\to \Gamma$, $u\otimes
v\otimes s\mapsto R\stdd(u,v,s)$,
\begin{equation}\label{curvatureuvs}
R\stdd(u,v,s):=(\std{}_u\circ
\std{}_v-\std{}_{\oR^\al\ra v}\circ \std{}_{\oR_\al\ra u}-\std{}_{[u,v]})(s)~.
\end{equation}
We now show that it is a left $A$-linear map
$\Vect^2(A)\otimes_A\Gamma\to \Gamma$ and relate it to the
curvature ${\stD}^{\!\!\!\!\!2}\,\,$.
To this end we extend the evaluation $\le~,~\re:
\TT^{0,r}\otimes_A\TT^{p,q}\to \TT^{p-r,q}$, of $\TT^{0,r}$ on $\TT^{p,q}$,  defined in \eqref{evTT},
to the evaluation of 
$\TT^{0,r}$ on $\TT^{p,q}\otimes_A\Gamma$,
\begin{equation}
\le~,~~\re 
: \TT^{0,r} \otimes_\AA\,\TT^{p,q}\otimes_A\Gamma\;\longrightarrow\;
\TT^{p-r,q}\otimes_A\Gamma~,~~\nu\otimes_A\,\tau\otimes_A s\mapsto \le\nu,\tau\otimes_A s\re:=\le\nu,\tau\re s~,
\label{oppopairing}
\end{equation}
that is a morphism in $\cat$. Similarly, from the duality between
$\Vect^r(A)$ and $\Omega^r(A)$ (submodules of  $\TT^{0,r}$ and
$\TT^{r,0}$) we have the evaluation
of $\Vect^r(A)$ on $\Omega^r(A)\otimes_A\Gamma$, that we still denote
$\le~,~~\re 
: \Vect^r(A)  \otimes_\AA\,\Omega^r(A)\otimes_A\Gamma\;\rightarrow\;\Gamma$.

\begin{theorem}\label{RD2}
Let $\VV$ be finitely generated and projective as $A$-module in $\cat$
and  let $\Gamma$ in $\cat$.  Consider a left connection $\std\in
{}_A\Con(\Gamma)$. 

i) 
The curvature $R\stdd$
defined in 
\eqref{curvatureuvs} satisfies, for all $u,v\in \Vect(A)$ and $s\in \Gamma$,
\begin{equation}
R\stdd(u,v,s)=
-\,\ii_ u\circ \ii_v\circ {\stD}^{\!\!\!\!\!2}\;(s)~.
\end{equation}

ii) The curvature $R\stdd\in {}_A\hom(\Vect^2(A)\otimes_A\Gamma,\Gamma)$. 
\begin{proof} Part i):
\begin{equation*}\begin{aligned}
~~~~~~~~~~~~~\,\ii_ u\circ \ii_v\circ {\stD}^{\!\!\!\!\!2}\;(s)&=\,\ii_ u\big(
                                                    \ii_{v}\circ
                                                    {\stD}(\std
                                                    s)\big)=\ii_u \big(\stDv(\std
                                                    s)-\stD(\std{}_v
                                                    s)\big)\\[.2em]
&=\ii_{[u,v]}(\std s)+\stDalpv (\ii_{\oR_\al\ra u}\std
  s)-\std{}_u(\std{}_v s)\\[.2em]
&=-(\std{}_u\circ
\std{}_vs-\std{}_{\oR^\al\ra v}\circ \std{}_{\oR_\al\ra
  u}s-\std{}_{[u,v]}s)=-R\stdd(u,v,s)
\end{aligned}
\end{equation*}
where in the second equality we added and subtracted $\stD\circ
\ii_v$, in the third we used the Cartan relation \eqref{CRID} of Theorem \ref{THM44}.
\\

\noindent Part ii):
Another expression for the curvature is 
\begin{equation}\begin{split}
R\stdd(u,v,s)&=-\,\ii_ u\circ \ii_v\circ {\stD}^{\!\!\!\!\!2}\;(s)=
-\,\ii_ u\le v\,,{\stD}^{\!\!\!\!\!2}\;s\re=-\le u\,,\le v, {\stD}^{\!\!\!\!\!2}\;s\re\re=
-\le
  u\otimes_A v, {\stD}^{\!\!\!\!\!2}\;s\re\\
&=-\frac{1}{2}\le u\wedge
  v, {\stD}^{\!\!\!\!\!2}\; s\re~,
\end{split}
\end{equation}
where in the last equality we used that $\Vect(A)\otimes_A\Vect(A)$ is
the direct sum of braided antisymmetric plus braided symmetric vector fields,
and that these latter have vanishing pairing with 2-forms. Since ${\stD}^{\!\!\!\!\!2}_{}\in 
{}_{A}\hom(\Gamma, \Omega^2(A)\otimes_A\Gamma)$
and $\le~\,,\,~\re :
\Vect^2(A)\otimes_A\Omega^2(A)\otimes_A\Gamma\longrightarrow
\Gamma$ we see that $R\stdd$ is well-defined
as  a map $\Vect^2(A)\otimes_A\Gamma\to \Gamma$. In other terms we can write
$$
R\stdd(u,v,s)=\frac{1}{2}R\stdd(u\wedge v\otimes_As)~.
$$
Moreover, $\le~,~\re$ is left $A$-linear because it is a morphism in
$\cat$, hence also $R\stdd$ is left $A$-linear: for every $a\in A,\; u,v\in
\Vect(A)$ and $s\in \Gamma$, $R\stdd(au,v,s)=a R\stdd(u,v,s)$. Finally, we
have $R\stdd\in {}_A\hom({\Vect^2(A)\otimes_A\Gamma},\Gamma)$ because $R\stdd$ transforms according to the $H$-adjoint action
$\RA^\cop$, indeed, for all $h\in H$, 
\begin{equation*}
\begin{split}
h\ra \big(R\stdd(u,v,s)\big)&=
-\frac{1}{2}\le h_{(1)}\ra u\,\wedge h_{(2)}\ra v , (h_{(4)}\RA^\cop
{\stD}^{\!\!\!\!\!2}\;)(h_{(3)}\ra s) \re\\[.2em]
&=
(h_{(4)}\RA^\cop R\stdd)(h_{(1)}\ra u,h_{(2)}\ra 
v,h_{(3)}\ra s)~,
\end{split}
\end{equation*}
i.e., $h\ra \big(R\stdd(u\wedge v\otimes_A s)\big)=(h_{(2)}\RA^\cop
R\stdd)(h_{(1)}\ra (u\wedge v\otimes_A s))$.
\end{proof}
\end{theorem}

An analogous theorem holds for right connections; it involves
considering the $A$-bimodule of forms $\hom_A(\VV,A)$. This is obtained via the isomorphism
${\cal D}^{-1}_{\VV, A}: \Omega (A)\!= \!{}_A\hom(\VV,A)\rightarrow
\hom_A(\VV,A)$, cf. \eqref{Dcallr}.
Equivalently, one can use the evaluation map 
\begin{equation}\label{VVtildeOm}
\langle~,~\rangle':=
\langle~,~\rangle\circ \tau^{}_{\Omega(A),\Vect(A)}:
\Omega(A)\otimes_A\Vect(A)\to A~.
\end{equation}
\subsection{Torsion}
In this subsection we set $\Gamma=\Vect(A)$ with $\VV$ finitely
generated and projective as $A$-module  in $\cat$. As usual
$\Omega(A)={}^*\Vect(A)={}_A\hom(\VV,A)$ is the right dual module. Consider the canonical element
$I:=\coev(1_A)\in \Omega(A)\otimes\Vect(A)$.
Notice that $I$ is invariant under the $H$-action: for all $h\in H,
h\ra I=\varepsilon (h)I$ (indeed, so is $1_A$, and $\coev: A\to \Omega(A)\otimes_A\Vect(A)$ is a
morphism in $\cat$).

Given a left  
connection $\std$, we define the associated torsion 2-form with values in vector
fields to be the tensor field
$$
\stD(I)\in \Omega^2(A)\otimes_A \Vect(A)~.
$$
We also define the torsion $T\stdd$ as the 
$\Bbbk$-linear map
$T\stdd:\Vect(A)\otimes\Vect(A)\to\Vect(A)$, $u\otimes
v\mapsto T\stdd(u,v)$,
\begin{equation}\label{Torsion-uv}
T\stdd(u,v):=\std{}_uv-\std{}_{\oR^\al\ra v}\:\!{\oR_\al\ra u}\,-\:\!{[u,v]}~.
\end{equation}
We show that it is a left $A$-linear map
$\Vect^2(A)\to  \Vect(A)$ and relate it to the
torsion 2-form $\stD(I)$.
\begin{theorem}\label{Tasd}
Consider a left connection $\std\in
{}_A\Con(\Vect(A))$. 
For all $u,v\in \Vect(A)$ we have
\begin{equation}
T\stdd(u,v)=
-\,\ii_ u\circ \ii_v\circ {\stD}\,(I)~,
\end{equation}
hence $T\stdd\in  {}_A\hom(\Vect^2(A),\Vect(A))$. 
\begin{proof}
The first of the coherence conditions \eqref{coerecoev}
between the evaluation and coevaluation maps
$\le~,~\re:\Vect(A)\otimes_A\Omega(A)\to A$ and $\coev:A\to
\Omega(A)\otimes\Vect(A)$ equivalently reads, for all $v\in \Vect(A)$, $\ii_v\:\!I=v$.
We then recall the definition of the covariant derivative in equation \eqref{defCD}, and compute, for all $u,v\in \Vect(A)$,
\begin{align}{}
~~~~~~~~~~~~~~~~~~~~~~~~~~\ii_ u\circ \ii_v\circ {\stD}\,(I)&=
\ii_ u\circ \stDv\,(I)- \ii_u\circ {\stD}\circ  \ii_v \,(I)\nonumber \\ 
&=\ii_{[u,v]}(I)+\stDalpv \circ \ii_{\oR_\al\ra u} (I)- \ii_u\circ 
  {\stD}(v)\nonumber\\
& =-\big({[u,v]}+\stDalpv (\oR_\al\ra u) - \std{}_uv\big)\nonumber\\
&=-T\stdd(u,v)
\end{align}
where in the second line we used the Cartan relation \eqref{CRID} of
Theorem \ref{THM44}.
The equality 
$$T\stdd(u,v)=-\ii_ u\circ \ii_v\circ{\stD}\,(I)=-\ii_ u\le v,
{\stD}\,I\re=-\le u\otimes_A v,{\stD}\,I\re=-\frac{1}{2}\le
u\wedge v,{\stD}\,I\re~$$ 
shows that $T\stdd$ is a well-defined map $\Vect^2(A)\to
\Vect(A)$;  we can hence write
$T\stdd(u,v)=\frac{1}{2}T\stdd(u\wedge v)$.
Left $A$-linearity of $T\stdd$
immediately follows from left $A$-linearity of $\le~,~\re$. In order to show that $T\stdd\in
{}_A\hom(\Vect^2(A),\Vect(A))$ it remains
to prove that
$T\stdd$ transforms according to the $H$-adjoint action
$\RA^\cop$, indeed, recalling the $H$-coinvariance of the canonical element
$I$ we have, for all $h\in H$, 
\begin{align*}h\/\ra \big(T\stdd(u,v)\big)&=\,
-\frac{1}{2}\le h_{(1)}\ra u\wedge   h_{(2)}\ra v , (h_{(4)}\RA^\cop 
{\stD})(h_{(3)}\ra I) \re\\ &=\,
-\frac{1}{2}\le h_{(1)}\ra u\wedge   h_{(2)}\ra v , (h_{(3)}\RA^\cop 
{\stD})(I) \re=
(h_{(3)}\ra T\stdd)(h_{(1)}\ra u,h_{(2)}\ra v)
\end{align*}
i.e., $h\ra \big(T\stdd(u\wedge v)\big)=(h_{(2)}\RA^\cop
T\stdd)(h_{(1)}\ra (u\wedge v))$.
\end{proof}
\end{theorem}
Notice that recalling the second definition in \eqref{hDD} we also have 
$h \RA^\cop T\stdd=
T^{\phantom{J\!}}_{h{{{\,\RA^{\!\: \cop\!\:}}}}}\stdd$. 
\\

An analogous result holds for right connections. In this case we consider
 $\Vect(A)$ right dual to $\Omega(A)$ with evaluation map as in
 \eqref{VVtildeOm} and coevaluation map $
\coev' :=
\tau^{}_{\Omega(A),\Vect(A)}\circ\coev: A\to
\Vect(A)\otimes_A\Omega(A)$, cf. equation \eqref{tevcoev}. 
 We further consider the canonical element
$I'=\coev'(1_A)\in \Vect(A)\otimes_\AA \Omega(A)$. 
The torsion of a right connection $\dst$ is then defined by
$\Dst(I')\in \Vect(A)\otimes_A\Omega^{2}(A)$.

\subsection{Sum of connections}
Let $\Gamma$ and $\widehat\Gamma$ be two modules in $\cat$ with
left (right) connections. The sum of these left (right) connections
is the left (right) connection on the tensor product module
$\Gamma\otimes_A\widehat\Gamma$.
\begin{theorem}[Sum of connections \cite{AS}] \label{defsumofD} Given connections $\dst \in \Con_A(\Gamma)$, $\widehat\dst \in
\Con_A(\widehat\Gamma)$ and  $\std \in {}_A\Con(\Gamma)$, $\widehat\std \in
{}_A\Con(\widehat\Gamma)$ on the modules $\Gamma$ and $\widehat\Gamma$
in $\cat$ their sums, respectively
defined by 
\begin{align}
\begin{split}
~~~~~~~~~~~~~~~~~
\end{split}
\begin{split}
\dst\oplus_\RR\widehat\dst :\,\Gamma\otimes_A\widehat\Gamma&\,\longrightarrow\, \Gamma\otimes_A\widehat\Gamma\otimes_\AA\Omega(A)\\[.2em]
 s\otimes_\AA\hat s&\,\longmapsto\, \tau_{23}\circ(\dst(s)\otimes_A\hat 
                                 s)+(\oR^\al\ra s) 
                                 \otimes_A(\oR_\al\RA\widehat\dst)(\hat 
                                 s)~,\nn\\[1em]
\std\,\tilde\oplus_\RR\widehat\std
  :\,\Gamma\otimes_A\widehat\Gamma&\,\longrightarrow\, \Omega(A)\otimes_\AA \Gamma\otimes_A\widehat\Gamma\nn\\
 s\otimes_\AA\hat s&\,\longmapsto\,  (\oR^\al\RA^\cop\std)(s) 
                                 \otimes_A(\oR_\al\ra \hat
                                 s)+\tau_{12}\circ(s\otimes_A\widehat\std(\hat
                                 s))~,
\end{split}
\end{align}
where 
$\tau_{23}:\Gamma\otimes_A\Omega(A)\otimes_A\widehat\Gamma\rightarrow
\Gamma\otimes_A\widehat\Gamma\otimes_A\Omega(A)$ and
$\tau_{12}:\Gamma\otimes_A\Omega(A)\otimes_A\widehat\Gamma\rightarrow
\Omega(A)\oA\Gamma\otimes_A\widehat\Gamma$ are the braiding
isomorphisms, are well-defined connections
$\dst\oplus_\RR\widehat\dst \in
{}_A\Con(\Gamma\otimes_A\widehat\Gamma)$ and $\std \,\tilde\oplus_\RR
\widehat\std\in {}_A\Con(\Gamma\otimes_A\widehat\Gamma)$. The sums
$\oplus_\RR$ and $\tilde\oplus_\RR$ are associative and $H$-equivariant.
\end{theorem}
\begin{corollary}\label{corsumconn}Let $\bigoplus_{n\in {\mathbb{N}}}\Gamma^{\otimes
    n}$  be the tensor algebra generated by
  $\Gamma$ in $\cat$.
Connections $\dst$ and $\std$ on $\Gamma$
  uniquely lift to connections on
$\bigoplus_{n\in {\mathbb{N}}}\Gamma^{\otimes
    n}$ that we still denote
  $\dst$ and $\std$ and that are given by the braided Leibniz rules
  of Theorem \ref{defsumofD}. 
In particular, connections on $\Vect(A)$ lift to
connections on $\TT^{0,\bullet}$ and similarly, connections on
$\Omega(A)$ lift to connections on $\TT^{\bullet,0}$.
\end{corollary}

Associated with the braided Leibniz rule for connections on 
 $\TT^{0,\bullet}$ and on $\TT^{\bullet,0}$
we have the braided Leibniz rule for covariant
derivatives on $\TT^{0,\bullet}$ and on $\TT^{\bullet,0}$.
The following corollary will be used in Section 6 in order to study
connections compatible with a metric structure on $A$, these are connections on the tensor product $\Omega(A)\otimes_\AA\Omega(A)$.  

\begin{corollary}\label{Leibcovderu}
For any connection $\std\in\Con(\Vect(A))$ and vector field  $u\in
\Vect(A)$ the covariant derivative
$\std{}_u:\Vect(A)\to \Vect(A)$ lifts to the covariant derivative
$\std{}_u:\TT^{0,\bullet}\to \TT^{0,\bullet}$ defined via the braided Leibniz
rule
\begin{equation*}
\begin{split}
\std{}_u(v\otimes_A z)&=({\oR}^{\al\!}\ra^\cop\std)_u v \oA \oR{}_{\al\!}\ra\!
z\;+\; \oR^{\al\!} \ra\!\!\: v\otimes_A\std{}_{\oR_\al \ra u} z\\[.4em]
&={\oR}^\al\!\ra\!\!\; (\std{}_{\oR_\beta\ra u}\oR_\gamma\ra v)
\oA \oR_\al\oR^{\beta\!\!\;}\oR^{\gamma\!} \ra \! z
\;+\; \oR^{\al\!} \ra\!\!\; v\otimes_A\std{}_{\oR_\al \ra u} z
\end{split}
\end{equation*}
where $({\oR}^\al\ra^\cop\std)_u:=\ii_u\circ ({\oR}^\al\ra^\cop\std)$.
\begin{proof}
The first addend in the first equality is straightforward; the second
addend follows by considering the identity,
for all
$\theta\in\Omega(A)$, $z'\in \TT^{0,\bullet}$, 
\begin{equation*}
\begin{split}
\ii_u\circ\tau_{12}(v\otimes_A\theta
\oA z')&\,=\,\le u,{\oR}^\beta\!\!\:\ra^\cop\theta\re \oR_\beta\!\!\:\ra v \oA z'
\,=\,{\oR}^\gamma \oR_\beta\!\!\: \ra\!\!\: v_{\,} \oR_{\gamma\!\!\:}\ra\!\!\: \le u, \oR^\beta\RA^\cop\theta\re\oA z'\\[.4em]
&=\oR^{\al}\!\!\:
\ra v\oA\ii_{\oR_{\al}\ra u}(\theta\oA z')
\end{split}
\end{equation*}
that is due to quasitriangularity of the $\RR$-matrix.
The second equality follows from triangularity of the $\RR$-matrix and
the identities
\begin{equation*}
\begin{split}
R_\al\ra\le\oR_\beta \ra u\,,\std\oR_\gamma\ra v\re\oA 
R^\al\oR^\beta\oR^\gamma \ra z&
\,=\,
\le R_{\al_{(1)\!}}\oR_\beta \ra u\,,R_{\al_{(2)\!}}\ra (\std\oR_\gamma\ra v)\re\oA 
R^\al\oR^\beta\oR^\gamma z\\[.2em]
&\,=\,\le u \,,R_\delta\ra(\std\:\!\oR_\gamma\ra v)\re\oA R^\delta\oR^\gamma
\ra z
\end{split}
\end{equation*}
and
$
{R_{\delta}}\ra(\std \:\!\oR_\gamma\ra v) \otimes_A R^\delta \oR ^\gamma\ra z\,=\,
{R_{\delta}}\ra(\std \;\!S^{-1\!}(R_\gamma)\ra v)
\otimes_A R^\delta R ^\gamma\ra z
\,=\,
(R_\al\ra^\cop\std\;\!)v\oA R^\al\ra
z\:\!,
$
both due to the quasitriangular structure.
\end{proof}
\end{corollary}

\section{Duality and Cartan structure equations for curvature and torsion}
We study the standard duality notions of dual left or right $A$-linear
(more generally graded $\OM$-linear) maps and dual connections in the
context of modules in $\catfp$.
This leads in Section 5.2 to the relation between the
curvature on a module $\Gamma$ and the associated 
curvature on the dual module $\stG$ and similarly for the torsion. These are
the  Cartan structure equations for curvature and torsion in
braided noncommutative geometry in a global coordinate independent formalism.
Globally defined curvature and torsion coefficients, with respect to a
pair of dual bases for
the finitely generated and projective module of vector fields, are introduced and the
Bianchi identities proven.

\subsection{Connections on dual modules}\label{dualinternalmorph}
Let $\Gamma$ be in $\catfp$, the dual module $\stG={}_A\hom(\Gamma, A)$ is  in $\catfp$.
The associated modules 
$\OM\otimes_\AA\Gamma$ and $\stG\otimes_\AA \OM$, defined in
\eqref{defOMM} and  \eqref{sotwsot},  are modules in
$\Catfp$, indeed $\stG\otimes_\AA \OM$ is right dual to 
$\OM\otimes_\AA\Gamma$ in $\Cat$ with evaluation and coevaluation maps 
that are the $\OM$-linear extensions of the ones of $\Gamma$;
explicitly, using a
pair of dual bases of $\Gamma$ as in \eqref{expcoev}, 
\begin{equation}\label{pairOMG}\!\!\!
\begin{split}
\le~,~\re&: (\OM\otimes_\AA 
\Gamma)\,\otimes_\OM (\stG\otimes_\AA \OM)\rightarrow
\OM\;,~~\le\theta\otimes_A s,\sts\otimes_A\eta\re=\theta\wedge
\le s,\sts\re\wedge\eta \\[.4em]
\coev&: \OM \longrightarrow (\stG\otimes_A \OM)\otimes_{\OM}(\OM\otimes_A
\Gamma)\;,~~\coev(\theta)=\theta\wedge \sts^i\otimes_\OM s_i~,
\end{split}
\end{equation}
where, recall \eqref{sotwsot}, $\theta\wedge \sts^i\otimes_\OM s_i=
\theta\wedge (\sts^i\otimes_\AA 1_\OM)\otimes_\OM (1_\OM\otimes_\AA s_i)=
(\oR^\al \ra \sts^i)
\otimes_A(\oR_\al\ra^\cop\theta)\otimes_\OM s_i$.

We study the duals of  $\Bbbk$-linear maps, of  morphisms in ${}^H\MMM$
($\Bbbk$-linear and $H$-equivariant maps) and of  internal morphisms
in $\Cat$ (graded left $\OM$-linear or right $\OM$-linear maps). 
\begin{definition}\label{def-dualL} The dual  (or adjoint) ${}^{*\!}\tilde
  L$ of a $\Bbbk$-linear
map 
$\tilde L :  
\OM\otimes_\AA\Gamma\to
\Omega^{\bullet+|\tilde L|}\!(A)\otimes_A \Gamma$, of degree $|\tilde
L|$, 
 is the right $\OM$-linear map ${}^{*\!}\tilde L: \stG\otimes_\AA \OM\to\stG\otimes\Omega^{\bullet+|\tilde L|}\!(A)$, of degree $|{}^{*\!}\tilde L|=|\tilde L|$, defined by 
\begin{equation}\label{LL'}
\le\sigma,{}^{*\!}\tilde L (\stsigma)\re=(-1)^{|\tilde  L||\sigma|}
\le\tilde L(\sigma), \stsigma\re
~,
\end{equation}
for all $\sigma \in
\OM\otimes_\AA \Gamma$  of homogeneous form
degree $|\sigma|$ and  $\stsigma\in \stG\otimes_\AA \OM$. 

Vice versa, the dual $L^*$ of a 
 $\Bbbk$-linear map $L:\stG\otimes_\AA \OM\to \stG\otimes_\AA
\Omega^{\bullet+|L|}\!(A)$ of degree $|L|$ 
is the graded left  $\OM$-linear map $L^*:  \OM\otimes_A\Gamma\to
\Omega^{\bullet+|\stL|}\!(A)\otimes_A\Gamma$ of degree $|L^*|=|L|$ defined by
\begin{equation}\label{L'L}
\le L^*(\sigma), \stsigma\re=(-1)^{|L||\sigma|}\le\sigma,L(\stsigma)\re~,
\end{equation}
for all $\sigma \in
\OM\otimes_\AA \Gamma$  of homogeneous form
degree $|\sigma|$ and  $\stsigma\in \stG\otimes_\AA \OM$.
\end{definition}
Notice that while ${}^*\tilde L$ is a right $\OM$-linear map, $L^*$ is
a graded left $\OM$-linear map, c.f. \eqref{gradedOMlinear}.

 By $\Bbbk$-linearity the map  ${}^{*\!}{\scriptstyle(\,~)} : \tilde L\mapsto
{}^{*\!}\tilde L$, defined in \eqref{LL'}  when $\tilde L$ has homogeneous
degree, extends to arbitrary $\Bbbk$-module maps $\tilde L :  
\OM\otimes_\AA\Gamma\to
\OM\otimes_A \Gamma$. Similarly we have the $\Bbbk$-linear map
${\scriptstyle(\,~)}{}^{*} :  L\mapsto L^*$ defined on arbitrary $\Bbbk$-module maps $L :  
\stG\otimes_\AA\OM\to\stG\otimes_A \OM$.
\begin{proposition}\label{dualfunc}
The $\Bbbk$-linear and grade preserving maps
${\,}^{*\!}{\scriptstyle(\,~)}$ and ${\scriptstyle(\,~)}{}^{\!*}$
restrict to $H$-equivariant maps
\begin{equation}\label{LtL}
{}^{*\!}{\scriptstyle(\,~)}: {}_\Bbbk\hom(\OM\otimes_\AA\Gamma , \OM\otimes_A \Gamma) 
\,\longrightarrow \,{}\hom_\OM(\stG\otimes_\AA 
 \OM , \stG\otimes\OM)~,~ ~\tilde L\mapsto 
{}^*\tilde L~,~
\end{equation}\\[-3em]
\begin{equation}\label{LtLp}
{}{\scriptstyle(\,~)}^{\!*}{}: 
{}\hom_\Bbbk(\stG\otimes_\AA 
 \OM , \stG\otimes\OM)
\,\longrightarrow \,
{}_\OM\hom(\OM\otimes_\AA\Gamma , \OM\otimes_A \Gamma) ~,~ ~ L\mapsto L^*~.~
\end{equation}
The map in \eqref{LtL} further restricts to an isomorphism in  $\Catfp$,
\begin{equation}\label{LtL'}
{}^{*\!}{\scriptstyle(\,~)}: {}_{\OM}\hom(\OM\otimes_\AA\Gamma , \OM\otimes_A \Gamma) 
\,\longrightarrow \,{}\hom_\OM(\stG\otimes_\AA 
 \OM , \stG\otimes\OM)~,
\end{equation}
with inverse 
$
{\scriptstyle(\,~)}^{\!*}: {}\hom_\OM(\stG\otimes_\AA 
 \OM , \stG\otimes\OM)
\,\longrightarrow 
{}_{\OM}\hom(\OM\otimes_\AA\Gamma , \OM\otimes_A \Gamma)
$,
which is the restriction of the map in \eqref{LtLp}.

\begin{proof}
We show that the map ${}^{*\!}{\scriptstyle(\,~)}$ in
\eqref{LtL} is  
$H$-equivariant: for all $h\in H$, $h\ra{}^*\tilde L={}^*(h\ra^\cop \tilde
L)$. By $\Bbbk$-linearity it is enough to prove $H$-equivariance on
elements $\tilde L$ of homogenous degree $|\tilde L|$. This is indeed
the case because for all $\sigma \in \OM\otimes_\AA \Gamma$  of
homogeneous form degree $|\sigma|$ and 
all $\stsigma\in \stG\otimes_\AA \OM$, we have
\begin{equation}\label{equi*}
\begin{split}
\le \sigma, (h\ra {}^*\tilde L)(\stsigma)\re&=\le\sigma,h_{(1)}\ra
({}^*\tilde L (S(h_{(2)})\ra\stsigma)\re=h_{(2)}\ra\le
S^{-1}(h_{(1)})\ra\sigma, {}^*\tilde L (S(h_{(3)})\ra\stsigma)\re\\[.4em]
&=(-1)^{|\tilde L||\sigma|}h_{(2)}\ra\le \tilde
L(S^{-1}(h_{(1)})\ra\sigma),S(h_{(3)})\ra\stsigma\re\\[.4em]
&=
(-1)^{|\tilde L||\sigma|}\le h_{(2)}\ra \tilde
L(S^{-1}(h_{(1)})\ra\sigma),\stsigma\re\\[.4em]
&=(-1)^{|\tilde 
  L||\sigma|}\le (h\ra^\cop \tilde L)(\sigma),\stsigma\re\\[.4em]
&=\le\sigma , {}^*(h\ra^\cop \tilde
L)(\stsigma)\re~.
\end{split}
\end{equation}
$H$-equivariance of ${\scriptstyle(\,~)}^{\!*}$ is proven substituting
in \eqref{equi*}
${}^*\tilde L$ and $\tilde L$ with $L$ and $L^*$, respectively.
\\

The restricted map
${}^{*\!}{\scriptstyle(\,~)}$ in \eqref{LtL'} is a morphism  in
$\Cat$ since for all internal morphisms 
$\tilde L\in 
{}_{\OM}\hom(\OM\otimes_\AA\Gamma , \Omega^{\bullet+|\tilde L|}\!(A)\otimes_A \Gamma)$
and forms  $\theta\in
\OM$ of homogeneous degree $|\theta|$, 
we have ${}^*(\theta\tilde L)=\theta\,{}^{*\!}\tilde L$
and 
${}^*(\tilde L\theta)={}^{*\!}\tilde L\:\!\theta$.
We prove for example the first relation, for all 
$\sigma \in \OM\otimes_\AA \Gamma$  of
homogeneous form degree $|\sigma|$ and  $\stsigma\in \stG\otimes_\AA
\OM$, we have 
\begin{equation*}
\begin{split}
\le\sigma,{}^{*\!\:\!}(\theta\tilde L)(\stsigma)\re
&=(-1)^{|\theta\tilde  L||\sigma|}\le(\theta\tilde L)(\sigma),\stsigma\re
=(-1)^{|\theta||\tilde L|} (-1)^{|\theta||\sigma|}(-1)^{|\theta\tilde L||\sigma|}
     \le \tilde L(\sigma \wedge\theta),\stsigma\re\\[.2em]
&= (-1)^{|\tilde L||\theta\wedge\sigma|}
     \le\tilde L(\sigma \wedge\theta),\stsigma\re
=\le\sigma ,\theta \,{}^{*\!}\tilde L(\stsigma)\re
=\le\sigma ,(\theta\, {}^{*\!}\tilde L)(\stsigma)\re
\end{split}
\end{equation*}
where we used the definition \eqref{LL'},   the bimodule structure of internal
morphisms in $\Cat$ given in  \eqref{bimodofOM}, then the definition \eqref{L'L} and
again   \eqref{bimodofOM}.

Finally, the morphism ${}^{*\!}{\scriptstyle(\,~)}$ in \eqref{LtL'} is
an isomorphism, with inverse ${\scriptstyle(\,~)}^{\!*}$, since  for all $L$ and $\tilde L$,
$({}^*\tilde L )^*=\tilde L$ 
and ${}^*_{} (L^*)=L$. This immediately follows from  the definitions
\eqref{LL'} and \eqref{L'L}.
\end{proof}
\end{proposition}
We have explicit expressions for the isomorphisms
${}^{*\!}{\scriptstyle(\,~)}$ and ${\scriptstyle(\,~)}^{\!*}$ in
$\Cat$   in terms of a
dual basis of $\Gamma$, i.e., of the coevaluation map
$\coev(1_{\OM})=\sts^i\otimes_As_i$.
For all internal morphisms 
$\tilde L\in 
{}_{\OM}\hom(\OM\otimes_\AA\Gamma , \Omega^{\bullet+|\tilde L|}\!(A)\otimes_A \Gamma)$,
$L\in \hom_\OM(\stG\otimes_\AA \OM, \stG\otimes_\AA 
\Omega^{\bullet+| L|}\!(A))$ and 
 $\sigma \in \OM\otimes_\AA \Gamma$  of
homogeneous form degree $|\sigma|$ and  $\stsigma\in \stG\otimes_\AA
\OM$, we have
\begin{equation}\label{*Lexplicit}{}^{*\!}\tilde L(\stsigma)=\sts^i\le \tilde 
L(s_i),\stsigma\re~,~~L^*(\sigma)=(-1)^{|L||\sigma|}\le \sigma,
L(\sts^i)\re s_i~.
\end{equation}
This is due to the identities
$
\le\sigma, \sts^i\le\tilde L(s_i),\stsigma\re\re=
\le\sigma, \sts^i\re\le\tilde
L(s_i),\stsigma\re=(-1)^{|\tilde L||\sigma|}\le\tilde L(\le\sigma,\sts^i\re
s_i),\stsigma\re=(-1)^{|\tilde L||\sigma|}\le\tilde
L(\sigma),\stsigma\re
=\le \sigma, {}^*\tilde L(\stsigma)\re
$,
where we first used
right $\OM$-linearity of the pairing, then its left $\OM$-linearity
and that of $\tilde L$. The expression for $L^*$ is similarly proven.

\begin{definition}\label{Defdualconn}
Let $\std\in {}_A\Con(\Gamma)$. The dual of the left connection $\stD$ is the $\Bbbk$-linear map 
${}^{*\!}(\stD): \stG \otimes_\AA \OM\longrightarrow
\stG\otimes_A\Omega^{\bullet+1\!}(A) $, defined by
\begin{equation}\label{dualstconn} 
\dd\le\sigma,\stsigma\re=\le \stD\sigma,\stsigma\re+(-1)^{|\sigma|}\le\sigma,{}^{*\!}(\stD)\stsigma\re~,
\end{equation}
for all $\sigma \in  \OM\otimes_A\Gamma$ of homogeneous form degree $|\sigma|$ and $\stsigma\in \stG\otimes_\AA \OM$.
\label{dualconnst}
\end{definition}
From the definition it follows  that ${}^{*\!}(\stD)$ is a right connection.
Indeed, for all $\sigma$, $\stsigma$ of homogeneous degree
$|\stsigma|$ and
$\vartheta\in\OM$, we have the identity
$\le \sigma , {}^{*\!}(\stD)(\stsigma\wedge \vartheta)\re=
\le 
\sigma, \,{}^{*\!}(\stD)\stsigma\wedge \vartheta\re+
(-1)^{|\stsigma|}\le 
\sigma, \stsigma\wedge\dd\vartheta\re$.

Vice versa, given a right connection $\dst \in {}_A\Con(\stG)$,
equation \eqref{dualstconn}, rewritten as
\begin{equation}\label{eqdualconnst}
\dd\le\sigma,\stsigma\re=\le (\Dst)^{*\,}\sigma,\stsigma\re+(-1)^{|\sigma|}\le\sigma,\Dst\,\stsigma\re~,
\end{equation}
{\sl defines} a left connection 
$(\Dst)^*$.
Obviously, the dual of a dual connection
is the initial connection.

If $\sigma=s\in \Gamma\subset \Gamma\otimes_A\OM$ and similarly, if
$\sigma'=s'\in\Gamma'\subset\OM\otimes_A \Gamma'$,
equations \eqref{dualstconn} and \eqref{eqdualconnst} read 
\begin{equation}\label{dddual}
\dd\le s, \sts\re=\le \std s, \sts\re+\le s, {}^*\std  \;\sts\re~~,~~~\dd\le
s,\sts\re=\le \dst^* s,\sts\re +\le s,\dst\; \sts\re~
\end{equation}
and define ${}^*\std\in\Con_A(\stG)$ in terms of $\std\in
{}_A\Con(\Gamma)$, and $\dst^*\in {}_A\Con(\stG)$ in terms of
$\dst\in\Con_{A}(\stG)$. 
Since the extensions $\stD$ and $\Dst$ of the
connections $\std$ and $\dst$ are uniquely determined by the Leibniz
rule, we have  ${}^*(\stD)=\Dstst$ and  $(\Dst)^*=\Dst_{\!\!\!\!\!\!{\phantom{I^I}}^*}$.

Using a dual basis we have the explicit expressions, for all $s\in
\Gamma$, $\sts\in \stG$,
$$
{}^*\std\;\sts=\sts^i\otimes_A\dd\le
s^i,\sts\re-\sts^i\otimes_A\le\std s^i,\sts\re~,~~~\dst^*s=\dd\le
s,\sts^i\re\otimes_A s_i-\le s,\dst\,\,\sts^i\re\otimes_A s_i~.
$$
For example, pairing the first expression with $s\in \Gamma$ and
using
$\le s,\sts^i\re\le\std s_i,\sts\re=\le\std(\le s,\sts^i\re s_i),\sts\re
$ $-$ $(\dd\le s,\sts^i\re)\le s_i,\sts\re =\le \std s,\sts\re- (\dd\le s,\sts^i\re)\le s_i,\sts\re $
we  obtain the first expression in \eqref{dddual}.  
\\

The difference of  two connections $\std, \std' \in {}_A\Con(\Gamma)$ is a left
$A$-linear map, the difference of their duals is the right $A$-linear
map 
${}^{*}\:\!\std - {}^{*}\:\!\std'=-{}^{\,*\!\!\:}(\std - \std')$
where on the right hand side we used the restriction to $\catfp$ of the isomorphism
${}^{*\!}{\scriptstyle(\,~)}$ of Proposition \ref{dualfunc}.
Since also $-^{\,*\!}{\scriptstyle(\,~)}$ is an isomorphism in $\catfp$  we immediately have
\begin{corollary}\label{dualconnaffineiso}
${}^{*\!}{\scriptstyle(\,~)}:{}_A\Con(\Gamma)\rightarrow
\Con_A(\stG)$ with  $\,-{}^{\,*\!}{\scriptstyle(\,~)}: 
{}_A\hom(\Gamma,\Omega(A)\otimes_A \Gamma )\rightarrow
\hom_A(\stG,\Gamma\otimes_A\Omega(A))$
 is an isomorphism of affine spaces over modules in $\catfp$.
\end{corollary}

There is a unique way of inducing connections on
duals of tensor product modules, indeed,  the sum of dual connections is the
dual of the sum of connections.
\begin{proposition}\label{staroplus} Consider the connections $\std \in {}_A\Con(\Gamma)$, $\widehat\std \in
{}_A\Con(\widehat\Gamma)$ and  
 $\dst \in \Con_A(\Gamma)$, $\widehat\dst \in
\Con_A(\widehat\Gamma)$ on the modules $\Gamma$, $\widehat\Gamma$
in $\catfp$. Let ${}^*\std \in \Con_A({}^*\Gamma)$, ${}^*\widehat\std \in
\Con_A({}^*\widehat\Gamma)$ and $\dst^* \in {}_A\Con(\Gamma^*)$, $\widehat\dst^* \in
{}_A\Con(\widehat\Gamma^*)$
 be the connections on the dual modules.
We have $${}^*(\std \:\tilde\oplus_\RR \:\!\widehat\std\:\!)={}^*\widehat\std\:\!\oplus_\RR\!
{}^*\std~,~~
(\widehat\dst\:\!\oplus_\RR\!\!\:\dst\;\!)^*=\dst^* \,\tilde\oplus_\RR \widehat\dst^*~,~~
$$ as connections in 
${}\Con_A({}^*\widehat\Gamma\otimes_A{}^*\Gamma)$ and in
${}_A\Con(\Gamma^*\otimes_A \widehat\Gamma^*)$, respectively.
\end{proposition}

We leave the proof of this proposition to the reader. It follows from
triangularity of the $\RR$-matrix, including the compatibility
\eqref{taudual}  of the braiding with the
dual braiding. The second equality follows from the first recalling
that the dual of a dual connection is the initial connection.

\subsection{Cartan structure equations and Bianchi identities}\label{CRcurv}
According to Corollary \ref{dualconnaffineiso} we have
\begin{lemma}\label{D'2D2'}Let $\Gamma$ be a module in $\catfp$. The 
dual of the curvature 2-form of a left connection $\std\in
{}_A\Con(\Gamma)$ is minus the curvature two form of the dual
connection ${}^*\std\in  \Con_A(\stG)$, i.e., 
${{}^*({\stD}^{\!\!\!\!\!2\,})}=-\,{\Dstst}^{\!\!\!\!\!\!\!2\;\:\;}$
Similarly, for a right connection $\dst\in \Con_A(\Gamma)$,
${({\Dst}^{\!\!\!\!\!2}\,)^*} =-\,\Dst_{\!\!\!\!\!\!{\phantom{I^I}}^*}^{\!\!\!\!\:\!\!2\:}$.
\begin{proof}
Use twice \eqref{dualstconn} and  ${}^*(\stD)=\Dstst$ to obtain,
for all $\sigma \in \OM\otimes_A  \Gamma$ of homogeneous form degree
$|\sigma|$ and $\stsigma\in \stG\otimes_\AA \OM$,
\begin{equation}\label{RDRD'}
0=\dd^2\le\sigma,\stsigma\re=\dd\big(\le
\stD\sigma,\stsigma\re+(-1)^{|\sigma|}\le\sigma,
\Dstst\stsigma\re\big)=\le
{\stD}^{\!\!\!\!\!2}\,\,\sigma\,,\stsigma\re+\le\sigma,
{\Dstst}^{\!\!\!\!\!\!2\:\:\,}\stsigma\re~.
\end{equation}
By definition, the dual of the curvature 2-form satisfies, cf. \eqref{LL'},  
$\le \sigma,
{}^*({\stD}^{\!\!\!\!\!2}\,\,)\stsigma\re=\le
   {\stD}^{\!\!\!\!\!2}\,\,\sigma\,,\stsigma\re$, hence 
${{}^*({\stD}^{\!\!\!\!\!2\,})}=-\,{\Dstst}^{\!\!\!\!\!\!\!2\;\:\;}$. The
second equality, 
${({\Dst}^{\!\!\!\!\!2}\,)^*}
=-\,\Dst_{\!\!\!\!\!\!{\phantom{I^I}}^*}^{\!\!\!\!\:\!\!2\:}$,
then follows  setting $\std=\dst^*$ and recalling that the double dual of
a connection (curvature) is the original connection (curvature).
\end{proof}
\end{lemma}

Let  $\Gamma$ and $\Vect(A)$ to be in $\catfp$. Their duals 
$\stG$ and $\Omega(A)$ are also in $\catfp$. According to Remark
\ref{newevcoev} the right dual of $\Vect(A)\otimes_\AA \Gamma$
is $\stG\otimes_\AA\Omega(A)$. Similarly (cf. also Proposition
\ref{onionprop}) the right dual of
$\Tau^{0,r}\otimes_\AA \Gamma$ is $\stG\otimes_\AA\Tau^{r,0}$, $r\in
\mathbb{N}$. In the next theorem we use the associated exact pairing 
$\le~,~\re 
: \TT^{0,r} \otimes_A\Gamma\otimes_\AA\,\stG\otimes_A\TT^{r,0}\longrightarrow A$
with $r=2$, and also the exact pairing in $\Catfp$ 
defined in \eqref{pairOMG}. 
\begin{theorem}[Second Cartan structure equation] \label{CartanII}
Let $\Gamma$ and $\VV$ be modules in $\catfp$ and let $\std\in {}_A\Con(\Gamma)$. For all
$u,v\in\Vect(A)$, $s\in \Gamma$, $\sts\in\stG$ we have
\begin{equation*}
\le R\stdd(u,v,s), \sts\re=\le u\otimes_A v\otimes_A s,
{\Dstst}^{\!\!\:\!\!\!\!\!\!2}\,\:\:\sts\re
\end{equation*}
or, equivalently, $\le R\stdd(u, v,s), \sts\re=\frac{1}{2}\le u\wedge
v\otimes_A s,
{\Dstst}^{\!\!\:\!\!\!\!\!\!2}\,\:\: \sts\re$.
\begin{proof}
First notice that for all $\vartheta\in \Omega^2\subset \TT^{2,0}(A), s\in
\Gamma, \sts\in\stG$, considering the pairing of
$\OM\otimes_A \Gamma$ with $\stG\otimes_A\OM$ defined in equation
\eqref{pairOMG} we have
\begin{equation}\label{pairpairs}
\begin{split}
\le \ii_u\circ \ii_v (\vartheta\otimes_A s),\sts\re&=\le \le u\otimes_A v,
\vartheta\re s,\sts\re= \le u\otimes_Av,
\vartheta\re\le s,\sts\re\\[.4em] &=
\le u\otimes_Av,
\vartheta \le s,\sts\re\re=
\le u\otimes_Av,
\le\vartheta\otimes_A s,\sts\re\re\\[.4em] &=
 \ii_u\circ \ii_v \le \vartheta\otimes_A s,\sts\re
\end{split}
\end{equation}
where in the second line we first used right $A$-linearity of
$\le~,~\re: \Vect^{\otimes 2}\otimes_A\Omega^{\otimes 2}\to A$ and
then left $\OM$-linearity of the pairing in \eqref{pairOMG}. 
Then, recalling also  Theorem \ref{RD2} and  Lemma \ref{D'2D2'}, we have
\begin{equation}
\begin{split}\le R\stdd(u,v,s),\sts\re&=
-\le\,\ii_ u\circ \ii_v\circ {\stD}^{\!\!\!\!\!2}\;(s),\sts\re=
-\,\ii_ u \circ \ii_v \le  {\stD}^{\!\!\!\!\!2}\;(s),\sts\re
=\,\ii_ u \circ \ii_v  \le s,
{\Dstst}^{\!\!\:\!\!\!\!\!\!2}\,\:\: \sts\re\\[.4em]
&=\le u\otimes_A v,\le s,
{\Dstst}^{\!\!\:\!\!\!\!\!\!2}\,\:\: \sts\re\re=
\le u\otimes_A v\otimes_A s,
{\Dstst}^{\!\!\:\!\!\!\!\!\!2}\,\:\: \sts\re
\end{split}\end{equation}
where in the last equality we used right $\OM$-linearity of the pairing
in \eqref{pairOMG}. The equivalent expression  $\le R\stdd(u, v,s), \sts\re=\frac{1}{2}\le u\wedge
v\otimes_A s,
{\Dstst}^{\!\!\:\!\!\!\!\!\!2}\,\:\: \sts\re$ trivially follows from
$u\wedge v=u\otimes_A v-\oR^\al\ra v\/\otimes_A \oR_\al\ra u$.
\end{proof}
\end{theorem}

\begin{theorem}[First Cartan structure equation]\label{CartanI}
Let $\VV$ be in $\catfp$ and $\std\in {}_A\Con(\VV)$. For all
$u,v\in\Vect(A)$, $\theta\in \Omega(A)$, we have
\begin{equation*}
\le T\stdd(u,v), \theta\re=-\le u\otimes_A v,
(\dd+\wedge\circ {}^{\!*}\std\;\!)\,\theta\re
\end{equation*}
or, equivalently,
$\le T\stdd(u,v), \theta\re=-\frac{1}{2}\le u \wedge v,
(\dd+\wedge\circ {}^{\!*}\std\;\!)\theta\re$. 
\begin{proof}
\begin{equation*}
\begin{split}
\le T\stdd(u,v), \theta\re
&=-\le\,\ii_ u\circ \ii_v\circ {\stD}(I),\theta\re
=-\,\ii_ u \circ \ii_v \le  {\stD}(I),\theta\re
=-\,\ii_ u \circ \ii_v (\dd \le I,\theta\re+\le I,{}^{*\!\;\!}\std\,\theta\re)\\[.4em]
&=-\,\ii_ u \circ \ii_v (\dd\theta+\wedge\circ  {}^{*\!\;\!}\std\,\theta)
=-\le u\otimes_A v , \dd\theta+\wedge\circ {}^{*\!\;\!}\std\,\theta\re
\end{split}\end{equation*}
where we used Theorem \ref{Tasd}, then \eqref{pairpairs} with
$\Gamma=\Vect(A)$, next the Definition \ref{Defdualconn} 
of dual connection with $\stG=\Omega(A)$, and in the second line 
$\OM$-bilinearity of the pairing 
$\le~,~\re:\OM\otimes_A \Vect(A)\otimes_A\Omega(A)\otimes_A\OM\to
\OM$, so that $\le I,\theta\re=\le\coev(1_{\OM}),\theta\re=
\le\omega^i\otimes_Ae_i ,\theta\re=\theta$ and $\le
I,\omega\otimes_A\eta\re=\le\omega^i\otimes_A
e_i,\omega\otimes_A\eta\re=\omega^i\wedge\le
e_i,\omega\re\wedge\eta=\omega\wedge \eta$ for all $\omega\in
\Omega(A), \eta\in \OM$.
The equivalent expression $\le T\stdd(u,v), \theta\re=-\frac{1}{2}\le u \wedge v,
(\wedge\circ {}^{\!*}\std+\dd)\theta\re$ trivially follows from
$u\wedge v=u\otimes v-\oR^\al\ra v\/\otimes\/ \oR_\al\ra u$.\end{proof}
\end{theorem}

In the proof we have shown that for all $\theta\in \Omega(A)$, $\le  {\stD}(I),\theta\re=(\dd +
\wedge\circ {}^{*\!\;\!}\std\;\!)\theta$.
This defines the {\it torsion on one forms}
\begin{equation}\label{tordualCart}
(\dd+\wedge\circ {\!\!\:}^*\std\;\!):
  \Omega(A)\to \Omega^2(A)
\end{equation} of an arbitrary right connection  ${}^*\std\in
\Con_A(\Omega(A))$ and shows that it is an internal morphism in
  $\hom_A(\Omega,\Omega^2(A))$. 
As for the curvature ${\Dstst}^{\!\!\:\!\!\!\!\!\!2}~~\!$,
there is a unique lift to the torsion
$(\dd\circ \wedge\:\!+\:\!\wedge\circ\Dstst)\in
\hom_\OM(\Omega(A)\otimes_A\OM,\OMP2{\,})$.

\begin{remark}In Theorems \ref{CartanII} and \ref{CartanI} we use the pairing 
$\le~,~\re 
: \TT^{0,2} \otimes_A\Gamma\otimes_\AA\,\stG\otimes_A\TT^{2,0}\;\to
A$. The pairing between 2-vector fields and 2-forms $\le~,~\re^{}_{\!\wedge} 
: \Vect^2(A) \otimes_A\Gamma\otimes_\AA\,\stG\otimes_A\Omega^2(A)\;\to
A$ is half the value of the restriction of the first pairing to
$\Vect^2(A)\subset \TT^{0,2}$ and $\Omega^2(A)\subset \TT^{2,0}$, so that the Cartan structure equations read $\le R\stdd(u, v,s), \sts\re=\le u\wedge
v\otimes_A     s,
{\Dstst}^{\!\!\:\!\!\!\!\!\!2}\,\:\: \sts\re^{}_{\!\wedge}$ and 
$\le T\stdd(u,v), \theta\re=-\le u \wedge v,
(\dd+\wedge\circ {}^{\!*}\std\,)\theta\re^{}_{\!\wedge}$. 
\end{remark}

Using a dual basis $\{e_i,\omega^i : i=1\ldots, n\}$ of $\VV$  
we define the curvature and torsion coefficients of a connection
$\std\in {}_A\Con(\VV)$, 
\begin{equation}
R^{}_{ijk}{}_{}^l:=\le R\stdd(e_i,e_j,
e_k),\omega^l\re~,~~T^{}_{ij}{}^l_{}:=\le T\stdd(e_i,e_j),\omega^l\re~
\end{equation}
and the curvature and torsion two forms (the signs are chosen to match the
commutative differential geometry case) 
\begin{equation}\label{Rsf-def}
\Rsf^{}_k{}_{}^l:= -\frac{1}{2}\omega^j\wedge\omega^i
R^{}_{ijk}{}_{}^l~,~~ \Tsf^l: =-\frac{1}{2}\omega^j\wedge\omega^i\:\! T^{}_{ij}{}^l~.
\end{equation}
Since $\omega^j\oA\omega^i\otimes_A\le e_i\oA e_j,
~~~\re$ is the identity map on   $\Omega(A)\oA\Omega(A)$ we have the equality
$\omega^j\otimes \omega^i\le e_i\oA e_j\oA e_k , 
{\Dstst}^{\!\!\:\!\!\!\!\!\!2}\,\:\:
\omega^l\re 
=
\omega^j\otimes \omega^i\le e_i\oA e_j,\le e_k , 
{\Dstst}^{\!\!\:\!\!\!\!\!\!2}\,\:\:
\omega^l\re\re 
=\le 
e_k,
{\Dstst}^{\!\!\:\!\!\!\!\!\!2}\,\:\:
\omega^l\re$ and similarly 
$\omega^j\otimes \omega^i\le e_i\oA e_j,
(\dd+\wedge\circ
{\!\!\:}^*\std\;\!)\omega^l\re=
(\dd+\wedge\circ
{\!\!\:}^*\std\;\!)\omega^l$. This leads to the coefficient expression
of the Cartan structure equations
$$
\le e_k,{\Dstst}^{\!\!\:\!\!\!\!\!\!2}\,\:\:
\omega^l\re=\frac{1}{2}\omega^j\wedge\omega^i R^{}_{ijk}{}_{}^l
=-\:\!\Rsf^{}_k{}_{}^l~,$$
$$
(\dd+\wedge\circ
{\!\!\:}^*\std\;\!)\omega^l=-\frac{1}{2}\omega^j\wedge\omega^i\:\! T^{}_{ij}{}^l= \Tsf^l~.
$$

We similarly define the coefficients one forms of the
connection ${}^*\std\in \Con_A(\Omega(A))$, dual to $\std\in
{}_A\Con(\VV)$,
$$\omega{}_k{}^l:=\le e_k,{}^*\std\omega^l\re~,
$$
so that, since $\omega^k\otimes_A\le e_k,
~\re$ is the identity map on   $\Omega(A)$,
  $${}^*\std\omega^l=\omega^k\otimes_A
\omega{}_k{}^l~.$$ In terms of these
coefficients we obtain
$$
{\Dstst}^{\!\!\:\!\!\!\!\!\!2}\,\:\: \omega^l=
\omega^k\oA(\dd\omega^{}_k{}{}_{}^l+\omega^{}_k{}^j_{}\wedge\omega^{}_j{}{}_{}^l)
$$
$$(\dd+\wedge\circ
{\!\!\:}^*\std\;\!)\omega^l=\dd\omega^l+\omega^j\wedge \omega^{}_j{}{}_{}^l$$
which, together with the identity $\omega^k\oA\le e_k,
{\Dstst}^{\!\!\:\!\!\!\!\!\!2}\,\:\: \omega^l\re={\Dstst}^{\!\!\:\!\!\!\!\!\!2}\,\:\: \omega^l$,
give the full coefficient expression of the Cartan structure equations
$$
\omega^k\oA(\dd\omega^{}_k{}{}_{}^l+\omega^{}_k{}^j_{}\wedge\omega^{}_j{}{}_{}^l)=\omega^k\otimes_A (-\!\:\Rsf^{}_k{}_{}^l)~,$$
$$
\dd\omega^l+\omega^j\wedge \omega^{}_j{}{}_{}^l=
\Tsf^l~.
$$

As in commutative differential geometry, applying
$\id_{\Omega(A)}\otimes_A \dd$ to the first equation and differentiating
the second  we readily obtain the  Bianchi identities,
\begin{equation}\label{Bianchiid}
\begin{split}
&\omega^k\oA (\dd\Rsf_k{}^l+\omega_k{}^j\wedge\Rsf_j{}^l-\Rsf_k{}^j\wedge\omega_j{}^l)=0~,\\[.4em]
&\dd\Tsf^l-\Tsf^j\wedge\omega_j{}^l=\omega^j\wedge\Rsf_j{}^l~.
\end{split}
\end{equation}
Notice that the commutator $[\omega,
\Rsf]_k^{~l}:=\omega_k{}^j\wedge\Rsf_j{}^l-\Rsf_k{}^j\wedge\omega_j{}^l$
in the first identity is not a braided commutator.

Similarly, if we consider  a connection $\std\in{}_A\Con(\Gamma)$ on a
module $\Gamma$ in  $\catfp$, by setting
$\omega^{}_k{}_{}^l:=\le s_k,{}^{*\!}\std\:\sts^l\re$ (with $\{s_i,\sts^i:
i=1\ldots, m\}$ a dual basis) we have the
Bianchi identity   $\sts^k\oA
(\dd\Rsf_k{}^l+\omega_k{}^j\wedge\Rsf_j{}^l-\Rsf_k{}^j\wedge\omega_j{}^l)=0$,  where $
\Rsf^{}_k{}_{}^l$  is defined in \eqref{Rsf-def} and
where now $R^{}_{ijk}{}_{}^l:=\le R\stdd(e_i,e_j,s_k),\sts^l\re$.
Equivalently,
$\le s_p,\sts^k\re \dd\Rsf_k{}^l+\omega_{p}^j\wedge\Rsf_j{}^l-\Rsf_p{}^j\wedge\omega_j{}^l=0$;
this last expression is obtained observing that if $P\in
\stG\otimes_A \OM$ then $\le s_l,\sts^k\re P_k=P_l$, where
$P_k:=\le s_k,P\re$, and using that 
$
\le s_k,{\Dstst}^{\!\!\:\!\!\!\!\!\!2}\,\:\:
\sts^l\re
=-\:\!\Rsf^{}_k{}_{}^l$.
\section{Riemannian geometry}\label{riemgeom}
We use the sum of connections (based on the tensor product of
internal morphisms in $\MMM^H$), the notion of dual connection and the
Cartan calculus results for the torsion in order to determine the
Koszul formula for the Levi-Civita connection
associated with a  pseudo-Riemannian metric tensor on the algebra
$A$. Hence its uniqueness, existence and explicit expression.
The Ricci tensor is canonically introduced leading to the notion of
noncommutative Einstein manifold. 
\subsection{Metric tensor}
Let  $A$ be a braided commutative $H$-module algebra with $H$ triangular and $A$-bimodule of braided derivations $\VV$ that is finitely generated and projective, hence in $\catfp$. Let $\Omega(A)={}^*\VV={}_A\hom(\VV,A)$  be the dual $A$-bimodule of forms. 
Proposition \ref{propo:projhomiso} and Theorem \ref{Thmfgp123} give two isomorphisms
in $\catfp$, that with slight abuse of notation are both denoted by
$\flat$, 
\begin{equation}\label{gbemoll}\begin{aligned}
\Omega(A)\otimes_A \Omega(A)&\stackrel{\flat_{}\,}{\longrightarrow}
                                     {}_A\hom(\Vect(A),\Omega(A))\stackrel{\flat}{\longrightarrow}
                                     {}_A\hom(\VV\oA\VV,A) \\
\g=\g^a\oA\g_a& \longmapsto ~\g^\flat:=\le\:\,,\g^a\re\oA \g_a
                                                                 \longmapsto ~\le\,{\mbox{-}}\oA{\mbox{-}}\,,\g\re=\le\,~,\le\,~,\g^a\re\g_a\re~
\end{aligned}
\end{equation}
(sum on the index $a$ understood).
\begin{definition}
A pseudo-Riemannian metric on $\VV$ is a braided symmetric element $\g\in\Omega(A)\otimes_A\Omega(A)$, i.e., $\g=\tau (\g)$, with associated internal morphism $\g^\flat\in {}_A\hom(\Vect(A),\Omega(A))$ that is an isomorphism. We also simply say that $\g$ is a metric on $A$.   
\end{definition}

Similarly to \eqref{gbemoll}, consider $
\flat':\VV\otimes_A \VV \rightarrow {}_A\hom(\Omega(A),\VV),
~\fgf^a\oA\fgf_a \mapsto
                           \le~\,,\fgf^a\re'\oA \fgf_a$,
(cf. Proposition \ref{propep'} and Remark \ref{ep'}) and let $\sharp':={\flat'}^{-1}$ be the
inverse isomorphism. The metric tensor $\g\in \Omega(A)\oA\Omega(A)$
induces the metric on $\OOm$ (or inverse metric tensor)  $\bar\g$, defined by
\begin{equation}\label{invmetric} \bar\g:={(\g^\flat)^{-1}}^{\!\:\sharp'}\!\in\VV\oA\VV~. 
\end{equation}Setting $\bar \g=\bar \g^b\oA\bar \g_b$ we
have ${\bar\g\!\:}^{\flat'}\!=\le \,~,\bar\g^b\re'\oA\bar\g_b$, and
the equality ${\bar\g\:\!}^{\flat'}=({\g^\flat})^{-1}$ reads
$$
\le\,~,\g^a\re\le\g_a,\bar\g^b\re'\oA\bar\g_b=\id_{\VV}~\,,~~
\le\,~,\bar\g^b\re'\le\bar\g_b,\g^a\re\oA \g_a=\id_{\OOm}~.
$$
Using dual bases for $\VV$ as a right and as a left rigid $A$-module in
$\catfp$, that is $\coev(1_A)=\omega^i\oA e_i$ and
$\coev'(1_A)={}^\al e_i\oA{}_\al\omega^i:=\oR^\al\ra e_i\oA\oR_\al \ra \omega^i
$ we have $\id_{\VV}=\le ~\,, \omega^i\re\oA e_i$ and $\id_{\OOm}=\le~\,,{}^\al e_i\re\oA{}_\al\omega^i$.
Then, recalling that  $\le~,~\re$ and  $\le~,~\re'$ induce the isomorphisms
$\OOm\otimes_A\VV\simeq {}_A\hom(\VV,\VV)$ and $ 
\VV\oA\Omega(A)\simeq {}_A\hom(\Omega(A),\Omega(A))$, we see that $\bar \g$ is the  inverse metric of $\g$ if and only
if 
\begin{equation*}\label{gbarg}
\g^a \le\g_a,\bar\g^b\re'\oA\bar\g_b=\omega^i\oA e_i~,~~
\bar\g^b\otimes_A \le\bar\g_b,\g^a\re \oA\g_a=
{}^\al e_i\oA{}_\al\omega^i~.
\end{equation*}
 
 \subsection{Levi-Civita connection}
We prove existence and uniqueness of a metric compatible and torsion
free connection establishing a noncommutative Koszul formula.

In this section we simplify the notation of the braiding via the
$\RR$-matrix action 
and set, for any  $w\in W$ with $W$ module in ${}^H\MMM$,  $\cat$ or $\Cat$,
${}^\al w=\oR^\al\ra w$ and ${}^{}_\al w=\oR_\al\ra w$ (for any $\alpha$). Hence, for
example, for all $u,v\in \Vect(A)$, $\theta\in \Omega(A)$, 
\begin{equation*}
\begin{split}
&{}^\al v\oA{}_\al u=\oR^\al\ra v\:\!\oA\oR_\al\ra u=\tau(u\otimes_A v)~,\\[.4em]
&{}^\beta\theta\otimes (\!\:{}^{}_\beta{}^*\std)=\oR^\beta\ra 
\theta\otimes(\oR_\beta\ra{}^*\std)~,~~({}^\al\std) \otimes {}_\al u
=(\oR^\al\ra^\cop \std) \otimes {}_\al u~.
\end{split}
\end{equation*}

Recall that considering sums of connections
(connections on tensor products) and dual connections a left connection $\std\in
{}_A\Con(\Vect(A))$ uniquely lifts to a left connection on
$\TT^{0,\bullet}(A)$ and to a dual right connection on $\TT^{\bullet,0}(A)$,
cf. Corollary \ref{corsumconn} and Proposition \ref{staroplus}. 
For example
on the metric tensor we have (cf. Theorem \ref{defsumofD})
$${}^{*}\std(\g)={}^{*}\std (\g^a\otimes_A \g_a)=\tau_{23}\circ ( {}^{*}\std(\g^a)\otimes_A\g_a)+\,{}^{\beta\;\!}\g^a\otimes_A({}_\beta{}^*\std)(\g_a)~.$$ 
 Similarly for a right connection $\dst\in
\Con_A(\Omega(A))$. Moreover, this latter is torsion free if its dual
 $\dst^*\in
{}_A\Con(\Vect(A))$ is torsion free, cf. Theorem \ref{CartanI} and \eqref{tordualCart}.
\begin{definition}Let $\g\in \Omega(A)\oA\Omega(A)$ be a
  pseudo-Riemannian metric on $A$. A right connection 
 $\dst\in \Con_A(\Omega(A))$ is metric compatible if it satisfies $\dst(\g)=0$.
A left connection
$\std\in {}_A\Con(\Vect(A))$ is  metric compatible if 
its dual ${}^{*}\std\in \Con_A(\Omega(A))$
is metric compatible.

A Levi-Civita connection is a  metric compatible and torsion free connection. 
\end{definition}

For ease of the reader in the statements of the next two theorems we  spell out the
condition that $\VV$ is a  module in $\catfp$.
\begin{theorem}[Uniqueness of Levi-Civita connection]\label{LCU}  Let
  $H$ be a triangular Hopf algebra,  $A$ a braided commutative
  $H$-module algebra, let the associated module in $\cat$ of braided derivations $\VV$ be finitely
  generated and projective, and let $\g\in \OOm\oA\OOm$ be a metric on $A$, where $\OOm={}^*\VV$.
If a torsion free metric compatible left connection $\std\in
{}_A\Con(\Vect(A))$ exists, it is unique.
\end{theorem}
\begin{proof} Assume $\std\in
{}_A\Con(\Vect(A))$ is a torsion free metric compatible
connection.  Applying the contraction operator to the identity
$$\dd\le v\oA z,\g\re=\le\std (v\oA z),\g\re+\le v\oA
z,{}^{*}\std \g\re= \le\std (v\oA z),\g\re$$ we have, for all
$u,v,z\in \Vect(A)$,
\begin{equation}\label{Luvzg}
\begin{split}\LLL_u\le v\oA z,\g\re&=
\le \std{}_u(v\oA z),\g\re\\[.4em]
&=\le\!\: {}^{\al\:\!\!}(\std{}_{{}^{\!}_\be u\,}{}_\gamma v)\oA{}_\al{}^{\beta\gamma\!\!\:}z 
\,+\, 
{}^\al v\oA \std{}_{{\!\!\:}^{\!}_\al u\,}z\,,\g\re\\[.4em]
 &= \le \:\! {}^{\al\:\!\!} (\std{}_{{}^\eta{}^{}_{\gamma\!\!\:} v}\,
{}_{\eta\beta}u + [{}^{}_\be u,{}^{}_\gamma v])\oA 
{}_\al{}^{\beta\gamma\!\!\:}z 
\,,\g\re
\,+\,
\le \!\:{}^\al v\oA\std{}_{{\!}^{\!}_\al u\,}z\,,\g\re \\[.4em]
&=\le {}^{\beta\gamma}z\oA \std{}_{{}^\eta{}^{}_{\gamma\!\!\:} v}\,
{}_{\eta\beta}u\,,\g\re\,+\, \le  {}^\al[{}^{}_\be u,{}^{}_\gamma v]\oA {}_\al{}^{\beta\gamma\!\!\:}z \,,\g\re
\,+\,
\le \!\:{}^\al v\oA\std{}_{{\!}^{\!}_\al u\,}z\,,\g\re \\[.4em]
&=\le {}^{\beta\gamma}z\oA \std{}_{{}_\beta{}^{\eta\!\!\:} v}\,
{}_{\gamma\eta}u\,,\g\re\,+\, \le  [ u, v]\oA z \,,\g\re
\,+\,
\le \!\:{}^\al v\oA\std{}_{{\!}^{\!}_\al u\,}z\,,\g\re 
\end{split}
\end{equation}
where in the second line we used Corollary \ref{Leibcovderu} for the
braided derivation rule of the covariant derivative; in the third line we used the torsion free condition $T(u,v)=\std{}_uv-\std{}^{}_{{{\!}^\eta
  v}}{{}_{\,\eta\;\!\!} u}\,-[u,v]=0$; in the fourth
the braided symmetry of the metric
and that the adjoint of the braiding on forms is the braiding on
vector fields, cf. equation \eqref{taudual}; in the last line we used the
Yang--Baxter equation (in the form
$\RR^{-1}_{23}\RR^{-1}_{13}\RR^{-1}_{12}=\RR^{-1}_{12}\RR^{-1}_{13}\RR^{-1}_{23}$)  and that
the braided bracket $[~,~]$ is $H$-equivariant.

The identity in \eqref{Luvzg} is $\Bbbk$-linear in $u,v,z$ and we rewrite it for the cyclically permuted elements
$u\otimes v\otimes z\mapsto {}^{\al\beta} z\otimes {}_\al u\otimes
{}_\beta v$ and $u\otimes v\otimes z\mapsto {}^\eta v\otimes{}^\gamma z\otimes
{}_{\gamma\eta} u$. We then subtract the second identity from the first and
add the third thus obtaining (after using  the Yang--Baxter equation,
the  braided symmetry of the metric and the braided antisymmetry of the
Lie bracket of vector fields)
\begin{equation}\label{Koszul}
\begin{split}
2 \le \!\:{}^\al v\oA\std{}_{{\!}^{\!}_\al u\,}z\,,\g\re =&\;
\LLL_u\le v\oA z,\g\re- 
\LLL_{{}_{}^\al v}\le {{}_\al u}\oA z,\g\re+
\LLL_{{}^{\al\beta\!\!\:} z}\le {{}_\al u}\oA {}^{}_\beta v,\g\re\\[.4em]
&\;-
 \le  [ u, v]\oA z \,,\g\re+
 \le  u\oA [v, z] \,,\g\re+
\le  [ u, {}^\beta z]\oA {}^{}_\beta v \,,\g\re~.
\end{split}
\end{equation}
The right hand side of this identity uniquely determines
the left hand side, that, in turn, because of the exactness of the pairing 
$\le~,~\re:\VV\otimes_A\Omega(A)\to A$ and the invertibility of $\g^\flat$,
uniquely determines the covariant derivative 
$\std{}_u:\Vect(A)\to \Vect(A)$ for all $u\in\Vect(A)$, i.e.,
$\std^\cd:\Vect(A)\otimes\Vect(A)\to \Vect(A)$,
$\std^\cd(u\otimes z)=\std{}_u(z)$. Recalling
Remark \ref{remcd} this proves uniqueness of the metric compatible and
torsion free connection $\std$.
\end{proof}
\begin{remark}\label{gequi}
If $\g$ is $H$-equivariant, i.e., if $h\ra\g=\varepsilon(h)\g$ for all $h\in
H$ then $\std$ is $H$-equivariant, $h\ra\std=\varepsilon(h)\std$.
This can be show by acting with $h$ on \eqref{Koszul}. Due to
$H$-equivariance of $\LLL$, $\le~,~\re$, $\g$ and $[~,~]$ the right hand
side is obtained replacing $u\otimes v\otimes z$ in \eqref{Koszul}
with  $h_{(1)}\ra u \otimes h_{(2)}\ra v\otimes h_{(3)}\ra z$ which
equals $2 \le \!\:{}^\al  (h_{(2)}\ra v)\oA\std{}_{{\!}^{\!}_\al
  (h_{(1)}\ra u)\,}( h_{(3)}\ra z)\,,\g\re$. On the other hand the $h$
action on the left hand side gives $2 \le (h_{(1)}\ra\!\:{}^\al v)\oA
(h_{(4)}\ra^\cop\std)^{}_{h_{(2)}\ra{\;\!}_\al  u\,}( h_{(3)}\ra
z)\,,\g\re$, cf. \eqref{copadjact2}. Comparison of these two
expressions shows that  $h\ra\std=\varepsilon(h)\std$ for any $h\in H$.
\end{remark}

\begin{remark} 
Using the braided symmetry of the metric
and that the adjoint of the braiding on forms is the braiding on
vector fields we can rewrite the left hand side of \eqref{Koszul} as
$$
\le \!\:{}^\al v\oA\std{}_{{\!}^{\!}_\al u\,}z\,,\g\re=\le {}^\be\big(\std{}_{{\!}^{\!}_\al u\,}z\big)\oA\!\:{}_\be{^\al}
v\,,\g\re=
\le ({}^\delta\std){}^{}_{{\!}^\be{\!}_{\al\!\:\!} u\,}{}^{\gamma\!} z\oA\!\:{}_{\delta\gamma\be}{^\al}
v\,,\g\re=
\le ({}^\delta\std){}^{}_{u\,}{}^{\gamma\!} z\oA\!\:{}_{\delta\gamma}v,\g\re
$$
where $ ({}^\delta\std){}^{}_{u\,} z=\le u, (\oR^\delta\ra\std) z\re$
and in the second equality we used that
$h\ra\big(\std{}_{{\!}^{\!}_\al u\,}z\big)=
h\ra\le {}_\al u, \std z\re=\le h_{(1)}\ra {}_\al u,  h_{(2)}\ra(\std z)\re$, for
all $h\in H$,  the braid property
$(\Delta\otimes\id)\mathcal{R}^{-1}=\mathcal{R}_{23}^{-1}\mathcal{R}_{13}^{-1}$
and equation \eqref{copadjact2}.
By $\Bbbk$-linearity in  $u,v,z$ and their arbitrariness, equation
\eqref{Koszul} holds also when rewritten for the element 
$u\otimes {}^\be v\otimes {}_\be z$ instead of $u\otimes v\otimes z$,
it reads
\begin{equation*}
\begin{split}
2\le ({}^\delta\std){}^{}_{u\,} z\oA\!\:{}^{}_{\delta}v,\g\re
=&~
\LLL_u\le z\oA v,\g\re- 
\LLL_{{}_{}^{\al\be} v}\le {{}_\al u}\oA {}_\be z,\g\re+
\LLL_{{}^{\al\!\!\:} z}\le {{}_\al u}\oA  v,\g\re\\[.4em]
&\;
- \le  [ u, {}^\be v]\oA {}_\be z \,,\g\re
- \le  u\oA [z,v] \,,\g\re
+\le  [ u,  z]\oA  v \,,\g\re~,
\end{split}
\end{equation*}
where  we simplified $\LLL_u\le {}^\be v\oA {}_\be z,\g\re=\LLL_u\le
z\oA v,\g\re$ and $[{}^\be v ,{}_\be z]=-[z,v]$.
The right hand side of this expression equals the right hand side of
the Koszul formula in \cite[eq. (6.65)]{Weber} (use $\le  [ u, {}^\be v]\oA {}_\be z \,,\g\re=
\le  {}^\gamma z\oA [{}_\gamma u, v]\,,\g\re=$ $-\le{}^\gamma
  z\oA  [ {}^\beta v, {}_{\be\gamma} u] \,,\g\re$
and rename $u,v,z$ as $X,Z,Y$). We remark that
the left hand side however differs, it equals that  in \cite[eq.
(6.65)]{Weber}, namely $2\le \std{}^{}_{u\,} z\oA v,\g\re$, only when the
action of the braiding on the connection is trivial, this is indeed the
case considered there, where the metric is $H$-equivariant
(hence in particular central, for all $a\in A$, $\g a={}^\al a\!\:{}_\al\g$ simply reads $\g a=a\g$)
and therefore the connection too is $H$-equivariant, cf. Remark \ref{gequi}. Thus the present result,
where we consider an arbitrary metric $\g$, generalizes to not necessarily
$H$-equivariant metrics the Koszul formula  obtained in \cite{Weber}.
\end{remark}

Let  $\kg:\Vect(A)
\otimes\Vect(A)\otimes \Vect(A)\to A$, $u\otimes v\otimes z\mapsto
\kg(u\otimes v\otimes z)$ be the $\Bbbk$-linear map defined by the right hand side
of equation \eqref{Koszul}.
The map $\kg$ is a $\Bbbk$-linear combination of compositions of the map
 $\le~\,,\g\re
:\Vect(A)\oA\Vect(A)\to A$, $v\oA z\mapsto \le~~,\g\re(v\oA
z)=\le v\oA z,\g\re$,
with the maps $\LLL: \Vect(A)\otimes A\to A$, $[~,~]:
\Vect(A)\otimes\Vect(A)\to \Vect(A)$, the braiding $\tau:
\Vect(A)\otimes\Vect(A)\to \Vect(A)\otimes\Vect(A)$ and the projection
$\pi: \Vect(A)\otimes\Vect(A)\to \Vect(A)\oA\Vect(A)$ to the balanced
tensor product over $A$. For example the first addend on the right hand side of
\eqref{Koszul} reads $\LLL_u\le v\oA
z,\g\re=\LLL\circ(\id\otimes\le~~,\g\re\circ\pi)\!\:(u\otimes
v\otimes z)$, where $\id$ stands for $\id_{\Vect(A)}$; the explicit
expression of $\kg$ in terms of these maps is
\begin{equation}\label{kgcomp}
  \kg=\LLL\circ(\id\otimes\le~~,\g\re\circ\pi)\circ(\id-\tau_{12}+\tau_{12}\circ
  \tau_{23})\,-\,\le~~,\g\re\circ\pi\circ
  \big([~,~]\otimes\id-\id\otimes
  [~,~]-([~,~]\otimes \id)\circ \tau_{23}\big)\:.\,
  \end{equation}
Existence of the Levi-Civita
connection is proven by studying the properties of this map.
\\[-.4em]

Recall that ${}_A\hom({\Vect(A) \otimes_A\Vect(A)\otimes  \Vect(A),
  A})\subset {}_\Bbbk\hom({\Vect(A) \otimes_A\Vect(A)\otimes  \Vect(A),
  A})$ is the submodule in $\catm$ of left $A$-linear maps; it is
not a module in $\cat$ because $\Vect(A) \otimes_A\Vect(A)\otimes \Vect(A)$ is not in $\cat$.
\begin{lemma}
The $\Bbbk$-linear map $\kg$ is a left $A$-linear map
in ${}_A\hom({\Vect(A) \otimes_A\Vect(A)\otimes  \Vect(A), A})$
 and satisfies the derivation property, for all $u, v, z\in \Vect(A), a\in
A$, 
\begin{equation}\label{kgder}
\kg(u\otimes_A v\otimes az)=\kg(u\otimes_A va\otimes z)+2\le{}^\al
v_{\!\;}\LLL_{{\!}^{}_\al\!\!\: u}(a)\oA  z,\g\re~.
\end{equation}
\begin{proof}
  We first show that $\kg\in {}_\Bbbk\hom({\Vect(A)
    \otimes_A\Vect(A)\otimes  \Vect(A), A})$, i.e., that $H$ acts on $\kg$ via the $\ra^\cop$ adjoint
action. Recall from \eqref{eqcompadjcop} that the composition of
internal morphisms is an internal morphism.
We prove that  $\kg$ carries the $\ra^\cop$ adjoint
action by showing that its components in \eqref{kgcomp}
carry the $\ra^\cop$ adjoint action.
The map $\le~\,,\g\re
:\Vect(A)\oA\Vect(A)\to A$, $v\oA z\mapsto \le~~,\g\re(v\oA
z)=\le v\oA z,\g\re$, is easily seen to be in
${}_\Bbbk\hom(\Vect(A)\oA\Vect(A), A)$, indeed, for all $h\in H$ and $v,z\in
\Vect(A)$, $$
h\ra(\le~~,\g\re(v\oA
z))=h\ra\le v\oA z,\g\re=\le h_{(1)}\ra (v\oA z),h_{(2)}\ra\g\re=
(h_{(2)}\ra^\cop\le~~,\g\re)(h_{(1)}\ra (v\oA
z))~.$$
The maps $\LLL, [~,~], \tau$ and $\pi$ are all $H$-equivariant and
hence can be seen as internal morphisms with trivial $\ra^\cop$
adjoint action (recall end of Section \ref{subsecMoHA}).
Thus $\kg$ is a composition of maps that carry the $\ra^\cop$ adjoint action and therefore
 $\kg\in {}_\Bbbk\hom({\Vect(A) \otimes\Vect(A)\otimes\Vect(A), A})$.
\\[-.4em]

Next we show that the map $\kg$ is well-defined on the 
balanced tensor product $\Vect(A) \otimes_A\Vect(A)\otimes  \Vect(A)$. 
The third addend in the right hand side of \eqref{Koszul} is well-defined because
$\LLL_{{}^{\al\beta\!} z}\le {{}_\al u}\oA {}^{}_\beta v,\g\re=
\LLL_{{}^{\al\!\!\:} z}\le {}^{}_{\,\al\!\!\;} (u\oA  v),\g\re$.
From the identity 
$\LLL_u(v)=[u,v]=-[{}^\al v,{}_\al u]=-\LLL_{{}^\al v}({}_\al u)$
and the braided Leibniz rule of the Lie derivative we have 
$$u\oA [v,z]+[u,{}^\beta z]\oA {}_\beta v=-\LLL_{{}^{\al\beta\!} z}({{}_\al u}\oA {}^{}_\beta v)=-\LLL_{{}^{\al\!}
  z}{\,}^{}_\al\!\;\!(u\oA v)$$ 
that implies that also the last two addends
of $\kg$ are well-defined on $\Vect(A) \otimes_A\Vect(A)\otimes
\Vect(A)$.
We are left to prove the equality
$\kg(ua\otimes v\otimes
z)=\kg(u\otimes av\otimes z)$ 
for the sum of the first, second and fourth addend in $\kg$. 
This is directly checked recalling that $\LLL_{av}=av=a\LLL_v$ on $A$
and using the compatibility between the Lie bracket and the $A$-module
structure of vector fields
$$[u,av]=\LLL_u(av)=\LLL_u(a) v+{}^\al a_{\,}\LLL_{{\!}_\al u}(v)
=\LLL_u(a) v+{}^\al a[{}_\al u,v]~,$$ 
$$
[au, v]=-[{}^{\al\beta}v,{}_\al a{\:\!}^{}_\beta
u]=-\LLL_{{}^{\al\beta}v}({}_\al a){}_\beta u -a\!\:\LLL_{{}^\beta v\:\!}{}^{}_\beta
  u=a[u,v]- \LLL_{{}^{\al\beta}v} ({}_\al a){}_\beta u~.$$
In order to show that $\kg\in {}_A\hom({\Vect(A)
  \otimes_A\Vect(A)\otimes \Vect(A), A})$  we are left to prove left
$A$-linearity of $\kg$.
This follows from left $A$-linearity of the
second plus fourth addend and of the third plus sixth addend in the
right hand side of \eqref{Koszul}.
\\[-.4em]

The derivation property is equivalent to 
\begin{equation}\label{kgder2}
\kg(u\otimes_A v\otimes az)\,=\,{}^{\delta\eta}
a_{\,}\kg({}^{}_\delta u\otimes_A {}^{}_\eta v\otimes
z) +
2\le{}^\al
v_{\!\;}\LLL_{{\!}^{}_\al\!\!\: u}(a)\oA  z,\g\re~.
\end{equation}
We use the braided Leibniz rule of the Lie derivative 
on covariant and contravariant tensors (cf. \eqref{LonOm}) in the
first two addends of $\kg$ in \eqref{Koszul}
and, since the metric is braided symmetric,  we 
also write the last addend of $\kg$ as 
$\le  [ u, {}^\beta z]\oA {}^{}_\beta v \,,\g\re=\le{}^\al v\otimes_A 
[{}_\al u,z],\g\re=\le{}^\al v\otimes_A 
\LLL_{{}_\al u}z ,\g\re$, thus obtaining the following expression,
\begin{equation}
\begin{split}\label{kg2}
\kg(u\otimes_A v\otimes z)=\:&\le {}^\al v\otimes_A {}^\beta z,\LLL_{{}_{\!\beta\al}
  u}{}\g\re+\le[u,v]\otimes_A z,\g\re
-\le u\otimes_A {}^\al z,\LLL_{{}_\al
  v}\g\re\\[.4em]
&+\LLL_{{}^{\al} z}\le{}_\al (u\otimes_A
v),\g\re+2\le{}^\alpha v\otimes_A \LLL_{{}_\al u}z,\g\re~. 
\end{split} 
\end{equation} 
We use this expression to compute the left hand side and the right
hand side of  \eqref{kgder2};
each of the first four addends in \eqref{kg2}
satisfies the homogeneous version of equation \eqref{kgder2}, for
example we have 
$$
\le {}^\al v\otimes_A {}^\beta (az),\LLL_{{}_{\!\beta\al}
  u}{}\g\re=
\le {}^\al v\otimes_A {}^\eta a {}^\gamma z,\LLL_{{}_{\!\gamma\eta\al}
  u}{}\g\re=
{}^{\delta\eta} a\le {}_\delta{}^\al v\otimes_A  {}^\gamma z,\LLL_{{}_{\!\gamma\eta\al}
  u}{}\g\re=
{}^{\delta\eta} a\le {}^\al{}_\eta v\otimes_A  {}^\gamma z,\LLL_{{}_{\!\gamma\al\delta}
  u}{}\g\re
$$
where in the last equality we used the Yang--Baxter equation  (in the form
$\RR^{-1}_{23}\RR^{-1}_{13}\RR^{-1}_{12}=\RR^{-1}_{12}\RR^{-1}_{13}\RR^{-1}_{23}$).
The last addend in \eqref{kg2} gives also an inhomogeneous term:
$\!\!\:
2\le{}^\alpha v\otimes_{\!\!\:A} \LLL_{{\!}_\al u}(az),\g\re\!\!\;=
2\le{}^\alpha v\otimes_A \LLL_{{}_\al u}(a)\,z,\g\re+
{}^{\delta\eta}a \;\!2\le{}_\delta{}^\alpha v\otimes_A \LLL_{{}_{\eta\al} u}z,\g\re
$,
thus proving  that $\kg$ in \eqref{kg2} satisfies \eqref{kgder2}.
\end{proof}
\end{lemma}

\begin{theorem}[Levi-Civita connection] \label{LCE}   
  Let  $H$ be a triangular Hopf algebra,  $A$ a braided commutative
  $H$-module algebra, let the associated module in $\cat$  of braided derivations $\VV$ be finitely generated and projective, $\OOm:={}^*\VV$ be the right dual module of one-forms and $\g\in\OOm\oA\OOm$ be a metric on $A$.
There exists a unique torsion free metric compatible left connection $\std\in
{}_A\Con(\Vect(A))$.

The explicit expression of the Levi-Civita connection in terms of the inverse metric tensor $\bar \g=\bar\g^b\otimes_A\bar\g_b$ (cf. \eqref{invmetric}), the dual basis $\{e_i,\omega^i : i=1\ldots, n\}$ of $\VV$ and the exact pairings $\le~,~\re: \VV\otimes_A\OOm\to A$ and $\le~,~\re':\OOm\otimes_A \VV\to A$ (cf. \eqref{tevcoev}) is, for all $v,z\in\VV$,
\begin{equation}\label{LEVICIVITA}
\std{}_v(z) =\frac{1}{2} \le \omega^i\kg({}^\al v\otimes_A {}_\al e_i\otimes z)\,,\, \bar\g^b\re'\,\bar\g_b
\end{equation}
where $\kg(u\otimes v\otimes z)\in A$ is given by the right hand side of equation \eqref{Koszul} for all  $u, v,z\in \VV$.
\begin{proof}
Uniqueness has been proven in Theorem \ref{LCU}. We are left to prove existence.
Since $\g$ is nondegenerate we define the map $\std^\cd:\VV\otimes\VV\to\VV)$
implicitly by, for all $u,v,z\in \VV$,  
\begin{equation}\label{LCexists}2\le  u\otimes_A
\std^{\cd}({} v\otimes z),\g\re=\kg({}^\al v\otimes_A {}_\al u\otimes z)~.
\end{equation}
Recalling that
$\le  u\otimes_A\std^{\cd}({} v\otimes z),\g\re=\le  u,\g^\flat(
\std^{\cd}({} v\otimes z))\re$, cf. \eqref{gbemoll},
the explicit expression is $\std^\cd=\frac{1}{2}{(\g_{}^\flat})^{-1}\circ \big(\kg\circ (\tau\otimes
\id_{\VV})\big)^\sharp$, where $\kg\circ (\tau\otimes
\id_{\VV})\in {}_A\hom(\VV\oA\VV \otimes \VV,A)$,
$$\sharp: \!\;\!{}_A\hom(\VV\oA\VV \otimes
\VV,A)
\rightarrow
{}_A\hom(\VV^{\otimes 2},\Omega(A)),~\tilde L\mapsto \tilde L^{\sharp\,};~  \tilde L^\sharp(v\otimes z)=\omega^i\tilde L(e_i\otimes_Av\otimes z)$$
is the isomorphism of
Theorem \ref{Thmfgp123} (equations \eqref{iso1e2}, \eqref{sharpmap}),
$(\g^\flat)^{-1}\in {}_A\hom(\Omega(A),\VV)$, 
and $$(\g_{}^\flat)^{\!-1\,}\circ:
{}_A\hom(\VV^{\otimes 2},\Omega(A))\to 
{}_A\hom(\VV^{\otimes 2},\VV)~, ~~\tilde L\mapsto
(\g_{}^\flat)^{\!-1}\circ\tilde L~.$$
This show that
$\std^\cd
\in {}_A\hom(\VV\otimes\VV,\VV)$.
The map $\std^\cd$ is furthermore a covariant derivative (as defined in Remark \ref{remcd})
because the derivation property \eqref{kgder} of $\kg$
implies the Leibniz rule $\std^\cd( v\otimes az)=\std^\cd(va\otimes
z)+\LLL_v(a)z$,
\begin{equation*}
  \begin{split}
\le  u\otimes_A
\std^{\cd}( v\otimes az),\g\re =\frac{1}{2}\kg({}^\al
v\otimes_A{}_\al u\otimes az)
&=\frac{1}{2}\kg({}^\al v \!\:{}^\be\! a\otimes_A {}_{\be\al} u\otimes z)+\le
u_{\!\;}\LLL_{v}(a)\oA  z,\g\re\\
&=\frac{1}{2}\kg({}^\al (v a)\otimes_A {}_\al u\otimes z)+\le
u_{\!\;}\oA\LLL_{v}(a) z,\g\re\\
&=\le  u\otimes_A
\std^{\cd}( v a\otimes z),\g\re + \le u_{\!\;}\oA\LLL_{v}(a)  z,\g\re~.
\end{split}
\end{equation*}
From  Remark \ref{remcd} it then follows that there exists a unique
connection $\std\in {}_A\Con(\VV)$ with the property 
$\std{}_v(z)=\std^\cd(v\otimes z)$, for all $v,z\in \VV$.

We now show that this connection $\std$ is torsion free.  For all $u, v, z\in \VV$ consider the
permutation $u\otimes v\otimes z\mapsto {}^{\gamma\delta}z
\otimes {}^\eta{}_{\delta} v\otimes {}_{\eta\gamma} u $.
On the one hand, from the Yang--Baxter equation and \eqref{LCexists}, 
\begin{equation*}
  \begin{split}
  \kg(u\otimes_A v\otimes z)-\kg({}^{\gamma\delta}z\otimes
{}^\eta{}_{\delta} v\otimes {}_{\eta\gamma} u)&=
\kg(u\otimes_A v\otimes z)-\kg({}^{\gamma\delta}z\otimes
{}_\gamma{}^{\eta} v\otimes {}_{\delta\eta} u)\\[.2em] &=
2\le {}^\eta v\otimes_A (\std{}_{{}_{\!\eta\!\!\:} u
}z-\std{}_{{\!\!\:}^\delta \!\!\:z\,}{}_{\delta\eta} u) ,\g\re~.
\end{split}
\end{equation*}
On the other hand, recalling the definition of $\kg$, and using the braided symmetry of the metric, the left hand side of the above expression equals $2\le [u,{}^\beta z]\otimes_A {}^{}_\beta v,\g\re=2\le {}^\eta v\otimes_A 
[{}_\eta u, z],\g\re$, thus showing that $\std$ is torsionless.

We similarly show the metric compatibility of $\std$.
From  \eqref{LCexists} and Corollary \ref{Leibcovderu} we have
$$\kg(u\otimes_A v\otimes z)+\kg(u\otimes_A {}^\gamma z\otimes {}_\gamma v)=2\le \std{}_u(v\otimes z),\g\re=2\LLL_u\le v\oA z,\g\re-2\le 
u\otimes_A v\otimes z , {}^*\std\g\re~.$$ From the definition of
$\kg$, cf.  \eqref{Koszul}, it easily follows that the left hand side of this expression simplifies to  $2\LLL_u\le v\otimes_A z,\g\re$, thus proving metric compatibility of the connection $\std$.

The explicit expression \eqref{LEVICIVITA} of the Levi-Civita connection follows from 
$\std{}_v(z)=\std^\cd(v\otimes z)= \frac{1}{2}{(\g_{}^\flat})^{-1}\big((\kg\circ (\tau\otimes
\id_{\VV}))^\sharp(v\otimes z)\big)$ recalling that $(\bar\g)^{\flat'}=(\g^\flat)^{-1}$. 
\end{proof}
\end{theorem}

\subsection{Ricci tensor, scalar curvature and Einstein manifolds}
There is a canonical notion of trace in a ribbon
category and a fortiori in a compact closed category. For an internal morphism $\tilde L\in {}_A\hom(\VV,\VV)$ in
$\catfp$ we have  {tr}$(\tilde L)= \le~,~\re'\circ
(\id_{\Omega(A)}\otimes_A \tilde L)\circ \coev$ that belongs to
${}_A\hom(A,A)$ in $\catfp$ because the coevaluation map $\coev: A\to \OOm\otimes_\AA \Vect(A)$ and the exact
pairing $\le~,~\re':\Omega(A)\otimes_A\VV\to A$
are morphisms in $\catfp$.
The trace is determined by its value on $1_A\in A$, which is $\,${tr}$(\tilde L)=\le \omega^i, \tilde
L( e_i)\re'$, where we used a dual
basis of $\Vect{(A)}$. The { Ricci tensor} is the trace of the
Riemann tensor given by (use  $A\oA\VV\oA\VV\simeq\VV\oA\VV$),
\begin{equation}\label{defofric1}
\begin{split}&{Ric}:=\le~\,,\,~\re'\circ (\id_{\Omega(A)}\otimes_A
R{\stdd})\circ(\coev\otimes_A\id_{\VV\otimes_A\VV}):\VV\otimes_A\VV\to A~,~~\\[.4em]
 &{Ric}(u,v)=\le \omega^i, R{\stdd} (e_i,u,v)\re'~.
\end{split}
\end{equation}

Since the curvature $R{\stdd}$ is left $A$-linear it follows  that 
${Ric}\in{}_A\hom(\VV\otimes_A\VV, A)$.
Similarly, the { scalar curvature tensor} ${S}$ is given by 
${S}={Ric}(\bar\g)$, where the inverse metric
$\bar\g\in \VV\otimes_A\VV$ has been defined in \eqref{invmetric}.

We also define an Einstein metric on $A$ to be a metric $\g$
proportional to its Ricci tensor,
\begin{equation}\label{defofeinst}
{Ric}=\lambda \le\,~\,
\,,\g\re~,~~(\lambda\in \Bbbk)~.
\end{equation}
This equation in ${}_A\hom(\VV\oA\VV,A)$ allows to study noncommutative Einstein manifolds.

\section{Examples\label{expex}}
Consider any of the examples discussed in Section \ref{exgeneral}
(cotriangular Hopf algebras, noncommutative manifolds via Drinfeld twists) 
together with a metric $\g$. Then there is a unique Levi-Civita
connection associated with $\g$, with Ricci and scalar curvature as studied above. 
In particular we can consider Einstein in vacuum equations on 
a cotriangular Hopf algebra $A$ for any  metric $\g$ on $A$.

In the following we first briefly show how to recover previous results
on noncommutative tori and (hyper)planes. Then we study the Riemannian
geometry of the Hopf algebra $A=K\otimes K$ with $K$ the four
dimensional Sweedler Hopf algebra with its basic cotriangular structure. We solve the noncommutative
Einstein manifold equations  with a metric that is non-central
and non-equivariant, its scalar curvature is $S=12$.
\subsection{Noncommutative torus}
Let $\Bbbk=\mathbb{C}$ and $\Theta=(\theta_{jk})_{j,k=1,\ldots n}$ be a real valued skew
symmetric matrix and $A$ the noncommutative $n$-dimensional
torus generated by the unitaries $U_j$, $j=1,\ldots n$ with commutation relations
$U_jU_k=e^{2\pi \theta_{jk}}U_kU_j$.
Consider the $n$ derivations $\partial_j: A\to A$,
defined by $\partial_j(U_k)=2\pi\ii \:\!\delta_{j k}U_k$ and $\partial_j(aa')=\partial_j(a)a'+a\partial_j(a')$, for any $a, a' \in A$. 
They generate the commutative and cocommutative Hopf algebra
$H$, the universal enveloping algebra of the Lie
algebra of the (commutative) $n$-torus. 
Considering the triangular structure $R=e^{2\pi i\theta_{jk}\partial_j\otimes
  \partial_k}$ of $H$ (belonging to a suitable topological
completion of $H\otimes H$) we have that $A$ becomes a
braided commutative $H$-comodule algebra.
The derivations $\partial_j$ generate the free $A$-module of
vector fields $\Vect(A)$. Since the action of $H$ on
$\partial_i$ is trivial, Koszul formula \eqref{Koszul} for the Levi-Civita connection
 simplifies to 
$
2 \le \partial_j\oA\std{}_{\partial_i\,}\partial_k\,,\g\re =
\partial_i\le \partial_j\oA \partial_k,\g\re- 
{\partial_j}\le  \partial_i\oA \partial_k,\g\re+
{ \partial_k}\le {\partial_i}\oA \partial_j,\g\re
$ 
and equals that in \cite{Rosenberg}. There, however, the connection
${\std}_u$ and the curvature tensor $R(u,v, z)$ were defined only for $u$ and
$v$ derivations of the algebra $A$ (such as elements of the
$\mathbb{C}$-linear span of the $\partial_j$'s), not for arbitrary
vector fields. 

One can proceed similarly with the quantum (hyper)plane  with
generators $x_j$, $j=1,\ldots n$ and commutation relations
$x_jx_k-x_kx_j= 2\pi\ii\:\!\theta_{jk}$. This is also an $H$-comodule algebra with $\partial_j(x_k)=\delta_{jk}$.
The expression of the Levi-Civita connection is as above and equals
that in \cite{NCG1}.
\subsection{Sweedler Hopf algebra and its tensor product}
We first study a differential calculus on the Sweedler four
dimensional Hopf algebra $K$ as an explicit case
of Example \ref{excotriang}, it is uniquely defined by the choice of the cotriangular structure. There is a whole family of cotriangular structures on $K$ (cf. e.g. \cite[Ex. 2 \S 10.1.1]{KS}), among the different ones we choose the basic one, defined in equation \eqref{ctrK} below.
There is no metric tensor on this differential geometry with
one dimensional tangent and cotangent space since the tensor product
of forms is automatically antisymmetric.
We then exemplify the previous section by studying the Riemannian
geometry of the tensor product $K\otimes K$. Differently from the
noncommutative torus and (hyper)plane cases, the free module of vector fields is
not generated by derivations of the algebra
$K\otimes K$,
we then use the general Koszul formula \eqref{Koszul} 
to determine the Levi-Civita connection.

\subsubsection{Differential calculus on Sweedler Hopf algebra $K$}\label{H4DC}
Let $K$ be the Sweedler Hopf algebra over a field $\Bbbk$ of
characteristic $0$. It is the
algebra generated by $x$ and $g$ and defining relations $g^2=1$,
$x^2=0$, $xg=-gx$.
A vector space basis is given by $(1,g,x,gx)$. The coalgebra structure is determined by
$\Delta(g)=g\otimes g$, $\Delta x=x\otimes g+1\otimes x$, with
antipode $S(g)=g, S(x)=gx$. 
Sweedler Hopf algebra is self dual. Let $U$ be a second copy of $K$,
we rename the generators $g,x$ as $\gamma, \psi$ (with the same
costructures). The relations $<\gamma, 1>=1, <\gamma,g>=-1, <1,g>=1, <\psi,x>=1$, all other
relations among the generators vanishing, extend to a Hopf algebra
pairing $<~,~>: U\otimes K\to \Bbbk$ proving self duality of $K$.

We consider the following triangular structure on $U$,
\begin{equation}\label{ctrK}\mathscr{R}=\mathscr{R}^\al\otimes{\mathscr{R}}_\al=\frac{1}{2}(1\otimes
  1+\gamma\otimes 1+1\otimes \gamma-\gamma\otimes \gamma)~
  \end{equation}
so that $K$ is cotriangular (cf. Example \ref{excotriang}).
The $(U,\mathscr{R})$-Lie algebra
${}_{\rm{inv}}\Vect(K)\simeq\gg$ of left invariant vector
fields is
 $ \gg=\{\chi\in U; \Delta\chi=\chi\otimes
  1+\overline{\mathscr{R}}^\al\otimes\overline{\mathscr{R}}_\al\ra\chi\}
$, cf. \eqref{gforbic}.
It is the linear span of
$$\chi:=\gamma\psi~.$$ Indeed we have $\Delta\chi=\chi\otimes 1+\gamma\otimes\chi
=\chi\otimes
1+\overline{\mathscr{R}}^\al\otimes\overline{\mathscr{R}}_\al\ra\chi$,
while $\overline{\mathscr{R}}^\al\otimes\overline{\mathscr{R}}_\al\ra
\gamma =1\otimes \gamma$ and $\overline{\mathscr{R}}^\al\otimes\overline{\mathscr{R}}_\al\ra
\psi=\gamma\otimes \psi$ imply that the subspace $\gg\subset U$ is one dimensional.

The bicovariant differential calculus associated with the cotriangular
structure of $K$ is therefore determined by the left invariant vector
field $\vvv\in {}_{\rm{inv}}\Vect(K)$, $\vvv(k):=k_1<\chi,k_2>$ for all
$k\in K$. Explicitly, from $<\chi,1>=0, <\chi,g>=0, <\chi,x>=1,
<\chi,gx>=-1 $ we obtain $\vvv(1)=0, \vvv(g)=0, \vvv(x)=1,
\vvv(xg)=g$. Accordingly with the non-cocommutative coproduct of $\chi$,
the left-invariant vector field $u$  is a braided derivation, e.g., $u(gx)=u(g)x-gu(x)$,
$u(xg\,xg)=u(xg)xg+xg\;\! u(xg)$. The braided commutator vanishes, $[\vvv,\vvv]=0$, since
$\overline{\mathscr{R}}^\al\ra\chi\otimes\overline{\mathscr{R}}_\al\ra\chi=-\chi\otimes\chi$
and $\chi\chi=0$. Let
$\omega\in {}_{\rm{inv}}\Omega^1(K)$ be the dual left invariant one-form:
\begin{equation}\label{dualityXomega}
\langle \vvv,\omega\rangle=1~.
\end{equation}  
The $K$-bimodule of one-forms is the free right module
$\Omega^1(K)\simeq {}_{\rm{inv}}\Omega^1(K)\otimes K$. From
$\dd k=\omega \vvv(k)$ for all $k\in K$ (cf. \eqref{defofbicdd}) we have
$$
\dd  g=0~,~~ \dd  x=\omega, ~~\dd (gx)=g\omega  ~.
$$
The bimodule commutation relations $\omega g=-g\omega$, $\omega x=-x\omega$
then read  $(\dd  x) g=-g \dd  x$, $(\dd  x) x=-x\dd  x$ and are
easily obtained from the exterior derivative on $x^2=0, gx=-xg$.
The space of right invariant one forms is the $\Bbbk$-linear span of
$\eta=\omega g=\dd x g$. Dually, from \eqref{dualityXomega} the free bimodule of vector fields is
determined by $\vvv\:\! g=-g\:\! \vvv$, $\vvv\:\! x=-x \:\!\vvv$, it
is the liner span of $u, gu, xu$; only $xu$ is a
derivation, it  does not generate the $K$-bimodule of vector fields.   
\\

Recall from Example \ref{cotriangularK} and \ref{exgeneral} that $K$ is a $U^\op\otimes
U$-module algebra with $R$-matrix
${\RR}=R^\alpha\otimes R_\alpha=(\id\otimes
\rm{flip}\otimes\id)(\mathscr{R}^{-1}\otimes\mathscr{R})$.
The braiding $\tau$  in 
$\Omega(K)\otimes \Omega(K)$ (a representation of the permutation group since
it squares to the identity) is with respect to this triangular
structure. 
We have $\tau(\omega\otimes_K \eta)=\eta\otimes_K\omega$ since
the left-invariant one-form $\omega$ trivialises the $R$-matrix of
$U^\op$ while the right-invariant one $\eta$ that of
$U$ (we thus recover Woronowicz braiding). By $K$-linearity it then follows that
$\tau(\omega\otimes_K\omega)=\tau(\omega\otimes_K\eta) g=\eta\otimes_K\omega g=-\omega\otimes_K\omega$.
The tensor product is therefore antisymmetric. This shows the first of the
following $K$-bimodule isomorphims
\begin{equation}\label{nog}
  \Omega(K)\wedge\Omega(K)\simeq \Omega(K)\otimes_K\Omega(K)\simeq K~.
\end{equation}
The last follows from the $\mathbb{Z}_2$-graded structure of the
bimodule relations $\omega g=-g\omega$, $\omega x=-x\omega$.

More generally, the contravariant tensor algebra is isomorphic to the exterior algebra,
$\Tau^{\bullet,0}=\bigoplus_{r\in\mathbb{N}} \Tau^{r,0}\simeq\OM$,
with bimodule isomorphisms $\Omega^{ 2r}(K)\simeq K$,
$\Omega^{2r+1}(K)\simeq \Omega(K)$.

Form \eqref{nog} we see that there is no symmetric two tensor and
hence no metric tensor on this
differential geometry with
one dimensional tangent and cotangent space and no top form.

\subsubsection{Riemannian geometry on $K\otimes K$}
Define the tensor product algebra $A:=K\otimes K$, where the two copies of $K$
mutually commute.  Let $g,x$ denote the generators of  $K\simeq K\otimes 1\subset
A$, and $g',x'$ those of  $K\simeq 1\otimes K\subset
A$, similarly for  $\gamma, \psi\in U\otimes 1$ and $\gamma', \psi'\in 1\otimes U$. The
triangular structure of ${\mathcal{U}}:=U\otimes U$ is canonically
given as
${\mathfrak{R}}={\mathfrak{R}}^\alpha\otimes {\mathfrak{R}}_\alpha
=\mathscr{R}^\alpha\mathscr{R}'^\beta\otimes
\mathscr{R}_\alpha\mathscr{R}'_\beta$,
whence the cotriangular one of $A$.  Then $A$ is an example of braided commutative  $H$-module algebra,
where  $H={\mathcal{U}}^\op\otimes {\mathcal{U}}$.
The $H$-Lie algebra is the direct sum of
one-dimensional Lie subalgebras: $\gg^{}_{A}=\gg{\;\!}\oplus{\;\!}\gg'\subset U\otimes U$, with
$\chi:=\gamma\psi\in\gg\subset U\otimes 1$,
$\chi':=\gamma'\psi'\in\gg'\subset 1\otimes U$ linear generators of
$\gg,\gg'$. Since $\overline{\mathfrak{R}}_\alpha\otimes
{\overline{\mathfrak{R}}}^\alpha\ra\chi=\overline{\mathscr{R}}_\alpha\otimes
  {\overline{\mathscr{R}}}^\alpha\ra\chi$, and similarly for
    $\chi'$, we have
    $[\chi,\chi']=\chi\chi'-\overline{\mathfrak{R}}_\alpha\ra\chi'\:
\overline{\mathfrak{R}}^\alpha\ra\chi=\chi\chi'-\chi'\chi=0$.
The associated left invariant vector fields are $\vvv, \vvv'$ and
hence have vanishing braided brackets. 
The dual left-invariant one-forms $\omega, \omega'$ commute with $1\otimes K$ and $K\otimes 1$ respectively and generate
the bimodule of one-forms 
$\Omega(A)=\Omega(K)\otimes K\,\oplus\, K\otimes
\Omega(K)$ as a free right (left) $A$-module. 
The exterior derivative is the sum $\dd\otimes 1 +
1\otimes \dd$. 
The bimodule of braided symmetric contravariant 2-tensors is freely generated (as a
right module) by $\omega\vee\omega':=\omega\otimes_A \omega'+\omega'\otimes_A\omega$.
The most general metric on $A$ is $\g=\omega\vee\omega' t$ with $t\in
A$ invertible, this is generally neither central nor $H$-equivariant.
\\

We are now set to study the Levi-Civita connection on the Hopf algebra 
$A$.
Since $A$  is generated by $x,g,x',g'$, we write $a=a( x,g,x',g')$ for
any $a\in A$ and define the algebra involutions
$$\hat a:=a( -x,- g,x',g')~,~~ \hat a':=a( x,g,-x',-g')~,~~\tilde a:=a(-x-g,-x',-g')~$$
so that $\hat a \omega =\omega a$, $\hat a'  \omega' =\omega' a$, $\hat a
\vvv=\vvv a$, $\hat a'  \vvv' =\vvv' a$.
We apply
Theorem \ref{LCE}  and obtain the following explicit expression of the
Levi-Civita connection.

\begin{theorem}
  The Levi-Civita connection  $\std\in  {}_A\Con(\Vect(A))$ of the
  arbitrary  metric tensor $\g=\omega\vee\omega' t$ ($t\in A$ invertible) on the Hopf algebra $A$
with bicovariant differential calculus defined above is given by 
\begin{equation}\label{exexV}
    \begin{split}
      \std: \Vect(A)&\longrightarrow \Omega(A)\otimes_A\Vect(A)\\
      \vvv&\longmapsto  -\omega'\otimes_A\big(\vvv \,\widetilde{\vvv'(t)t^{-1}\,}\!+
      \vvv'\,\widetilde{\vvv(t)t^{-1}}\big)\\
            \vvv'&\longmapsto  -\omega \otimes_A\big(\vvv \,\widetilde{\vvv'(t)t^{-1}}+ \vvv'\,\widetilde{\vvv(t)t^{-1}}\big)~.
    \end{split}\end{equation}
\end{theorem}
\begin{proof}
The Koszul formula
\eqref{Koszul} gives 8 equations, for $u,v,z$ respectively 1) $\vvv, \vvv,
\vvv$, 2) $\vvv, \vvv', \vvv$, 3) $\vvv, \vvv, \vvv'$, 4) $\vvv, \vvv',
\vvv'$ and the corresponding ones with $\vvv\leftrightarrow \vvv'$.
The right hand side of the first equation vanishes since
$\langle\vvv,\omega'\rangle=0$. The left hand side reads
$-2\langle\vvv\otimes_A\std{}_\vvv\vvv,\g\rangle=
-2\langle\vvv,\langle\std{}_\vvv\vvv,\omega'\rangle\omega\rangle
t=2\widehat{\langle\std{}_\vvv\vvv,\omega'\rangle}t$, hence
$\langle\std{}_\vvv\vvv,\omega'\rangle=0$ that is $\std{}_\vvv\vvv=a\vvv$
for some $a\in A$. Similarly, the second equation implies $\std{}_\vvv\vvv=b\vvv'$
for some $b\in A$. Together the first two equations therefore imply $$\std{}_\vvv\vvv=0~.$$
Similarly, the third equation reads
$-2\widehat{\langle\std{}_\vvv\vvv',\omega'\rangle}t=2\vvv(t)$ that implies
$\std{}_\vvv\vvv'=-\widehat{\vvv(t)t^{-1}}\vvv'+a\vvv$ for some $a\in A$.
Together with the fourth, which implies
$\std{}_\vvv\vvv'=-\widehat{\vvv'(t)t^{-1}\;}{\phantom{\!\!\!\!\!\!\!\!\!\!\!\!L^H}}^{\mbox{$'$}}\;\!\vvv+b\vvv$ for some $b\in
A$, we then have
$\std{}_\vvv\vvv'=-\widehat{\vvv'(t)t^{-1}\;}{\phantom{\!\!\!\!\!\!\!\!\!\!\!\!L^H}}^{\mbox{$'$}}\;\!\vvv-\widehat{\vvv(t)t^{-1}}\vvv'$,
i.e., $\std{}_\vvv\vvv'=-\vvv\;\!\widetilde{\vvv'(t)t^{-1}}-\vvv'\;\!\widetilde{\vvv(t)t^{-1}}$.
The remaining four equations, corresponding to $\vvv\leftrightarrow \vvv'$, give
$\std{}_{\vvv'}\vvv'=0$ and $\std{}_{\vvv'}\vvv=\std{}_\vvv\vvv'$. This
proves expression \eqref{exexV}.
  \end{proof}
\begin{remark}  The dual right connection on one-forms (readily obtained form \eqref{dddual}) reads:
$
      {^*\!\!\;}\std: \Omega(A)\longrightarrow \Omega(A)\otimes_A\Omega(A)
,~^*\std(\omega)=\omega\vee\omega' \,{\,\widehat{\vvv'(t)t^{-1}\;}}{\phantom{\!\!\!\!\!\!\!\!\!\!\!\!L^H}}^{\mbox{$'$}}
,~
^*\std(\omega')=\omega\vee\omega'
\widehat{\,\vvv(t)t^{-1}\,}.
$ The torsionless and the metric
compatibility conditions can then be easily checked.
\end{remark}

We next compute the curvature, Ricci tensor and scalar curvature. We set $b:=-\widehat{\vvv'(t)t^{-1}\;}{\phantom{\!\!\!\!\!\!\!\!\!\!\!\!L^H}}^{\mbox{$'$}}$,  $b':=-\widehat{\vvv(t)t^{-1}}$, so that
  $\std{}_\vvv\vvv'=\std{}_{\vvv'}\vvv=b\vvv-b'\vvv'$. Then
  $\std{}_\vvv(\std{}_\vvv\vvv')=u(b)b+u(b')u'+\widehat{b'}(bu+b'u')$, and
  similarly
  $\std{}_{\vvv'}(\std{}_{\vvv'}\vvv)=u'(b')b'+u'(b)u+\widehat{b}{\phantom{\!\!\!\!\!\!\!\!L^H}}^{\mbox{$'\!\!$}}(bu+b'u')$. Since
  $[u,u]=[u',u']=[u,u']=0$ the expression of the
  Riemann tensor simplifies and the only nonvanishing components are:
  \begin{equation}
    \begin{split}&\mbox{$\frac{1}{2}$}R(u,u,u')=R(u,u',u)=-R(u',u,u)=u(b)u+u(b')b'+\widehat{b'}(bu+b'u')\\ &\mbox{$\frac{1}{2}$}R(u',u',u)=R(u',u,u')=-R(u,u',u')=u'(b')u'+u'(b)b+\widehat{b}{\phantom{\!\!\!\!\!\!\!\!L^m}}^{\mbox{$'\!$}}(bu+b'u')~.\end{split}
 \end{equation}
From \eqref{defofric1} (and the pairing $\le \omega,u\re'=-\le
u,\omega \re=-1$) the Ricci tensor is given by 
\begin{equation}\label{ricciexp}Ric(u,u)=\widehat{u(b')}'+\widetilde{b'}\widehat{b'}'~,~~Ric(u,u')=-2\widehat{u(b)}-2b'\!\;\widehat{b}-\widehat{u'(b')}'-{b}\widehat{b'}'
\end{equation}
with $Ric(u',u')$ and $Ric(u',u)$ obtained from the above by exchanging $u\leftrightarrow u',
b\leftrightarrow b', {}\widehat{\phantom{{}_{nn}}}\leftrightarrow{}{\widehat{\phantom{{}_{nn}\!}}'}$.

 The inverse metric is $\g^{-1}=-t^{-1}(u\otimes_Au'+u'\otimes_A u)$ (use
$\g^\flat(u)={\widehat{t}}{\phantom{_{\!}\!\!\!m}}^{\mbox{$'\!\!$}}
\omega'$, $\g^\flat(u')=\widehat{t}\omega$ and
$\le~~,{\g^{-1}}^a\re'\otimes_A \g^{-1}_a={\g^\flat}^{-1}$ where
$\g^{-1}={\g^{-1}}^a\otimes_A \g^{-1}_a)$.
Then by left $A$-linearity of the Ricci tensor the scalar curvature reads
$S=Ric(\g^{-1})=3t^{-1}(\widehat{u(b)}+\widehat{u'(b')}'+b'\;\!\widehat{b}+b\widehat{b'}')$.
\\

The Einstein metric conditions \eqref{defofeinst} for
$\g=\omega\vee\omega' t$ are immediate from \eqref{ricciexp}. A
solution with proportionality factor $\lambda=-6$ is given
by setting $t=1+x+x'$ (so that
$b=-1+x-x'+2xx'$, $u(b)=1+2x'$, $u'(b)=-1+2x$ and $b', u'(b'), u(b')$ 
are obtained by exchanging $x\leftrightarrow x'$).

The Einstein metric
$\g=\omega\vee\omega' (1+x+x')$ is neither central nor
equivariant. Its Levi-Civita connection is not a bimodule connection
(there is no generalised braiding $\sigma$, cf. e.g. \cite[Def. 3.66]{BM}, that satisfies $\sigma(u\otimes
dg)=\std (u\,g)-\std(u)g$, indeed  $\dd g=0$
while $\std (u\,g)-\std(u)g\not=0$).
The scalar curvature is $S=12$.
\\

  \noindent{\large{\bf{Acknowledgments}}}\\
The author would like to thank the organizers of the conference 
``Noncommutative manifolds and their symmetries'' dedicated to Giovanni
Landi on the occasion of his 60th birthday (Scalea 16-20.9.2019) and
the guest of honour for the opportunity to present the initial key findings of this research.

The author acknowledges partial support from INFN, CSN4, Iniziativa
Specifica GSS, from INdaM-GNFM and from Universit\`a del Piemonte Orientale.  
\\[.5em]
\indent The author states that there is no conflict of interest.
Data sharing not applicable to this article as no datasets were generated or analysed during the current study.

\end{document}